\documentclass[oper,nonblindrev]{informs4} % current default for manuscript submission
\usepackage{natbib}
\usepackage{multirow}
\usepackage{booktabs}
 \bibpunct[, ]{(}{)}{,}{a}{}{,}%
 \def\bibfont{\small}%
\TheoremsNumberedThrough     % Preferred (Theorem 1, Lemma 1, Theorem 2)
\EquationsNumberedThrough    % Default: (1), (2), ...
\ECRepeatTheorems
\usepackage{amsmath}
\usepackage{graphicx}
\usepackage{bbm}
\usepackage{amssymb} 
\usepackage{amsfonts}
\usepackage{datetime}
\usepackage{mathtools}
\usepackage{comment}
\usepackage{xcolor}
\usepackage{threeparttable}
\usepackage{algorithm}
\usepackage{algorithmic}
\usepackage{caption}
\usepackage{subcaption}
\usepackage{endnotes}
\usepackage[colorinlistoftodos]{todonotes}
\definecolor{links}{RGB}{204,36,29}
\usepackage[colorlinks=true,breaklinks=true,bookmarks=true,urlcolor=links,citecolor=links,linkcolor=links,bookmarksopen=false,draft=false]{hyperref}

\newcommand{\sumj}{\sum_{ j\in[N]}} % product sum
\newcommand{\sumi}{\sum_{i\in [N]}} % product sum
\newcommand{\sumt}{\sum_{t \in [T]}} % period sum
\newcommand{\summ}{\sum_{m\in [M]}} % memory sum

\def \x {\boldsymbol{x}}

\def \boldb {\boldsymbol{\beta}}

\def \h {\boldsymbol{h}}
\def \w {\boldsymbol{w}}
\def \v {\boldsymbol{v}}

\def \bbmr {\mathbbm{R}}
\def \r {\boldsymbol{r}}
\def \y {\boldsymbol{y}}
\def \a {\boldsymbol{a}}
\def \u {\boldsymbol{u}}
\def \e {\boldsymbol{e}}

\makeatletter
\newcommand*{\rom}[1]{\expandafter\@slowromancap\romannumeral #1@}
\makeatother

\def \z {\boldsymbol{z}}
\def \w {\boldsymbol{w}}
\def \ind {\mathbbm{1}}

\def\argmax{\mathop{\mathrm{argmax}}}

\DeclareMathOperator{\rev}{Rev}

\DeclareMathOperator{\conv}{conv}
\DeclareMathOperator{\conc}{conc}

\DeclareMathOperator{\pers}{pers}

\DeclareMathOperator{\Z}{\mathbb{Z}}

\DeclareMathOperator{\for}{ for }
\DeclareMathOperator{\st}{s.t.}
\DeclareMathOperator{\MNL}{MNL}

\def\mcirc{\mathop{\circ}}

 \renewcommand{\theARTICLETOP}{}
 
\newcommand{\huan}{\textcolor{black}}  
\newcommand{\re}{\textcolor{black}} 

\begin{document}
%%%%%%%%%%%%%%%%

% Outcomment only when entries are known. Otherwise leave as is and 
%   default values will be used.
%\setcounter{page}{1}
%\VOLUME{00}%
%\NO{0}%
%\MONTH{Xxxxx}% (month or a similar seasonal id)
%\YEAR{0000}% e.g., 2005
%\FIRSTPAGE{000}%
%\LASTPAGE{000}%
%\SHORTYEAR{00}% shortened year (two-digit)
%\ISSUE{0000} %
%\LONGFIRSTPAGE{0001} %
%\DOI{10.1287/xxxx.0000.0000}%

% Author's names for the running heads
% Sample depending on the number of authors;
% \RUNAUTHOR{Jones}
% \RUNAUTHOR{Jones and Wilson}
% \RUNAUTHOR{Jones, Miller, and Wilson}
% \RUNAUTHOR{Jones et al.} % for four or more authors
% Enter authors following the given pattern:
%\RUNAUTHOR{}

% Title or shortened title suitable for running heads. Sample:
% \RUNTITLE{Bundling Information Goods of Decreasing Value}
% Enter the (shortened) title:
% \RUNTITLE{}

% Full title. Sample:
% \TITLE{Bundling Information Goods of Decreasing Value}
% Enter the full title:
\TITLE{Assortment Optimization Under History-Dependent Effects}

% Block of authors and their affiliations starts here:
% NOTE: Authors with same affiliation, if the order of authors allows, 
%   should be entered in ONE field, separated by a comma. 
%   \EMAIL field can be repeated if more than one author
\ARTICLEAUTHORS{%
\AUTHOR{Taotao He}
\AFF{Antai College of Economics and Management, Shanghai Jiao Tong University, \href{mailto:hetaotao@sjtu.edu.cn}{hetaotao@sjtu.edu.cn}, \href{https://taotaoohe.github.io/}{https://taotaoohe.github.io/}}
\AUTHOR{Yating Zhang}
\AFF{Antai College of Economics and Management, Shanghai Jiao Tong University, \href{mailto:ytzhang20@sjtu.edu.cn}{ytzhang20@sjtu.edu.cn}}
\AUTHOR{Huan Zheng}
\AFF{Antai College of Economics and Management, Shanghai Jiao Tong University, \href{mailto:zhenghuan@sjtu.edu.cn}{zhenghuan@sjtu.edu.cn}}
} % end of the block

\ABSTRACT{% 
This paper examines how to plan multi-period assortments when customer utility depends on historical assortments. We formulate this problem as a nonlinear integer programming model and show it is NP-hard in the presence of a negative history-dependent effect (such as a satiation effect). We build solution methodologies for obtaining global optimal solutions \re{under a general setting that the history-dependent effects could be a mixture of positive and negative}. We propose using a lifting-based framework to reformulate the problem as a mixed-integer exponential cone program that state-of-the-art solvers can solve. We also design a sequential revenue-ordered policy and show that it solves our problem to optimality in polynomial time when historical assortments positively affect customer utility (such as an addiction effect). Additionally, we identify an optimal cyclic policy for an asymptotic regime, and we also relate its length to the customer's memory length. Finally, we present a case study using a catering service dataset, showing that our model demonstrates good fitness and can effectively balance variety and revenue.}

%  We formulate this problem as a nonlinear integer programming model and propose solution methodologies for different regimes of this problem. First, we present a sequential revenue-ordered policy and show that it is optimal when historical assortments positively affect customer utility (such as an addiction effect). Second, we show that the problem is NP-hard in the presence of a negative history-dependent effect (such as a satiation effect). \re{To tackle the general case where the history-dependent effects could be a mixture of addiction and satiation,} we propose using a lifting-based framework to reformulate the problem as a mixed-integer exponential cone program that state-of-the-art solvers can solve. Additionally, we identify an optimal cyclic policy for an asymptotic regime, and we also relate its length to the customer's memory length. Finally, we present a case study using a catering service dataset, showing that our model demonstrates good fitness and can effectively balance variety and revenue.

% Sample
%\KEYWORDS{deterministic inventory theory; infinite linear programming duality; 
%  existence of optimal policies; semi-Markov decision process; cyclic schedule}

% Fill in data. If unknown, outcomment the field
\HISTORY{This version: \today}

\maketitle
%%%%%%%%%%%%%%%%%%%%%%%%%%%%%%%%%%%%%%%%%%%%%%%%%%%%%%%%%%%%%%%%%%%%%%

% Samples of sectioning (and labeling) in MSOM
% NOTE: (1) \section and \subsection do NOT end with a period
%       (2) \subsubsection and lower need end punctuation
%       (3) capitalization is as shown (title style).
%
%\section{Introduction.}\label{intro} %%1.
%\subsection{Duality and the Classical EOQ Problem.}\label{class-EOQ} %% 1.1.
%\subsection{Outline.}\label{outline1} %% 1.2.
%\subsubsection{Cyclic Schedules for the General Deterministic SMDP.}
%  \label{cyclic-schedules} %% 1.2.1
%\section{Problem Description.}\label{problemdescription} %% 2.

\section{Introduction}
% structure
% \begin{itemize}
%     \item potential revenue opportunity by incorporating history dependent effects into assortment
% \end{itemize}

Many companies sell products to customers who make repeat purchases or visit over multiple periods. For example, corporate cafeterias provide workday luncheons for employees; \re{online grocery flash sell platform plans assortments for customers who repeatedly visit the online store and seek deep discounts}. Customers' utility may vary over time and may divert to products offered by competitors if they feel ``bored''. Therefore, these companies must carefully plan their daily assortments over time to keep their customers satisfied. 

\huan{A typical example is the corporate dining industry. Corporate dining is an essential but often overlooked industry. It was reported that a working person will spend an average of \$2500 per year on working-day lunches~\citep{ezcaterreport}. Though corporate cafeterias are convenient and usually offer subsidized prices, they often face complaints of boredom due to the need for more dish variety. To address this, cafeteria managers manually diversify their menus based on their experience and industry knowledge. However, without a careful understanding of how previous menus impact customer demand, manual changes may be inefficient.}  

Many studies have explored and validated that customers' historical experiences may affect their willingness to buy or retain a product \citep{bowden2009process,lemon2016understanding}. For example, a customer may initially feel satisfied with a product, but this satisfaction may wane over time. This utility decay phenomenon is also observed in fashion products \citep{caro2014assortment}. On the other hand, a customer may also prefer to consistently consume ``traditional and familiar" food that could maintain or increase her utility~\citep{edenred}. Today, with the increasing collection of transactional data, companies seek to leverage these data to estimate customers' utility changes and build data-driven assortment planning systems. For instance, new technology companies have emerged that revolutionize traditional catering services by analyzing transactional data to estimate customers' preferences and offering optimized daily cafeteria menus~\citep{fooda}. 

%For further details and references in this regard, the reader is referred to Fooda's Cafeteria Solution\footnote{\url{https://www.fooda.com/cafeteria-solutions-corporate-dining}} and Zerocater's Corporate Cafeteria Solution.\footnote{\url{https://zerocater.com/corporate-cafeterias/}} 

% Fitness centers must plan, in advance, diverse workout classes to maintain their clientele \citep{lei2023content}. 

Besides corporate cafeterias, analogous situations arise in various other sectors. Flash sales platforms must organize their promotional assortments ahead of time, based on the products suppliers provide~\citep{MartnezdeAlbniz2020UsingCD}. For instance, Meituan, a leading Chinese online shopping platform, manages a grocery flash sale channel. Due to the lead time of contract negotiations with suppliers, Meituan needs to decide, in advance, assortments for daily deals over a fixed planning horizon. Given that customers frequently visit the channel for discounted groceries, planning must consider historical offerings. More details about Meituan's flash sale campaign are in the ecompanion~\ref{sec:ec:meituan}. Managers of cafeterias and flash sales encounter a common multi-period planning problem that decides how to offer identical assortments to all customers over a fixed planning horizon. In practice, these assortments are usually predetermined before the planning horizon, possibly due to limits (such as lead time) from the supplier side.

%fresh food, dairy, beverages, etc.  offer about 20 fresh products at deep discount on sale every day with a limited time window. Meituan schedule the flash sale assortments in advance based on contracts with third parties \citep{MartnezdeAlbniz2020UsingCD}. Attracted by the low prices, customers may view the assortment in the flash sale campaign every day. More details about Meituan's flash sale campaign are in Appendix~\ref{eg:meituan}.  

% Live stream sellers or traditional TV home sellers also build a calendar containing various types of goods on different days to attract fans consecutively watching and avoid their fatigue due to repeated live shows of the same product.\footnote{We give an example of the live stream show calendar in Appendix~\ref{eg:calendar}.}

In this paper, we focus on the above applications and investigate the corresponding \textit{ multi-period assortment planning problem considering history-dependent utilities} to maximize the long-run average revenue.
The history-dependent effects bring new challenges for assortment planning because of the interplay between cross-product and cross-period effects. Cross-product effects refer to the influence of other products in an assortment on purchase decisions. For instance, customer substitution behavior is commonly modeled using discrete choice models, like the multinomial logit (MNL) model~\citep{mcfadden1974}. Cross-period effects refer to the temporal influences on customer utility, i.e., the utility of a product in a given period depending on its offering history. Given both types of effects, the offering decisions for multiple products over multiple periods are made simultaneously. The number of candidate assortments, therefore, scales exponentially with the number of products and the planning horizon. This can be enormous for even a practically-sized problem. Moreover, customer behavior under both types of effect makes the decision more complex because, for a given period, the purchase probability is a nonlinear function---specifically, a softmax function---of the history-dependent utility functions. 

\begin{figure}[hbt]
   \FIGURE
   {\includegraphics[width=0.8\linewidth]{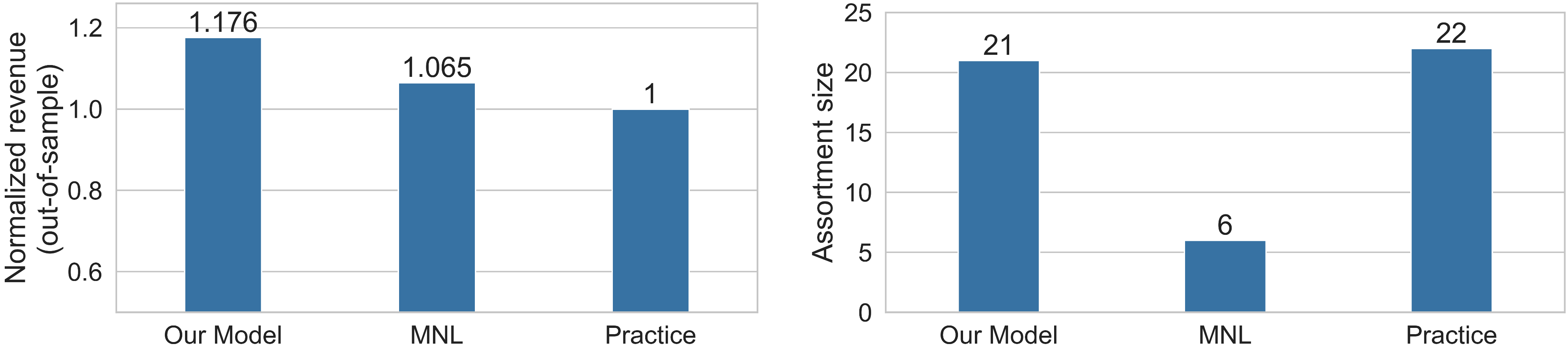}}
   {\re{The average out-of-sample results of three assortment plans}\label{fig:intro:case}}
  {Our model is the history-dependent model. We normalize the revenue generated by the current practice to one. The assortment size is the number of product types offered in the planning horizon.}
\end{figure}

Despite the complexity of assortment planning with history-dependent effects, the benefit of considering such effects is significant. Using the data from a corporate cafeteria with which we cooperated, we evaluate assortments derived from a history-dependent model (our model, estimated in  Section~\ref{sec:case:study}), from the MNL model, and the current practice observed in the data (as unit one for normalization). The results shown in Figure \ref{fig:intro:case} indicate that our model achieves the highest revenue in the out-of-sample test, \re{about $10.4\%$ more than the revenue generated by the MNL model.} Our model also produces assortments that have greater variety than those obtained by the MNL model. The current practice decided by managers resulted in assortments with a greater variety but the lowest revenue. See Section \ref{sec:case:study} for details of the case study.

\subsection{\re{Description of Results}}
We study a multi-period assortment planning problem in which a firm facing both cross-product and cross-period effects needs to decide what to offer in each period to maximize the average expected revenue. We model customer choice in a given period using a variant of the MNL model in which the utility of each product linearly depends on its past offering decisions. In our model, an assortment decision for a given period affects both current and future purchase behavior in a nonlinear fashion. In other words, our assortment planning problem is a nonlinear integer programming problem.

We show that this problem is NP-hard under the presence of negative history-dependent effects and propose exact solution methodologies \re{for a general setting that the history-dependent effects could be a mixture of addiction and satiation}. We propose a lifting-based reformulation framework to obtain a mixed-integer exponential cone program (MIECP) formulation. The key idea behind our formulation is to lift the original nonconvex constraints into a higher-dimensional space where the lifted constraints are simpler to model. This reformulation allows us to utilize prevalent techniques in the optimization literature, such as the McCormick envelope, the perspective formulation, and the Lov\'asz extension, to model the lifted structures. Numerical studies indicate that our formulation efficiently achieves optimality and yields higher revenue than heuristic policies. We also identify an optimal sequential revenue-ordered (RO) policy if history-effects are all positive, i.e., addiction effects. This characterization indicates that our problem with the addiction effect is solvable in polynomial time. For problems with a long-term planning horizon, theoretically, we show that a cyclic policy is optimal and find that, under a certain non-overlapping condition, the optimal cycle length equals the customer memory length plus one. Computationally, finding such policies is NP-hard but also admits an MIECP formulation. Moreover, the non-overlapping condition leads to a mixed-integer linear programming (MILP) formulation that is stronger than the conic one.

We conduct a case study using transactional data from a corporation cafeteria operated by a high-tech catering service company in China. The cafeteria primarily serves workday lunches for employees who work nearby. We estimate customers' utility function and then design weekly day-by-day lunch menus. The estimation shows that a satiation effect does exist, but it only lasts for a short period (i.e., two or three days). We demonstrate that our model can provide menus that balance variety and revenue, while current industry practices sacrifice revenue, and the classic MNL model yields menus with limited variety and low revenue.

We highlight the contributions of our paper as follows: 
\begin{itemize}
    \item[1.] Our model, which we refer to as \textsc{HAP} (i.e., history-dependent assortment problem),
simultaneously incorporates the cross-product and cross-period effects into the multi-period assortment planning process. \re{The cross-period history-dependent effects could be all positive, all negative, and mixed across products and periods.} We demonstrate, via a case study of a corporate cafeteria, that our model has several benefits. First, our result suggests that the service provider should provide a greater variety of products if satiation effects exist. When satiation effects last longer, it is beneficial to offer more types of products across periods. Furthermore, it is recommended that a product be offered only once within the memory length to prevent a utility decrease due to satiation effects, particularly in the case of products with high prices. This observation naturally leads to a ``non-overlapping" cyclic assortment planning policy, which we explore in Section~\ref{sec:two:cycle}, where we establish its optimality. Most importantly, our history-dependent policy effectively balances variety and the revenue objective, resulting in assortments that achieve the highest revenue in out-of-sample tests while maintaining the same variety as found in practical scenarios.
    
    %our history-dependent assortment model provides assortments that generate the highest revenue (out-of-sample) while maintain the same level of product variety in practice. 

    %First, our model has better prediction performance than one without history-dependent effects. Second, our model effectively balances assortment variety and revenue, which lifts current revenue by more than $19\%$.
    % We conduct a case study on a catering service dataset to design weekly menus for a company cafeteria. We show that our model is well-fitting and effectively balances variety and revenue, a key concern for any catering service operator. We also demonstrate the impact of the negative history-dependent effect on a menu's variety and offering pattern. highlight benefits of our model.
    
    \item[2.] We characterize, in terms of history-dependent effects, the conditions under which a sequential revenue-ordered policy globally solves our HAP model. In general, it can be shown that finding a global optimal solution is NP-hard. To tackle this difficulty, we develop an MIECP formulation of our \textsc{HAP} model by using lifting and coupling ideas from the optimization literature, including perspective formulation, the McCormick envelope, and the Lov\'asz extension. \re{In particular, we modify the Lov\'asz extension to model attraction values under general history-dependent effects.} Our formulations, \re{together with cutting-plane algorithms}, outperform an MILP formulation of HAP obtained using a prevalent strategy in the mixed-integer nonlinear programming (MINLP) literature and the original nonlinear integer formulation HAP solved by the state-of-the-art MINLP solver~\texttt{SCIP}. 
    %We further design a double-cutting-plane algorithm for the issue that the number of constraints in MIECP increases fractionally with memory length to strengthen its computation performance.
    \item[3.] We establish the \re{asymptotic} optimality of cyclic policies, and finding such policies allows an MIECP formulation. Our study reveals a connection between the length of an optimal cyclic policy and the length of customer memory. In particular, we characterize the condition under which the optimal cyclic length is the customer memory length plus one. This condition, in addition, leads to an MILP formulation for finding \re{optimal cyclic policies}, which is, theoretically, tighter than the previous conic formulation.
    % \item[4.] \re{We design two algorithms to strengthen the computation of our formulations. First, we design a double-cutting-plane algorithm to address the issue that the number of constraints in MIECP increases fractionally with memory length. Second, we construct a projected-cutting-plane algorithm to implement the strong MILP model, which includes can solve instances with $1000$ products within $800$ seconds.}
\end{itemize}
  
\subsection{\re{Literature Review}}\label{sec:literature}
Extensive research has incorporated history-dependent effects into operations problems, including pricing and promotion~\citep{bi2020promotion,chen2023customized}, service design~\citep{bernstein2022intertemporal,lei2023content}, resource allocation~\citep{adelman2013dynamic}, retention management~\citep{kanoria2023managing}, etc. Studies show that ignoring history-dependent effects may lead to substantial losses. For instance, considering the reference price can reduce $5.82\%$ prediction error of sales \citep{wang2018prospect}, and the optimal bundle accounting for customer satiation improves the revenue by more than $4.5\%$ \citep{gurlek2023optimal}. 

History-dependent effects are receiving increasing attention in assortment planning. One stream of work is to design personalized assortment planning for online platforms. The history-dependent effect derives from repeated interactions between the firm and customers, which breaks the traditional single-interaction assumption in personalized assortment optimization literature~\citep{el2023joint}. Based on the easily accessible purchase history of a customer, these studies consider two settings: (i) customers have stable repeated interactions with the firm, e.g., \cite{chen2023assortment} optimize adaptive assortment decisions based on information inferred from a customer's historical purchases; (ii) a customer can leave or churn if he or she is not satisfied. Hence, many studies focus on designing delicate exploration-exploitation recommendation algorithms that may limit exploration budgets \citep{bastani2022learning,sumida2023optimizing}.

For non-personalized assortment planning, the history-dependent effect is captured by a dynamic choice/demand function of the entire market. The effect may derive from attributes of a product itself, e.g., declining utility. Such intertemporal connections lead to a nonlinear and high-dimensional problem. Two settings are considered: (i) an adaptive assortment policy, e.g., \cite{caro2020managing} consider a content release problem where the attraction of content decays and the number of followers stochastically evolves; (ii) a non-adaptive assortment plan for the entire horizon made before the first period, e.g., \cite{caro2014assortment} optimize release schedule of fast fashion products with exponentially declining attraction before period one due to a long lead time. Our work contributes to the second type of literature. Our model can represent the model in \cite{caro2014assortment} and adds more flexibility (see details in Section~\ref{sec:discuss}), which brings new challenges and needs additional methods.

\re{Our paper contributes to the growing literature on developing new solution methods for assortment planning. The three primary approaches include the revenue-ordered policy, approximation algorithms, and integer programming for obtaining exact solutions.} Since~\cite{talluri2004revenue} showed that the optimal assortment under the MNL model is revenue-ordered, this policy has been extensively studied in various contexts~\citep{wang2017consumer,sumida2023optimizing}. For example,~\cite{gao2021assortment} prove that the revenue-ordered assortment policy is optimal under the MNL model when the offered assortment is gradually revealed.~\cite{xu2023assortment} exploit the revenue-ordered property to solve a multistage assortment optimization problem. We contribute to this literature by showing that a sequential revenue-ordered policy is optimal for our problem under the addiction effect. 

\re{A stream of studies focuses on designing approximation algorithms for assortment plannings. For instance, \cite{feldman2015capacity} and \cite{desir2020mc} design constant factor approximation algorithms for capacitated assortment optimization under the nested logit model and the Markov chain choice model, respectively. Later, \cite{desir2022capacitated} propose a fully polynomial time approximation scheme (FPTAS) for capacitated assortment problems under the mixed MNL, Markov chain, and nested logit choice models. FPTAS are also proposed to solve assortment optimization under the exponential choice model~\citep{aouad2023exponomial}, with online dynamic behavior~\citep{liu2020assortment,feldman2023display,aouad2024click}, and with multi-purchase behavior~\citep{jasin2024assortment,luan2025joint}.}

% including the fully polynomial time approximation schemes (FPTAS)~\citep{aouad2024click,jasin2024assortment}, PTAS, polynomial factor approximation algorithms \citep{caro2014assortment,gallego2023constrained}, and constant factor approximation algorithms~\citep{aouad2019approximation,chen2023assortment} 

\re{Our approach is mostly related to the third direction --- using integer programming techniques to obtain global optimal assortment decisions.} A stream of studies establishes MILP models for assortment planning under various choices models, including the ranking-based choice model~\citep{bertsimas2019exact}, the mixture of Mallows model~\citep{desir2021mallows}, the decision forest choice model~\citep{akchen2023assortmentoptimizationdecisionforest}, and the mixture-of-nested-logit model~\citep{fan2024assortmentvideogame}. In particular,~\cite{chen2024integer} conduct polyhedral studies for the MNL-based assortment planning in the context of quick-commerce. Mixed-integer conic programming has also been used to model choice behavior.~\cite{sen2018conic} propose a mixed-integer second-order cone programming formulation for assortment problems under the mixed MNL model, and \cite{chen2022offline} use this idea to solve location-dependent offline-channel assortment planning in omnichannel retailing.~\cite{akccakucs2023exact} develop various mixed-integer exponential cone programs to solve a share-of-choice product design problem where customers follow a logit-based choice model. Our solution approach is inspired by formulations in~\cite{akccakucs2023exact}. However, due to the interaction of cross-product and cross-period effects, additional ideas are required for our problem, including such relaxation techniques as McCormick envelope~\citep{mccormick1976computability} and the Lov\'asz extension~\citep{lovasz1983submodular}. 
% With these techniques, we obtain an equivalent mixed-integer exponential cone formulation widely used in many recent studies \citep[e.g.][etc]{zhu2022joint,chen2023exponential}.

Last, we comment that our paper contributes to a growing set of applications of conic programming/discrete optimization in designing operations systems, including joint inventory-location problems~\citep{atamturk2012conic}, battery swap networks design~\citep{mak2013infrastructure}, service delivery scheduling~\citep{kong2013scheduling}, process flexibility design~\citep{simchi2012understanding,yan2018design}, electric vehicle charging~\citep{chen2023exponential}, proactive policing~\citep{he2023proactive}, network flow~\citep{simchi2019constraint}, and nonlinear resource allocation~\citep{he2024discrete}. 

% Our paper is also related to an arising topic on applying convex programs for operation management problems, including facility location design~\citep{shen2024data,xu2024ensemble}, scheduling~\citep{kong2013scheduling,zhang2022schedule}, pricing~\citep{yan2022representative,he2024discrete}, and robust optimization~\citep{simchi2019constraint,zhu2022joint,natarajan2023distributionally}. Our study applies lifting framework and obtains an mix integer exponential conic exact formulation for a class of temporal assortment planning problems. Our formulation can obtain global optimal assortments efficiently
% \citep[e.g.][etc]{zhu2022joint,chen2023exponential}.

% ~\cite{bertsimas2019exact} propose a new integer optimization formulation of the product line design problem when the choice model is ranking-based, and solve it using Benders decomposition. When customers choose according to the mixture of Mallows model, ~\cite{desir2021mallows} exploit structures of choice probabilities to derive an MIP formulation for the assortment optimization problem.~\cite{akchen2023assortmentoptimizationdecisionforest} develop mixed-integer programming formulations for assortment optimization under the decision forest model.
% ~\cite{fan2024assortmentvideogame} derive a MIP formulation for an assortment optimization problem arising in a class of video games where the choice model is a constrained mixture-of-nested-logit model.

\subsection{\re{Structure}}
The rest of the paper is organized as follows. In Section~\ref{sec:model},  we formally define the assortment planning problem and establish its computational complexity. In Section~\ref{eq:solutions}, \re{we propose a mixed-integer exponential cone formulation that can handle any mixed addiction-satiation history-dependent effects. We prove that, under positive effects, a revenue-ordered policy is optimal.} In Section~\ref{sec:two:cycle}, we present our results on cyclic policies. We present a case study using real data in Section~\ref{sec:case:study} and report numerical results on synthetic data in Section~\ref{sec:numerical}. We conclude the paper in Section~\ref{sec:conclusion}. All proofs are given in the ecompanion.

\section{Model}\label{sec:model}
\re{In Section~\ref{sec:sub:model}, we define our multi-period assortment planning problem with history-dependent effects and establish its computational complexity. In Section~\ref{sec:discuss}, we discuss the model features in detail.}

\subsection{\re{Incoporating History-Dependent Effects into Choice Models}}\label{sec:sub:model}   
Consider a firm that sells $N$ products to customers. Let $[N]:=\{1, \ldots, N\}$ denote the set of all products, and let $[N]^+:= [N] \cup \{0\}$, where the index $0$ denotes the outside option. For each product $ i \in [N]$, let $r_i$ denote the revenue obtained from selling an instance of product $i$. There is a total of $T$ planning horizons, and the firm needs to decide what to offer in each period. More specifically, for each period $t \in [T]$, let $\x^t = (x^t_{1}, \ldots, x^t_{N}) \in \{0,1\}^N$ be a binary decision variable modeling the subset of products offered in period $t$, i.e., $x^t_i = 1$ if and only if product $i$ is offered in period $t$. We assume that in each period, one unit mass of customers arrives, and the customers select an item from the offered products or leave without a purchase. For example, in a setting such as corporate cafeterias, the firm serves relatively stable customers (e.g., most customers are employees who work nearby), and customers make purchase decisions on a daily basis.

We modify the MNL model to encapsulate the impact of the historical assortments on customer choices. Each product has a utility $U_i$, which depends on attributes of the product and environment. Following \cite{guadagni1983logit,wang2018prospect} and \cite{xie2024personalized}, we assume that the utility of one product is a linear function of its attributes and historical offerings. Specifically, the utility of product $i$ in period $t$ affinely depends on its historical offering records within a memory length $M$:
\begin{equation*} \label{eq:utility}\tag{\textsc{Utility}}
U_i^t = \beta_i^0 + \sum_{m\in [M]} x^{t-m}_i \beta_{i}^m , 
\end{equation*}
where $M \in \Z_{+}$ denotes the customer memory length. For simplicity in the presentation, we assume throughout this paper that $x_i^k:=0$ for $k\leq 0$. The parameter $\beta^0_i$ is the base utility of the product $i$ determined by its attributes. The parameter $\beta_{i}^m$ is the incremental utility if product $i$ is offered in the $m^{\text{th}}$ period before the current time. Let $\boldb := (\beta_i^m)_{i\in[N],m\in[M]}$ denote the history-dependent effect matrix. \re{The value of $\beta_{i}^m$ can be positive or negative. A positive $\beta_{i}^m$ indicates that exposure to product $i$ increases its utility, acting as a stimulus \citep{fox1997modeling}. While a negative $\beta_{i}^m$ suggests a satiation effect, leading to customer fatigue from repeated offerings. The history-dependent effects can vary across products and periods.} In our planning model, the assortments are fixed for each period. Therefore, we do not consider personalized assortments or model individual-level utility functions.

%For instance, some products may exhibit addictive behavior, while others show a trend of satiation. A product initially has an addictive effect but then later shifts to satiation.} 

For a given assortment plan $\x:=(\x^1, \cdots, \x^T)$, we use $\pi^t_i(\x^t, \ldots, \x^{t-M})$ to denote the purchase probability of product $i$ in period $t$. Under the MNL model, the purchase probability of product $i$ in period $t$ is 
\[
\begin{aligned}
\pi^t_{i}(\x^t, \ldots, \x^{t-M})&= 
    \dfrac{ x^t_i\exp\bigl(\beta_i^0 + \sum_{m\in [M]} x^{t-m}_{i} \beta^m_{i}\bigr) }{1+ \sum_{j\in [N]}x^t_j \exp\bigl(\beta_j^0 + \sum_{m\in [M]} x^{t-m}_{j} \beta^m_{j}\bigr) } \quad \for i \in [N], \\
\pi^t_{0}(\x^t, \ldots, \x^{t-M}) &= 
    \dfrac{1  }{1+ \sum_{j\in [N]}x^t_j\exp\bigl(\beta_j^0 + \sum_{m\in [M]} x^{t-m}_{j} \beta^m_{j}\bigr) },
\end{aligned}
\]
where we normalize the utility of the outside option to 0. The purchase probabilities depend on two components: the current assortment $\x^t$ and the past assortments $\x^{t-1},\dots,\x^{t-M}$. The former dependence captures cross-product effects such as the substitution effect, while the latter dependence models cross-period effects. {The purchase probability can also be viewed as the market share of product $i$ in period $t$ and presents the effect of historical assortments on sales of the entire market~\citep{batsell1985new}.} 

With these definitions, assortment planning under history-dependent effects can be expressed as the following fractional binary program, called the history-dependent assortment problem: 
\begin{equation}\label{eq:total}
 \max \Biggl\{ \sumt  \frac{1}{T} \sumi r_i \pi^t_i(\x^t, \ldots, \x^{t-M}) \Biggm| \x^t  \in \mathcal{X} \cap \{0,1\}^N \, \for t \in [T]  \text{ and }  \x \in \mathcal{P}
    \Biggr\} , \tag{\textsc{HAP}}
\end{equation}
where $\mathcal{X}$ (resp. $\mathcal{P}$) is a system of linear inequalities that models cross-product constraints (resp. cross-period constraints). The cross-product constraints often include resource restrictions within a single period, such as cardinality and capacity constraints, and address various operational considerations in the assortment planning literature \citep{sen2018conic,chen2024integer}. Cross-period constraints describe the intertemporal relationships of assortments across different periods, such as non-overlapping assortments over consecutive periods, which occur in multi-period assortment planning problems \citep{liu2020assortment, chen2023assortment}.

\re{This model is generally difficult to solve. We show in the following proposition that ~\eqref{eq:total} is NP-hard under very simple settings if history-dependent effects are negative.} 
%In particular, when products have identical revenue and the planning horizon is two, we can reduce a NP-hard problem, the \textit{3/4-Partition} problem, to our problem \eqref{eq:total}.
\begin{proposition}\label{prop:nphard}
\eqref{eq:total} is NP-hard even when the planning horizon is two, the memory length is one, and history-dependent effects are negative.  
\end{proposition}

\subsection{Model Discussion}\label{sec:discuss}
We outline here the main features of our model, justify our assumptions, and identify relevant scenarios. 

It is apparent that our formulation requires a firm to present the same assortment to all customers in each period before the planning horizon begins. In other words, the firm is unable to take advantage of personalized assortments and adaptive policies, which are common practices in online assortment optimization such as Amazon's recommendation system. However, those practices are not suitable for the case of corporate cafeterias, where it is common to treat employees with the same menu, and it is impractical to dynamically adjust assortments due to the lead time in both procurement and cooking processes~\citep{van2006planning}. Observations from our partnered catering company indicate that the corporate cafeteria plans the lunch menus on a weekly basis. Moreover, our model assumption is also reasonable for flash sale channels (i.e., predetermine assortments for several days due to the negotiation time with suppliers).

Our model describes the history-dependent effect through a linear utility function with historical offerings as features. The linear utility function assumption is common in the logit-based choice model to capture impacts such as product attributes \citep{guadagni1983logit} and historical prices \citep{wang2018prospect}. Moreover, the linear utility function covers multiplication in the attraction value in literature. For example, our model with negative history-dependent effects represents the model in \cite{caro2014assortment}. Our model offers additional flexibility by allowing for \re{positive, negative, and mixed} history-dependent effects and offering products in any period.

In the linear utility function, we incorporate the historical offerings to capture the history-dependent effect. The utility function captures the direct effect of historical assortments that are critical decisions at a firm level. Moreover, we only consider the direct history-dependent effect, i.e., the historical offering of a product will only affect the product's future utility itself. We do not include the cross-product history-dependent effect, i.e., the historical offering of a product might affect the utility of other products, because this effect is usually at a much smaller scale \citep[e.g.,][]{vilcassim1999investigating,dube2005differences, mehta2007investigating}. In our empirical estimation, the cross-product history-dependent effect is approximately one-tenth of the direct effect and is significant only when $M=1$. We believe that our model efficiently captures the primary effect and maintains model tractability. 
 
We manage the history-dependent effects across periods through the memory length $M$, inspired by the patience level from dynamic pricing studies~\citep{popescu2007dynamic,liu2015optimal}. This patience level indicates how many periods a customer will wait for a price reduction. Consequently, it is logical to assume that customers might retain a limited memory length for a product.  
  
Our model does not include inventory decisions. In cafeteria settings, for instance, procurement choices can be made subsequent to assortment planning. Additionally, inventory rollover is not considered in the study because cooked dishes cannot be kept for the following day due to hygiene standards. To simplify the analysis, this study will concentrate on the assortment planning problem with cross-product and cross-period effects, which is already a challenging problem.

% \subsection{Hardness Result}\label{sec:model:complexity} 

\section{\re{Solution Approaches}}\label{eq:solutions}
\re{In this section, we aim to solve the history-dependent assortment problem~\eqref{eq:total} to its global optimality. In Section~\ref{sec:reformulation}, we propose a mixed-integer exponential cone formulation for the general case, and, in Section~\ref{sec:model:RO}, we characterize conditions under which a sequential revenue-ordered policy is optimal.}

\subsection{\re{Integer Programming Approach}}\label{sec:reformulation}
 Due to the combinatorial and nonconvex nature of~\eqref{eq:total}, state-of-the-art global optimization solvers such as \texttt{BARON}~\citep{tawarmalani2005polyhedral} and  \texttt{SCIP}~\citep{BestuzhevaEtal2021OO} fail to solve practical-sized instances---see Section~\ref{sec:numerical:speed} for computational evidence. Nevertheless, we focus on developing techniques for obtaining mixed-integer programming formulations of~\eqref{eq:total}, which will allow us to leverage modern mixed-integer program solvers such as \texttt{GUROBI} and \texttt{MOSEK}. Then, in Section~\ref{sec:two:cycle}, we will use these techniques to find asymptotically optimal cyclic policies.

Our formulation introduces the following additional variables to represent the customer choice behavior. For each period $t \in [T]$, we use $\rho^t$ to denote the no-purchase probability for period $t$:
\begin{equation*} 
\rho^{t} = \frac{1}{1 + \sum_{j \in [N]}x^t_{j} \exp\bigl(\beta_j^0 + \sum_{m\in [M]} x^{t-m}_{j} \beta^m_{j}\bigr)},  \label{eq:non-prob}
\end{equation*}
and for each product $i \in [N]$, we use $y^t_{i}$ to denote the choice probability for product $i$ in period $t$:
\begin{equation*}
     y^t_{i} = \frac{x^t_{i} \exp\bigl(\beta_i^0 + \sum_{m\in [M]} x^{t-m}_{i} \beta^m_{i}\bigr)}{1 + \sumj x^t_{j} \exp\bigl(\beta_j^0 + \sum_{m\in [M]} x^{t-m}_{j} \beta^m_{j}\bigr)}.  \label{eq:choice-prob}
\end{equation*}
In what follows, we refer to $\rho^t$ as the no-purchase probability variable and $y^t_i$ as the choice probability variable. Using these definitions, we can reformulate ~\eqref{eq:total} as follows:
\begin{subequations}\label{eq:total-1}
\begin{alignat}{3}
    \max_{\x,\y,\boldsymbol{\rho}} \quad & \frac{1}{T} \sumi  \sumt r_iy_{i}^t \notag \\
        \st \quad & \x^t \in \{0,1\}^{N } \cap \mathcal{X} \text{ and }  \x \in \mathcal{P}&& \for \, t \in [T]  \notag \\
    &\rho^t + \sumi y_{i}^t = 1  \quad && \for\, t \in [T] \notag \\
    & y^t_i = \rho^tx^t_{i} \exp\Bigl(\beta_i^0 + \sum_{m\in [M]} x^{t-m}_{i} \beta^m_{i}\Bigr) \quad  && \for \, t \in [T] \text{ and } i \in [N]. \tag{\textsc{Choice}} \label{eq:choice}
\end{alignat} 
\end{subequations}
The objective function is a non-negative weighted combination of choice probability variables. The second constraint ensures that the sum of the choice probability and the no-purchase probability equals one for each period. The last constraint represents the choice probability for product $i$ in period $t$, and thus we refer to it as constraint ~\eqref{eq:choice}. This reformulation encapsulates the nonconvexity of problem ~\eqref{eq:total} in constraint~\eqref{eq:choice}. Hence, it suffices to derive a convex representation of constraint~\eqref{eq:choice}.  

The main component of constraint \eqref{eq:choice} is the composition of a univariate exponential function and an affine utility function, defined as follows:
\[
\alpha_{i}(z^1, \ldots, z^{M} ) := \exp\biggl(\beta_i^0 + \sum_{m\in [M]}  \beta^m_{i} \cdot z^{m}\biggr) \qquad \for\; (z^1,\ldots, z^{M}) \in \{0,1\}^M.
\]
For a given assortment plan $\x:=(\x^1, \ldots, \x^T)$, we interpret \re{$\alpha_i(x_i^{t-1}, \ldots, x_i^{t-M})$} as the attractiveness of product $i$ in period $t$. Therefore, we will henceforth refer to  $\alpha_i(\cdot)$ as the attraction value function. To handle the discrete nature of the attraction value function, we will use the notion of a \textit{continuous extension}. Given a function $f: \{0,1\}^n \to \bbmr$,  a function $g:[0,1]^n \to \bbmr$ is called a continuous extension of $f(\cdot)$ if $g(\cdot)$ is continuous  and $g(\z) = f(\z)$ for every $\z \in \{0,1 \}^n$. In addition, if an extension function $g(\cdot)$ is convex (resp. concave), we say that it is a convex (resp. concave) extension. It is common for a discrete function to admit many convex (resp. concave) extensions. A natural idea is to use the tightest convex/concave extension. Since the discrete domain $\{0,1\}^n$ is the vertex set of the continuous domain $[0,1]^n$, the pointwise largest convex (resp. smallest concave) extension of $f(\cdot)$ is the \textit{convex (resp. concave) envelope}, denoted as $\conv(f)(\cdot)$ (resp. $\conc(f)(\cdot)$), of $f(\cdot)$ over the box $[0,1]^n$~\citep[see][Theorem 6]{tawarmalani2002convex}. 
%Recall that the convex (resp. concave) envelope of $f(\cdot)$ is referred to as the tightest under (resp. over)-estimator of $f$ over $[0,1]^n$~\citep{rockafellar1997convex}.
 
\def \vf {f}

\subsubsection{\re{A Lifted Representation of~\eqref{eq:choice}}}\label{section:envelope-characterization}
  
In this subsection, we focus on modeling constraint~\eqref{eq:choice}. The idea is to lift the original nonconvex constraint~\eqref{eq:choice} into a higher-dimensional space where the lifted constraint is dramatically simpler to represent than the original constraint. Specifically, this will allow us to invoke prevalent techniques, such as the McCormick envelope~\citep{mccormick1976computability}, the perspective formulation~\citep{gunluk2010perspective}, and the Lov\'asz extension~\citep{lovasz1983submodular}, to model the lifted nonconvex structures. 

% We elaborate on this in the following Sections~\ref{section:pers} and~\ref{section:supermodular}. 

 \textbf{Lifting for Simplicity} The challenge in lifting is to identify a proper lifting space. To lift, we adapt a construction that~\cite{akccakucs2023exact} use to solve the logit-based product design problem, which we briefly review here. Consider a nonconvex constraint given as $\lambda + \lambda f(\x) \leq 1$,  where $\lambda \geq 0, \x \in  \bbmr^n$, and $f(\cdot)$ is a convex function. This constraint can be relaxed into a convex constraint as follows. First, we multiply  and divide $\x$ by $\lambda$ to obtain
\[
\lambda + \lambda f\Bigl(\frac{\lambda \cdot\x}{\lambda}\Bigr) \leq 1,
\] 
where $\lambda f(\frac{ \lambda \cdot \x}{\lambda}) = 0$ when $\lambda=0$. Next, we replace $\lambda \x$ with a new variable vector $\u$. This leads to 
\[
\u = \lambda \x  \quad  \text{and} \quad  \lambda + \lambda f\Bigl(\frac{\u}{\lambda} \Bigr) \leq 1,
\]
where the first constraint can be linearized using the McCormick envelope~\cite{mccormick1976computability} and the second  constraint is convex, since the function $\lambda f(\u/\lambda)$ is the \textit{perspective function} of the convex function $f(\cdot)$~\citep{rockafellar1997convex}. 
Formally, given a function $f: D \to \bbmr$, where $D$ is a bounded subset of $\bbmr^n$,  the \textit{perspective function} of $f$, denoted as $\pers(f)$, is defined as 
\[
\pers(f)(\lambda,\x) =  \begin{cases}
\lambda	f(\frac{\x}{\lambda})  & \text{ for } \lambda >0 \text{ and } \x \in \lambda \cdot D \\
0 & \text{ for } \lambda = 0\text{ and }\x = \boldsymbol{0}. 
\end{cases}
\]
Perspective reformulations have been widely used in solving nonconvex optimization problems~\citep{gunluk2010perspective}; see~\cite{bestuzheva2023computational} for a review of applications and a detailed computational study of perspective reformulations. 
 
We now present our lift of the constraint~\eqref{eq:choice}. Note that the choice probability of product $i$ in period $t$ is computed by scaling its attraction value function $\alpha_{i}(\cdot)$ with a nonlinear factor $\rho^t x_{i}^t$, that is, 
\[
y^t_i = (\rho^tx^t_i) \cdot \alpha_i(x_i^{t-1}, \ldots, x_i^{t-M}).
\]
First, we introduce a new variable $\gamma_{i}^t$ to represent the nonlinear factor $\rho^t x_{i}^t$. Second, we use the new variable $\gamma_{i}^t$ to scale the historical assortment decisions $x^{t-m}_i$ with $m \in [M]$, and we introduce the variable $z_{im}^{t}$ to represent the scaled variable $\gamma_{i}^t x^{t-m}_i$. Using the additional variables  $\gamma_{i}^t$ and $\z_{i}^t = (z_{i1}^{t}, \ldots, z_{iM}^{t})$ and the definition of a perspective function,  we decompose constraint~\eqref{eq:choice} into constraints
\begin{subequations}\label{eq:extended}
\begin{alignat}{3}
& \gamma_{i}^t = \rho^t x_{i}^t \label{eq:extended-1} \\
& z_{im}^t = \gamma_{i}^t x^{t-m}_{i} \quad \for \, m \in [M] \label{eq:extended-2}\\
&y^t_i = \pers(\alpha_{i})( \gamma^t_i , \z_{i}^t).\label{eq:extended-3} 
\end{alignat}
\end{subequations}
Constraint~\eqref{eq:extended-1} lifts the nonlinear scaling factor $\rho^t x_{i}^t$ to $\gamma_{i}^t$; constraint~\eqref{eq:extended-2} lifts the scaled historical assortment decision $\gamma_{i}^t x_{i}^{t-m}$ to $z_{im}^t$; and constraint~\eqref{eq:extended-3} perspectifies the attraction value function in the space of the newly introduced variables $\gamma_{i}^t$ and $\z_{i}^t$. 

To obtain a mixed-integer convex formulation of constraint ~\eqref{eq:choice} in the lifted space, we individually relax each of the three constraints in~\eqref{eq:extended}.  Specifically, we use the McCormick envelope~\citep{mccormick1976computability} to relax constraint~\eqref{eq:extended-1}. To apply the McCormick envelope to our setting, we derive a lower bound $\rho^t_L$ and upper bound $\rho^t_U$ on the no-purchase probability variable $\rho^t$, and impose the following constraints: 
\begin{equation}\label{eq:relax1}
\begin{aligned}
 \gamma_{i}^t &\leq \min\bigl\{\rho^t_L x_{i}^t + \rho^t-\rho^t_L,\ \rho^t_U x_{i}^t  \bigr\}  \\
 \gamma_{i}^t & \geq \max \bigl\{ \rho^t_L x_{i}^t,\ \rho^t_U x_{i}^t + \rho^t - \rho^t_U \bigr\}.
\end{aligned}
\end{equation}
Similarly, using $0$ (resp. $\rho^t_U$) as a lower (resp. upper) bound on $\gamma^t_i$, the McCormick envelope reformulates constraint~\eqref{eq:extended-2} as follows: 
\begin{equation}\label{eq:relax2}
\begin{aligned}
 z_{im}^t &  \leq \min \bigl\{ \gamma_{i}^t , \rho^t_U x^{t-m}_i  \bigr\} && \quad \for \, m \in [M]  \\
 z_{im}^t & \geq \max \bigl\{ 0,\ \rho^t_U x^{t-m}_i + \gamma^t_i - \rho^t_U \bigr\} && \quad \for \, m \in [M].
\end{aligned}
\end{equation} 

Finally, we represent constraint~\eqref{eq:extended-3} by using two continuous extensions of the attraction value function $\alpha_i(\cdot)$. On one side, a natural continuous extension $\tilde{\alpha}_i: [0,1]^{M} \to  \bbmr$ can be obtained from $\alpha_i(\cdot)$ by relaxing the binary assortment variables to continuous ones. This continuous function is convex since it is the composition of a convex function and an affine function. Using this extension, we relax  constraint~\eqref{eq:extended-3} to 
\begin{equation*}
    y^t_{i} \geq \pers(\tilde{\alpha}_{i})(\gamma^t_{i},\z^t_{i}).
\end{equation*}
This is a convex inequality since $\tilde{\alpha}_{i}(\cdot)$  is a convex function and so is its perspective function~\citep{rockafellar1997convex}. Moreover, this constraint can be represented as an exponential cone constraint \citep{ben2001lectures} (see more details in the ecompanion \ref{thm:proof:miecp}):
\begin{equation}\label{eq:conv-pers}
    y_i^t \ge  \exp(\beta^0_i) \cdot \gamma^t_i \cdot \exp \biggl( \frac{ \summ \beta_i^m z_{im}^t}{\gamma^t_i}\biggr),
\end{equation}
where the right-hand-side function can be regarded as a convex under-estimator of the choice probability variable $y^t_i$. On the other side, an exact representation of~\eqref{eq:extended-3} requires us to overestimate the choice probability variable. To do this, we replace $\alpha_i(\cdot)$ with its tightest concave extension or, equivalently, concave envelope. Using the concave envelope, we arrive at another convex relaxation of~\eqref{eq:extended-3}: 
\begin{equation}\label{eq:conc-pers}
    y_{i}^t \leq \pers\bigl(  \conc(\alpha_{i})\bigr)(\gamma^t_i, \z^t_{i}) .
\end{equation}
Next, we show that an explicit linear description of this constraint can be obtained from perspectifying a variant of the Lov\'asz extension of $\alpha_i(\cdot)$~\citep{lovasz1983submodular}.

\re{\textbf{The Concave Envelope of $\alpha_{i}(\cdot)$}}
%The supermodularity allows us to use the \textit{Lov\'asz extension} to describe concave envelope of $\alpha_i(\cdot)$~\citep{lovasz1983submodular,tawarmalani2013explicit,atamturk2022submodular}. 
Our construction is based on the Lov\'asz extension of a set function. The Lov\'asz extension is an extension of a function $f$ defined on $\{0,1\}^n$ to a function $\hat{f}$ defined on $[0,1]^n$. For each subset $S$ of $[n]$, let $\chi^S\in \{0,1\}^n$ be its indicator vector---that is, the $i^{\text{th}}$ coordinate of $\chi^S$ is $1$ if and only if $i \in S$. Observe that every vector $\x \in [0,1]^n$ can be expressed uniquely as $\x =  \lambda_0  \chi^{T_0} + \lambda_1  \chi^{T_1} + \cdots + \lambda_n \chi^{T_{n}} $,
where $\lambda_k \ge 0$ for $k = 0,1, \ldots,n $ and $\sum_{k =0}^n\lambda_k =1$, and $\emptyset = T_0 \subseteq T_1 \subseteq \cdots \subseteq T_n = [n]$. Thus, 
\[
\hat{f}(\x) = \lambda_0 f(\chi^{T_0}) + \lambda_1 f(\chi^{T_1}) + \cdots + \lambda_n f(\chi^{T_n}) 
\]
is a well-defined extension of the function $f$ (called the Lov\'asz extension of $f$) on the continuous domain $[0,1]^n$. The Lov\'asz extension serves as a bridge between supermodularity and concavity: $f$ is supermodular if and only if its Lov\'asz extension $\hat{f}$ is concave~\citep{lovasz1983submodular}. Moreover, when the Lov\'asz extension $\hat{f}$ is concave, it is expressible as the pointwise minimum of affine functions~\citep{tawarmalani2013explicit}.

% \begin{figure}[hbtp]
% \centering\includegraphics[width=0.8
% \linewidth]{Fig/switch_h.png}
%     \caption{\re{Adjust the nested vectors $(\chi^0,\chi^1, \chi^2)$ from $\bigl((0,0)^T, (0,1)^T,(1,1)^T\big)$ (left) to $\bigl((0,1)^T, (0,0)^T,(1,0)^T\bigr)$ (right) when $\beta_i^1<0$ and $\beta_i^2>0$ }}\label{fig:switch} 
% \end{figure} 

\re{When the history-dependent effects are negative, the attraction value function $\alpha_i(\cdot)$ is supermodular, and its concave envelope is given by its Lov\'asz extension~\citep{tawarmalani2013explicit}. The key idea behind the closed-form description of the Lov\'asz extension is to interpolate the attraction value function affinely over the nested historical assortments $\chi^{T_0}, \chi^{T_1}, \ldots, \chi^{T_M}$, where $T_0 \subseteq T_1 \subseteq \cdots \subseteq T_M$. These nested assortments can be thought of as a path through the set of all possible past assortments $\{0,1\}^M$ going from $(0,0, \ldots,0 )$ to $(1,1, \ldots, 1) $, with each step along one coordinate direction.  For example, when the memory length is two, the nested assortments $(0,0),(1,0)$, and $(1,1)$ correspond to the blue path in Figure~\ref{fig:switch1}. This path starts at $(0,0)$, moves along the first coordinate, then the second, and finally reaches $(1,1)$. To treat the general case of mixed effects, we modify the nested assortments---or equivalently the path---as follows: Instead of moving from $(0,0, \ldots,0 )$ to $(1,1, \ldots, 1)$, we start at $\chi^{I}$ and end at $\chi^{M \setminus I}$, where $I$ denotes the indices of positive effects. For instance, consider mixed effects with $\beta_i^1>0$ and $\beta_i^2\le 0$. The blue path in Figure~\ref{fig:switch2}, starting from $(1,0)$ and ending at $(0,1)$, corresponds to the past assortments that we will use in this case.}

% \re{
% One of the key steps in constructing Lov\'asz extension is the nested indicator vectors $\chi^{T_0}, \chi^{T_1}, \ldots, \chi^{T_n}$, where $T_0 \subseteq T_1 \subseteq \cdots \subseteq T_n$. In geometry, the nested vectors form a path in the cube $\{0,1\}^n$, starting from $\boldsymbol{0}$, passing through nested vertices, and finally terminating at $\boldsymbol{1}$. For instance, when the memory length is $2$ and the historical assortments satisfies $x_i^{t-1}\ge x_i^{t-2}$, the three nodes $(0,0),(1,0)$, and $(1,1)$ form a blue path in Figure~\eqref{fig:switch1} and their convex combination can represent all nodes in the lower right triangle. In the general case, the collection of vectors to present the historical assortments are not nested, but there still exists a path. We obtain such a path by switching elements in the positive dimensions of the nested vectors. More specifically, in our example, for the mixed effects $\beta_i^1>0$ and $\beta_i^2\le 0$, we switch the first element of $(0,0),(1,0),(1,1)$ from 1 (resp. 0) to 0 (resp. 1), and obtain $(1,0),(0,0),(0,1)$ forming the blue path in Figure~\eqref{fig:switch2}. Meanwhile, we switch the first dimension of the historical assortments from $x_i^{t-1}$ to $1-x_i^{t-1}$, as shown in Figure~\eqref{fig:switch2}. Then, the convex combination of the three nodes $(1,0),(0,0)$ and $(0,1)$ can present the switched historical assortments $(1-x_i^{t-1},x_i^{t-2})$ in the lower left triangle in Figure~\eqref{fig:switch2}, which satisfies $x_i^{t-1}\ge x_i^{t-2}$.}

\begin{figure}[hbtp]
\caption{\re{Illustration of switched nested assortments under two different history-dependent effects}}\label{fig:switch}
\hspace{0.1cm}
  \begin{subfigure}[b]{0.4\textwidth}
  \centering
  \caption{\re{Satiation effects: $\beta_i^1 \le 0 , \beta_i^2 \le 0$ } }
    \label{fig:switch1} 
   \resizebox{0.5\textwidth}{!}{\begin{tikzpicture}[scale=2,  >=stealth]

    % 绘制正方形格子
    \draw[thick] (0,0) -- (1,0) -- (1,1) -- (0,1) -- cycle;

    % 绘制蓝色粗箭头（更粗的线条）
    \draw[->, line width=2pt, blue] (0,0) -- (1,0); % 从 (0,0) 到 (1,0)
    \draw[->, line width=2pt, blue] (1,0) -- (1,1); % 从 (1,0) 到 (1,1)

    % 标注点（经过的点变为蓝色，坐标标签加粗并变为蓝色）
    \filldraw[blue] (0,0) circle (1.5pt) node[below left, font=\bfseries, blue] {$(0,0)$};
    \filldraw[blue] (1,0) circle (1.5pt) node[below right, font=\bfseries, blue] {$(1,0)$};
    \filldraw[blue] (1,1) circle (1.5pt) node[above right, font=\bfseries, blue] {$(1,1)$};
    \filldraw[black] (0,1) circle (1pt) node[above left] {$(0,1)$}; % 未经过的点保持黑色

    % 标注维度
    \node at (0.5, -0.2) {$x_i^{t-1}$}; % x_i^{t-1} 维度
    \node at (-0.2, 0.5) {$x_i^{t-2}$}; % x_i^{t-2} 维度

\end{tikzpicture}}
  \end{subfigure} 
\hspace{0.6cm}
  \begin{subfigure}[b]{0.4\textwidth}
  \centering 
   \caption{\re{Mixed effects: $\beta_i^1 > 0, \beta_i^2\le 0 $} }
    \label{fig:switch2} 
    \resizebox{0.5\textwidth}{!}{\begin{tikzpicture}[scale=2, >=stealth]

    % 绘制正方形格子
    \draw[thick] (0,0) -- (1,0) -- (1,1) -- (0,1) -- cycle;

    % 绘制蓝色粗箭头（更粗的线条）
    \draw[->, line width=2pt, blue] (1,0) -- (0,0); % 从 (0,0) 到 (1,0)
    \draw[->, line width=2pt, blue] (0,0) -- (0,1); % 从 (1,0) 到 (1,1)

    % 标注点（经过的点变为蓝色，坐标标签加粗并变为蓝色）
    \filldraw[blue] (0,0) circle (1.5pt) node[below left, font=\bfseries, blue] {$(0,0)$};
    \filldraw[blue] (1,0) circle (1.5pt) node[below right, font=\bfseries, blue] {$(1,0)$};
    \filldraw[blue] (0,1) circle (1.5pt) node[above left, font=\bfseries, blue] {$(0,1)$};
    \filldraw[black] (1,1) circle (1pt) node[above right] {$(0,1)$}; % 未经过的点保持黑色

    % 标注维度
    \node at (0.5, -0.2) {$1-x_i^{t-1}$}; % x_i^{t-1} 维度
    \node at (-0.2, 0.5) {$x_i^{t-2}$}; % x_i^{t-2} 维度

\end{tikzpicture}} 
  \end{subfigure}
\end{figure}
\re{Now, we formally define the historical assortments used in our construction. For each product $i \in [N]$, let $I_i$ denote the indices of positive effects, that is, $I_i := \{m \mid \beta_i^m>0\}$. Let $\Omega$ denote all permutations of past periods $[M]$. For each product $i \in [N]$ and a permutation $\sigma \in \Omega$, we define $M+1$ historical assortments, denoted as $\h^\sigma_{i,0}, \h^\sigma_{i,1}, \ldots, \h^\sigma_{i,M} $, as follows: $\h^\sigma_{i,0} = \chi^{I_i}$ and for each $k  \in [M]$}
\begin{equation}\label{eq:historyassortment}
\re{\h^\sigma_{i,k} = \begin{cases}
    \h^\sigma_{i,k-1} + \e_{\sigma(k)} & \text{ if  }  \sigma(k) \notin  I_i\\
     \h^\sigma_{i,k-1} -  \e_{\sigma(k)} & \text{ if } \sigma(k) \in I_i,\tag{\textsc{\re{SwitchNested}}}
\end{cases}}
\end{equation}
\re{With this definition, we are ready to describe the concave envelope of $\alpha_i(\cdot)$ which, after perspectification, yields an explicit linear description for the constraint~\eqref{eq:conc-pers}.}

\begin{proposition}\label{prop:conc-alpha}
\re{Constraint~\eqref{eq:conc-pers} is equivalent to the following system of linear inequalities.}
\re{\begin{equation}\label{eq:pers-extension}
     y^t_i \leq \alpha_i(\h^\sigma_{i,0} ) ( \gamma^t_i - \tilde{z}_{i\sigma(1)}^t) + \sum_{k \in [M]} \alpha_i(\h^\sigma_{i,k} ) (   \tilde{z}_{i\sigma(k)}^t - \tilde{z}_{i\sigma(k+1)}^t)   \qquad  \text{ for } \sigma \in \Omega\\
\end{equation}where $\tilde{z}_{i\sigma(M+1)}^t = 0$, $\tilde{z}^t_{i\sigma(k)} =    z^t_{i\sigma(k)}$ if $\sigma(k) \notin  I_i$, and $\tilde{z}^t_{i\sigma(k)} =  \gamma_i^t -  z^t_{i\sigma(k)}$ if $ \sigma(k) \in I_i$. } 
\end{proposition}

\subsubsection{\re{The Final Formulation}}
Using the linear constraints~\eqref{eq:relax1},~\eqref{eq:relax2}, and~\eqref{eq:pers-extension} and the exponential cone constraint~\eqref{eq:conv-pers} to represent~\eqref{eq:choice}, we obtain an MIECP formulation of problem~\eqref{eq:total} as follows:     
\begin{equation}\label{eq:Conic}
\begin{aligned}
    \max_{\x,\y,\boldsymbol{\rho},\boldsymbol{\gamma},\z} \quad &  \frac{1}{T} \sumt \sumi  r_iy_{i}^t  \\
    \st \quad &   \x^t \in \{0,1\}^{N } \cap \mathcal{X} \text{ and } \x \in \mathcal{P} && \for \, t \in [T]  \\
    &\rho^t + \sumi y_{i}^t = 1  \quad && \for\, t \in [T]  \\
    &\eqref{eq:relax1},~\eqref{eq:relax2},~\eqref{eq:conv-pers},~\eqref{eq:pers-extension} \quad && \for \, t \in [T] \text{ and } i \in [N].
\end{aligned}     \tag{\textsc{Conic}}
\end{equation}

\begin{theorem}\label{them:decomposition}
\eqref{eq:Conic} is a mixed-integer exponential cone formulation of problem~\eqref{eq:total}.  
\end{theorem}
% \begin{remark}
%     Note that the number of inequalities in 
% \end{remark}
\begin{remark}\label{rmk:M=2}
   In this remark, we consider the case where the memory length is less than two. In this case, we can derive the tightest convex extension, also known as the \textit{convex envelope}, of $\alpha_i(\cdot)$. Using the convex envelope instead of the natural continuous extension, we obtain the following tighter formulation than the one given by constraint \eqref{eq:conv-pers} for $M=2$. \re{If $\boldb_i$ are all positive or negative, we have} 
\begin{subequations}\label{eq:conv-pers-linear1}
    \begin{alignat}{3}
         y_i^t\ge  &   \alpha_i(\boldsymbol{0}) (\gamma_i^t-z_{i1}^t-z_{i2}^t) +  \alpha_i(\boldsymbol{e}_1) z_{i1}^t+  \alpha_i(\boldsymbol{e}_2)  z_{i2}^t     \\
y_i^t  \ge &  \alpha_i(\boldsymbol{e}_1)(\gamma_i^t - z_{i2}^t) + \alpha_i(\boldsymbol{e}_2)(\gamma_i^t - z_{i1}^t) +  \alpha_i(\boldsymbol{1})  (z_{i1}^t +z_{i2}^t -  \gamma_i^t),
    \end{alignat}
\end{subequations}
\re{ and if $\boldb_i$ are mixed, the convex envelope is}
 \re{\begin{subequations}\label{eq:conv-pers-linear2}
    \begin{alignat}{3}
y_i^t \ge & \alpha_i(\boldsymbol{0}) (\gamma_i^t - z_{i1}^t)  + \alpha_i(\boldsymbol{e}_1) (z_{i1}^t  -z_{i2}^t ) + \alpha_i(\boldsymbol{1})  z_{i2}^t  \\
y_i^t \ge & \alpha_i(\boldsymbol{0}) (\gamma_i^t - z_{i2}^t)  + \alpha_i(\boldsymbol{e}_2) (z_{i2}^t  -z_{i1}^t )+ \alpha_i(\boldsymbol{1})  z_{i1}^t  
    \end{alignat}
\end{subequations}}where $\boldsymbol{e}_m\in\{0,1\}^2$ is the unit vector with one as the $m^{\text{th}}$ coordinate and zero everywhere else, and $\boldsymbol{1}=(1,1)$. For $M = 1$, constraints~\eqref{eq:conv-pers-linear1} and~\eqref{eq:conv-pers-linear2}  reduce to $y_i^t \ge  \bigl(  \alpha_i(1) - \alpha_i(0)\bigr) z_{i1}^t + \alpha_i(0)\gamma_i^t$. 
In addition to the tightness, another advantage of constraint~\eqref{eq:conv-pers-linear1} and~\eqref{eq:conv-pers-linear2} over the conic constraint~\eqref{eq:conv-pers} is that they are linear inequalities. Therefore, replacing~\eqref{eq:conv-pers} in~\eqref{eq:Conic} with~\eqref{eq:conv-pers-linear1} and~\eqref{eq:conv-pers-linear2}, we obtain an MILP formulation of \eqref{eq:total}. This allows us to solve~\eqref{eq:total} using a commercial MILP solver---such as \texttt{GUROBI}---which has advantages over its conic counterparts in terms of practical scalability and numerical stability, as illustrated in Section~\ref{sec:numerical:speed}. \hfill \Halmos
\end{remark}

We conclude this section with a discussion of how our formulation relates to existing formulations. Recently, a large amount of research has focused on developing strong formulations for binary linear fractional programming problems and MNL-based assortment planning problems~\citep[e.g.][]{sen2018conic,mehmanchi2019fractional,atamturk2020submodularity,kilincc2023conic,he2024convexification, chen2024integer}. To apply these results in solving problem \eqref{eq:total}, we need to linearize the nonlinear functions appearing in denominators and numerators---specifically, the attraction value functions. A prevalent linearization technique replaces each attraction value function by its multilinear extension, called the Fourier Expansion in~\cite{o2014analysis}, and then linearizes the resulting multilinear functions using the recursive McCormick relaxation~\citep{mccormick1976computability,khajavirad2023strength}.
We combine this linearization strategy with binary linear fractional programming formulation techniques to obtain an exact formulation of problem~\eqref{eq:total}, which we refer to as the multilinear-extension-based formulation. The complete derivation is provided in the ecompanion~\ref{sub:multi:detail}. Later, in Section~\ref{sec:numerical:speed}, a computational experiment shows that our formulation \eqref{eq:Conic} significantly outperforms the multilinear-extension-based formulation.

Our problem setting is related to the logit-based share-of-choice product design (SOCPD) studied in~\cite{akccakucs2023exact}, which aims to maximize the share of $K$ customer types by selecting $n$ binary design attributes for a single product, that is, 
\[
\max_{\a \in \mathcal{A} \subseteq \{0,1\}^n } \quad \sum_{k=1}^K\lambda_k \cdot \frac{\exp(u_k(\a))}{1 + \exp(u_k(\a))},
\]
where for each customer type $k$, $u_k(\a)$ is an affine utility function for the product, the fraction term models the purchase probability of a customer of type $k$, and $\lambda_k$ is the fraction of customers who belong to type $k$. Our model shares some similarities with SOCPD since the exponential-affine attraction value functions appear in both problems. However, there are two major differences. First, since we consider assortments of multiple products, our formulation~\eqref{eq:Conic} requires both convex and concave extensions of the attraction value function to achieve exactness. In contrast, since SOCPD treats only one single product, it suffices to use the natural convex extension of the exponential-affine function in that case. Second, in our model, the attraction value function is multiplied by an assortment decision variable to model whether the attractiveness of a product is activated in the current period. This requires us to adapt the perspectification trick used in~\cite{akccakucs2023exact} to allow for nonlinear scaling, as discussed in Section~\ref{section:envelope-characterization}.

\subsection{\re{Optimality of Revenue-Ordered Policies}}\label{sec:model:RO}
% \re{In this subsection, we identify a special case in which~\eqref{eq:total} can be solved in polynomial time and build connections with the classic revenue-ordered assortment structure in literature.} 
 
% \re{Theorem~\ref{thm:seq:rev:opt}, shows that if the history-dependent effects are non-negative and the cross-product and cross-period constraints are absent then a \textit{sequential revenue-ordered} policy is optimal. Such a RO policy is obtained by sequentially maximizing the revenue in each period with the product utility updated using past assortments. Details of the sequential RO policy is given in Algorithm~\ref{alg:greedy}.}

% are a family of nested assortments selected according to their revenue-decreasing order. 

A prevalent strategy for solving assortment optimization is to construct Revenue-Ordered (RO) assortments, which is an assortment including products with revenue higher than a threshold. Specifically, RO structures and their variants have been shown to be optimal for solving assortment optimization problems under the MNL choice model~\citep{talluri2004revenue}, sequential MNL choice models \citep{gao2021assortment}, MNL choice models with network effects~\citep{wang2017consumer}, etc. Inspired by the simplicity and optimality of RO structures, we would like to understand how it works in our problem.

We find that when there are no cross-product and cross-period constraints, if the history-dependent effect is non-negative, that is, $\boldb \ge 0$, \re{the optimal assortment in each period is an RO assortment, and the revenue threshold increases when $T$ increases.} Such an assortment planning can be obtained by a \textit{sequential revenue-ordered} policy. The sequential RO policy is obtained by sequentially selecting the best RO assortment based on the currently updated product utility and maximizing revenue in the current period. This policy is adaptive while ignoring the impact of the current assortment decision on future revenue. Details of the sequential RO policy are shown in Algorithm~\ref{alg:greedy}.  Because there are no cross-product and cross-period constraints, the computation time required to find the optimal RO assortment is $\mathcal{O}(N)$ in each period \citep{talluri2004revenue,liu2008choice}. Thus, the total computation time of the sequential RO policy is $\mathcal{O}(NT)$. 

\begin{algorithm} 
    \caption{Sequential revenue-ordered policy  of~\eqref{eq:total}}\label{alg:greedy}
    \begin{algorithmic}[1]
        \STATE Input:  $\boldsymbol{r},\{\beta_i^0\}_{i\in[N]}, \boldb,M,N,T$  
        \FOR {$t = 1, 2,\dots,  T $} 
        \STATE Update history assortments: $\x^{t-1},\dots,\x^{t-M}$, and set $\x^{t-m} = \boldsymbol{0}$ if $t-m\le 0$.  
        \STATE Compute the best revenue-ordered assortment $\x^t = \argmax \limits_{\boldsymbol{z}\in \{0,1\}^N}  \sumi r_i \pi^t_i(\boldsymbol{z}, \x^{t-1}, \ldots, \x^{t-M})$ 
        \ENDFOR 
        % \STATE Set $\x^k = \x^{\min\{T,M+1\}}$ for $k = \min\{T,M+1\}+1, \ldots,T$
        \STATE Output: $(\x^1,\dots, \x^T)$
    \end{algorithmic}
\end{algorithm}

We are now in a position to state the main result of this subsection, Theorem~\ref{thm:seq:rev:opt}, characterizing conditions on which the sequential RO policy is optimal. 

\begin{theorem}\label{thm:seq:rev:opt}
 Assume the absence of cross-product and cross-period constraints. Then, the sequential-revenue-ordered policy solves~\eqref{eq:total} if the history-dependent effects are non-negative, that is, $\boldb\ge \boldsymbol{0}$.
\end{theorem}

% \re{Theorem~\ref{thm:seq:rev:opt} identifies optimality condition when the sequential RO policy is optimal. Although this policy ignores the impact of the current assortment decision on future revenue, it is optimal when the history-dependent effects are non-negative. Its computation time is only $\mathcal{O}(NT)$ since the computation time required to fined an optimal RO assortment in each period is $\mathcal{O}(N)$ \citep{talluri2004revenue}.} \re{Theorem~\ref{thm:seq:rev:opt} also characterizes the optimal assortment structure when $\boldb\ge \boldsymbol{0}$. Its proof indicates that the assortments are RO in each period and nested with decreasing size over the time horizon.}

% \re{Theorem~\ref{thm:seq:rev:opt} contributes to the literature that constructs RO assortments for solving unconstrained assortment optimization problem. Studies have shown that RO structures and their variants are optimal under the MNL model~\citep{talluri2004revenue}, sequential MNL models \citep{gao2021assortment}, and MNL models with network effects~\citep{wang2017consumer}. Theorem~\ref{thm:seq:rev:opt} indicates that the RO structure still exists under non-negative history-dependent effects. Moreover, Proposition~\ref{prop:hold:condition} shows that for an optimal assortment planning of \eqref{eq:total}, the set of products that are offered across all periods has a RO structure.}

For the case of negative history-dependent effects, which is the computationally intractable regime of our model, as we argued in Proposition~\ref{prop:nphard}, the sequential policy fails to solve our model. However, the RO structure still exists under the negative and even mixed history-dependent effects and appears in a different form. For an optimal assortment planning of model \eqref{eq:total}, we argue that the set of products that are offered across all periods has an RO structure.

\begin{proposition}\label{prop:hold:condition} Assume the absence of cross-product and cross-period constraints. Then, there exists an optimal assortment $\x = (\x^1, \ldots, \x^T)$ of~\eqref{eq:total} such that $\bigl\{i \in [N] \bigm| \sum_{t \in T}x^t_i \ge 1 \bigr\}$ is revenue-ordered.   
\end{proposition}

% We conclude this section with an example to numerically illustrate that the optimality gap of RO-type policies can be significant. This motivates us to study the global optimality of \eqref{eq:total}
% in Section~\ref{sec:reformulation}. 

\section{Cyclic Policies: Optimality and Computation}\label{sec:two:cycle}

In this section, we consider the situation in which a firm plans assortments for a long period. We establish the asymptotic optimality of cyclic policies in Section~\ref{sec:general:cycle} and formulate the cyclic policy problem as an MIECP in Section~\ref{sec:formulation:cycle}. In Section~\ref{sec:cyclength}, we characterize conditions under which the length of optimal cyclic policies is the memory length plus one. \re{In Section~\ref{sec:bf:model}, we show that such conditions lead to a tighter MILP formulation for finding optimal $(M+1)$-cyclic policies.}

 \subsection{Optimality of Cyclic Policies}\label{sec:general:cycle} 
The assortment planning problem with an infinite horizon can be formulated as follows. The firm can select any assortment sequence (or policy) $\x =\{\x^t\}_{t \in \Z_+}$ with elements in $\{0,1\}^N$. As before, the firm may face cross-product and cross-period constraints. Let $\mathcal{X}$ (resp. $\mathcal{P}^{\infty}$) denote the feasible set of cross-product (resp. cross-period) constraints. Then, the firm aims to find an optimal assortment sequence for the following problem:
\begin{equation}\label{eq:asymptotic}
  \sup_{\x}\biggl\{ \lim_{T\to \infty}  \frac{1}{T}  \sum_{t \in [T]}  \sum_{i \in [N]} r_i\pi_i^t(\x^t,\dots,\x^{t-M}) \biggm| \x^t \in \mathcal{X} \cap \{0,1\}^N \; \for t \in \Z_+,\ \x \in \mathcal{P}^\infty \biggr\}.\tag{\textsc{Infty}}
\end{equation}
An assortment policy $\x=\{\x^t\}_{t \in \Z_+}$ is \textit{cyclic} if there exists a positive integer $L$ such that $\x^{t+L} = \x^t$ for all $t \in \Z_+$. The smallest $L$ for which this holds is called the cycle length of $\x$. A cyclic policy $\x$ can be represented by a finite sequence of assortments $(\x^1, \ldots, \x^L )$, where $L$ is the cycle length of $\x$. In contrast, whenever we say that a finite sequence of assortments $(\x^1, \ldots, \x^L)$ is cyclic, we are referring to the policy in which this finite sequence of assortments is repeated infinitely often. 

For a given cycle of length $L$, we are interested in finding a cyclic policy of length $L$ that maximizes the long-term average revenue. This is equivalent to searching for a sequence of assortments $(\x^1, \ldots, \x^L )$ of length $L$ that maximizes the average revenue by repeating the sequence. Such a sequence can be defined as an optimal solution of a variant of~\eqref{eq:total}. To define the variant, we modify the attraction value function so that it depends on past assortments in the cycle. Given a position $t$ in the cycle of length $L$, we use a function $\tau(\cdot|t)$ to track positions within the memory length $M$:
\begin{align*}
\tau(m|t) = \begin{cases}
     {(t-m \text{ mod } L)}  & \text{ if } (t-m \text{ mod } L) > 0  \\
     {L} & \text{ if } (t-m \text{ mod } L) =0.
\end{cases}   
\end{align*}
The attraction value of product $i$ at position $t$ can be specified as follows: 
\[
\alpha_i^t(x^1_i, \ldots, x^L_i ) = \exp\biggl( \beta_i^0 + \sum_{m \in [M]} \beta^m_ix^{\tau(m \mid t) }_i\biggr). 
\]
A cyclic policy is called an \textit{$L$-cyclic policy} if it is optimal for the following variant of problem~\eqref{eq:total}:
\begin{equation}\label{eq:cyclic-policy}
   %\begin{aligned}
       \max \Biggl\{ \frac{1}{L} \sum_{t \in [L]} \sum_{i \in [N]} \frac{r_ix^t_i\alpha^t_i(x^1_i, \ldots, x^L_i)}{1+ \sum_{j \in [N]}x^t_j\alpha^t_j(x^1_i, \ldots, x^L_i)} \Biggm|
        \x^t \in \mathcal{X} \cap \{0,1\}^N \for t \in [L] ,\ \x \in \mathcal{P}^L \Biggr\},\tag{\textsc{Cycle}}
   %\end{aligned} 
\end{equation}
where $\mathcal{X}$ denotes cross-product constraints and $\mathcal{P}^L$ denotes cross-period constraints. In addition, we assume that a finite sequence of assortments satisfying $\mathcal{P}^L$ can be repeated infinitely often to obtain a cyclic policy that satisfies the cross-period constraint $\mathcal{P}^\infty$ in~\eqref{eq:asymptotic}. 

The main result of this subsection relates to problems~\eqref{eq:asymptotic} and~\eqref{eq:cyclic-policy}. To establish this result, we represent~\eqref{eq:asymptotic} as a directed graph, which we refer to as the \textit{assortment graph}. Each node of the graph denotes a list of $M$ assortments, indexed from the latest to the earliest assortments in memory. Each arc represents a memory transition driven by an assortment decision. Specifically, an arc starting at a node should end at a node that drops the earliest assortment of the starting node and adds a new assortment. Hence, each arc encodes an assortment decision in a period, and its weight is defined as the single-period revenue generated by the new assortment. With such a graph, we can show that~\eqref{eq:asymptotic} is equivalent to finding the maximum mean cycle of the assortment graph \citep{karp1978characterization}, whose cycle length is denoted as $L^*$. Given $L^*$, we can use \eqref{eq:cyclic-policy} to compute the cyclic policy. We formalize this observation as the following theorem.

\begin{theorem}\label{thm:general:cyclic}
Given an instance of~\eqref{eq:asymptotic}, there exists a positive integer $L^*$ such that $L^*$-cyclic policy is optimal to~\eqref{eq:asymptotic}.
\end{theorem} 
% A detailed proof is provided in the ecompanion. Formulation techniques developed for solving~\eqref{eq:total} can be used to obtain an MIECP formulation for~\eqref{eq:cyclic-policy}; see details in Section~\ref{sec:formulation:cycle}

\subsection{Finding Optimal Cyclic Policies} \label{sec:formulation:cycle}

Formulation techniques developed for solving~\eqref{eq:total} can be used to obtain a mixed-integer exponential cone formulation for~\eqref{eq:cyclic-policy}. Since the main discrepancy between the two models is the definition of attraction value functions, we only need to modify constraints in~\eqref{eq:Conic} that involve attraction value functions and historical assortments. In particular, we replace $x^{t-m}_i$ in~\eqref{eq:relax2} by $x_i^{\tau(m|t)}$ and $\alpha_i(\cdot)$ in~\eqref{eq:pers-extension} by $\alpha^t_i(\cdot)$:
\re{\begin{subequations}\label{eq:conv-pers-modify}
    \begin{alignat}{3}
        & \max \bigl\{ 0,\ \rho^t_U x^{\tau(m|t)}_i + \gamma^t_i - \rho^t_U \bigr\}  \leq z_{im}^t   \leq \min \bigl\{ \gamma_{i}^t , \rho^t_U x^{\tau(m|t)}_i  \bigr\} &&\quad \for \, m \in [M]  \label{eq:modify:z}\\
       & y^t_i  \leq \alpha_i^t(\w^\sigma_{i,0} ) ( \gamma^t_i - \tilde{z}_{i\sigma(1)}^{t}) + \sum_{k \in [M]} \alpha_i^t(\w^\sigma_{i,k} ) (   \tilde{z}_{i\sigma(k)}^{t} - \tilde{z}_{i\sigma(k+1)}^{t}) &&\quad  \for\, \sigma \in \Omega, \label{eq:modify:h}
    \end{alignat}
\end{subequations}}\re{where $\w^{\sigma}_{i,k}\in \{0,1\}^L$, and its $\tau(m|t)^{\text{th}}$ element equals the $m^{\text{th}}$ element of $\h_{i,k}^{\sigma}$ defined in~\eqref{eq:pers-extension} and others are zero.} This replacement yields the following mixed-integer exponential cone formulation of~\eqref{eq:cyclic-policy}:
\begin{equation}\label{eq:cyclic-policy-conic}
\begin{aligned}
\max_{\x,\y,\boldsymbol{\rho},\boldsymbol{\gamma},\boldsymbol{\z}} \quad & \frac{1}{L}  \sum_{t \in [L]}  \sumi r_iy_{i}^t  \\
    \st \quad & \x^t \in \mathcal{X} \cap \{0,1\}^{N }  \text{ and } \x \in \mathcal{P}^{L}   && \for  t \in [L] \\
    &\rho^t + \sumi y_{i}^t = 1   && \for\, t \in [L]  \\
    & ~\eqref{eq:relax1}~\eqref{eq:conv-pers}~\eqref{eq:conv-pers-modify} && \for \, t \in [L] \text{ and } i \in [N]. \\
\end{aligned}     \tag{\textsc{Cycle-Conic}}
\end{equation}\eqref{eq:cyclic-policy-conic} has its own value because it is an efficient way to generate a feasible assortment planning with a large planning horizon. Especially, the retailer may pre-determine a cycle length $L$ based on experience or industry knowledge and use \eqref{eq:cyclic-policy-conic} to derive the optimal policy under a given $L$. %\re{When all $\beta_i$ are positive and no cross-period constraints, the optimal cycle length is 1. When some $\beta_i$ are negative, it is not straightforward to determine the optimal cycle length. In the next subsection, we characterize the optimal cycle length under a common non-overlapping assortment assumption.}
 
\subsection{Cycle Length Characterization}\label{sec:cyclength}
It is intuitive that, in order to avoid the negative impact on a product's utility if the product exhibits a satiation effect, a product should not be offered consecutively. In particular, we say that a cross-period constraint on a policy $\x$ is \textit{non-overlapping} if $M+1$ adjacent assortments $(\x^t, \x^{t+1},\dots,\x^{t+M})$ are non-overlapping, that is, for each product $i \in [N]$,  $x^t_i + x^{t+1}_i +\dots +x^{t+M}_i\leq 1$ for every $t \in \Z_+$. This non-overlapping condition is common in assortment optimization literature~\citep{liu2020assortment,chen2023assortment}. \re{For instance, \cite{li2024should} identify a similar non-overlapping condition on customer orders to plan assortments for a local warehouse.} Under the non-overlapping condition, any cyclic policy with length $L\le M $ must have empty sets in some periods. We do not consider such cyclic policies. For cyclic policies with length $L>M$, we find that this condition allows us to relate the length of an optimal cyclic policy for solving~\eqref{eq:asymptotic} to the customers' memory length as follows. 
\begin{theorem}\label{thm:nonoverlap:condition}
Assume that $M < N$. An $(M+1)$-cyclic policy is optimal for problem~\eqref{eq:asymptotic} if the cross-period constraint is non-overlapping. \end{theorem} 
%\begin{remark}\label{rmk:cycle:general}
\re{We remark that Theorem~\ref{thm:nonoverlap:condition} holds for mixed satiation-addiction effects, because its proof replies on a graph whose structure unaffected by signs of $\boldb$.}
%\end{remark}

The condition $M<N$ guarantees that there at least exists a non-overlapping assortment policy without empty assortments. Without the non-overlapping condition, the optimal cycle length varies as the parameters change. In particular, a slight fluctuation in the input parameters, like the base utility of a product, could change the optimal cycle length. We provide an example to illustrate this situation in~\ref{eg:cycle:change}. 

The base utility and history-dependent effect are estimated from real data and lie in confidence intervals. Hence, given the fact that the optimal cycle length is sensitive to input parameters, finding a cycle length with robust performance may be more critical than characterizing the optimal cycle length in practice. A positive result is that although the $(M+1)$-cyclic policy is not necessarily optimal for~\eqref{eq:asymptotic} without the non-overlapping constraint, it still provides a good approximation, as illustrated in the following example. 
 
% Theorem~\ref{thm:general:cyclic} indicates that for a given instance of~\eqref{eq:asymptotic}, once we know the optimal cycle length $L^*$, it suffices to solve~\eqref{eq:cyclic-policy} with the optimal cycle length $L^*$. The infinite horizon problem~\eqref{eq:asymptotic} is decomposed to a finite one and becomes solvable. Furthermore, if $L^*$ is relatively small, we can efficiently solve a small-scale~\eqref{eq:cyclic-policy} and generate the $L^*$-cyclic policy. Hence, a natural question is: what is the relationship between input parameters $M$, $\boldb$, and $\boldsymbol{r}$, and the optimal cycle length? If it is hard to characterize $L^*$, is there any good approximation of $L^*$? Next, we will show although we can not characterize $L^*$ in all cases, we find $L^* = M+1$ in a special case. For more general cases, we show $(M+1)$-cyclic policy still has performs well by a numerical study.

% Example~\ref{eg:hard:character} shows the hardness to characterize $L^*$ in general cases. The good news is that the following example indicates $M+1$ is a good approximation of $L^*$ because $(M+1)$-cyclic policy has good numerical performance in various cases.
\begin{example}
We numerically test the performance of the $(M+1)$-cyclic policy {under satiation effects}. We fix the number of products $N=30$ and vary the planning horizon $T\in\{10,30,50\}$. The revenue and the base utility of each product are randomly generated from uniform distributions with ranges $[1,10]$ and $[-1,1]$, respectively.  For a weak (W) (resp. strong (S)) satiation effect, we set the history-dependent effect $\beta_i^m$ using a uniform distribution $U[-1,0]$ (resp. $U[-2,-1]$). For each instance, we compute the revenue gap between the $(M+1)$-cyclic policy and the best feasible solution, defined as $100\% \times \frac{R_{\mathsf{fea}} - R_{\mathsf{cycle}}}{R_{\mathsf{fea}}}$, where $R_{\mathsf{cycle}}$ is the revenue achieved by the $(M+1)$-cyclic policy and $R_{\mathsf{fea}}$ is the revenue of the best feasible solution obtained by~\eqref{eq:Conic} within $7200$ seconds. 
\begin{table} 
\centering
\caption{\re{Revenue gap of $(M+1)$-cyclic policies under satiation effects with $N=30$ and no constraints}}\label{tab:cyclic}
\begin{tabular}{ccccccccc}
\cline{1-4}\cline{6-9}
$\boldb$ & \textbf{T} & $ M=1$  ($\%$) & $ M=2$ ($\%$) &  \text{   }   &  $\boldb$ & \textbf{T} & $ M=1$  ($\%$) & $ M=2$ ($\%$) \\ \cline{1-4}\cline{6-9}
W &10   & 0.51  & 1.69  & \text{   } &         S &10   & 0.39  & 4.45            \\
W &30    & 0.16     & 0.60  & \text{   } &           S &30   & -0.06  & -0.50          \\ 
W &50    & 0.08     & 0.52   &\text{   } &          S &50   & -0.52  & 0.46         \\  \cline{1-4}\cline{6-9}
\end{tabular}
\end{table}
Table \ref{tab:cyclic} records the average gap of five randomly generated instances of each parameter configuration. A negative gap indicates that the revenue generated by the $(M+1)$-cyclic policy is higher than the revenue generated by the best feasible solution achieved within 7200s. Overall, the $(M+1)$-cyclic policy performs well. Moreover, as the length of the planning horizon increases, the revenue gap decreases. \hfill \Halmos\end{example}
% The table shows that $(M+1)$-cyclic policy performs well for all cases, especially under a large planning horizon. For instance, the average gap is less than $0.6\%$ when $T=30$. The revenue gap could be negative because the feasible solution obtained within $2$ hours is even smaller than that of the $(M+1)$-cyclic policy. It indicates that computing the $(M+1)$-cyclic policy requires little computation time when $M$ is small.   \Halmos

% The non-overlapping condition provides an additional benefit in problem formulation when cycle length is small. Next, we provide a  model formulation when $M=1$ with non-overlapping requirement. 

\subsection{\re{A Strong Formulation for Finding Non-overlapping $(M+1)$-Cyclic Policies}}\label{sec:bf:model}
\re{ In this subsection, we propose a bound-free formulation for the $(M+1)$-cyclic policy, which is provably tighter than its conic counterpart~\eqref{eq:cyclic-policy-conic} whose computational performance depends on the bounds ($\rho^t_L$ and $\rho^t_U$) of the no-purchase probability variable $\rho^t$. } 
%, and, in Section~\ref{sec:numerical:cycle}, we demonstrate its computational efficacy

%\re{We start by reformulating~\eqref{eq:cyclic-policy} into a simpler form. 
%which is a special case of the assortment optimization problem studied in~\cite{caro2014assortment}.
%}
% When $M = 1$, attraction value functions are univariate functions over the binary domain $\{0,1\}$, and, thus, can be represented using an affine function. For each product $i \in [N]$, let  $u_i := \exp(\beta^0_i) $ and $d_i := \exp(\beta^0_{i} + \beta^1_{i}) - \exp(\beta^0_{i})$. Then, for each period $t\in[T]$, the attraction value function can be linearized as follows: 
% \[
% \alpha^1_i(x^2_i) = u_i + d_ix^{2}_i \quad \text{ and } \quad \alpha^2_i(x^1_i ) = u_i + d_ix^{1}_i   \quad \for \, (x^1_i,x^2_i) \in \{0,1\}^2,
% \]
% where $u_i$ represents the base attraction value of product $i$, and $d_i$ is the incremental attraction value of product $i$ conditional on whether it is offered in another period.
\re{From the non-overlapping constraint $x_i^t +  x^{\tau(1 \mid t)}_i + \dots +  x^{\tau(M\mid t)}_i \leq 1$, we have:}
\[
\re{x_i^t \alpha_i^t(x^1_i, \ldots, x^L_i ) = x_i^t \exp\biggl( \beta_i^0 + \sum_{m \in [M]} \beta^m_ix^{\tau(m \mid t) }_i\biggr) = x_i^t \exp ( \beta_i^0) \qquad \text{ for } t \in [M+1].}
\]
\re{Let $u_i = \exp(\beta_i^0)$. ~\eqref{eq:cyclic-policy} is equivalent to }
\re{\begin{equation}\label{eq:2cycMNL}
\begin{aligned}
    \max \quad & \frac{1}{M+1}  \sum_{t \in [M+1]}  \sumi \frac{r_i u_ix^t_i }{1+ \sumj u_jx^t_j} \\
    \st \quad  & x_i^t + \summ x^{\tau(m \mid t)}_i  \leq 1 &&   \for i \in [N] \\
    &\x^t \in \mathcal{X} \cap \{0,1\}^{N} &&   \for t \in [M+1].
\end{aligned}\tag{\textsc{(M+1)-Cyclic}}
\end{equation}}Note that this reformulation reveals a connection between~\eqref{eq:cyclic-policy} and the assortment problem studied in~\cite{caro2014assortment} in which a retailer must decide, in advance, the release date of each product in a given collection over a selling season. More specifically, \re{when $M=1$}, \eqref{eq:2cycMNL} is exactly the case where a retailer faces two planning horizons and sells products with a single period life-cycle, and, in each period, customers make a purchase decision on offered products according to the MNL model. 

Next, we introduce variables for building a strong formulation of~\eqref{eq:2cycMNL}. For each period $t \in [M+1]$, we use $\rho^t$ and $u_i\cdot \gamma^t_i$ to represent the no-purchase probability and choice probability of product~$i$, respectively,
\begin{equation}\label{eq:CC}
\rho^t = \frac{1}{1+\sum_{j\in [N]}u_jx^t_j} \quad \text{ and } \quad \gamma_i^t = \frac{x^t_i}{1+\sum_{j\in [N]}u_jx^t_j} \qquad \text{ for } i \in [N].  \tag{\textsc{CC}}
\end{equation}
Note that this is known as the Charnes-Cooper transformation~\citep{charnes1962programming}. For each period \re{$t \in [M+1]$}, it maps all possible assortments $\x^t \in \{0,1\}^N$ into all possible choice probabilities $(\rho^t,u_1\gamma^t_1, \ldots,  u_N\gamma^t_N )$, where $(\rho^t, \boldsymbol{\gamma}^t)$ belongs to 
\begin{equation}\label{eq:choice-prob-two-cycle}
    \rho^t + \sum_{j \in [N]}u_j\gamma^t_j = 1 \quad \text{ and } \quad  \gamma^t_j  \in \{0, \rho^t\} \, \for j \in [N].
\end{equation}
Using these relations, we are ready to present a new formulation of~\eqref{eq:2cycMNL}:
\re{\begin{subequations}\label{eq:2cycMILP-boundfree}
\begin{alignat}{3}
    \max \limits_{\x, \boldsymbol{\rho},\boldsymbol{\gamma}, \boldsymbol{\Gamma}}\quad & \frac{1}{M+1} \sum_{t \in [M+1]}  \sumi  r_i  u_i \gamma_{i}^t \notag \\ 
 \st \quad & x_i^t + \summ x^{\tau(m \mid t)}_i  \leq 1  \text{ and }   \x^t \in \mathcal{X} \cap\{0, 1\}^N && \for i \in [N] \text{ and } t\in [M+1] \\
        &  \rho^t + \sum_{i \in [N]}u_i \gamma^t_i = 1 \text{ and } 0 \leq  \gamma^t_i \leq \rho^t \quad  && \for  i \in [N] \text{ and } t\in [M+1]\label{eq:perfect-1} \\
 &x^t_i = \gamma^t_i + \sum_{j \in [N]}u_j \Gamma^t_{ij} && \for i \in [N] \text{ and }  t\in [M+1] \label{eq:perfect-2}  \\
        & \max\{0, \gamma^t_i + \gamma^t_j -  \rho^t\} \leq \Gamma^t_{ij}  \leq \min \{\gamma^t_i, \gamma^t_j\} \quad && \for  i \in [N] \, ,  j \in [N] \setminus\{i\}, \text{ and } t\in [M+1] \label{eq:perfect-3}\\
        & \Gamma^t_{ii} = \gamma^t_i  && \for i \in [N] \text{ and }  t\in [M+1].\label{eq:perfect-4}
\end{alignat}
\end{subequations}}The constraint~\eqref{eq:perfect-1} is obtained from relaxing the constraints from~\eqref{eq:choice-prob-two-cycle}.  The constraint~\eqref{eq:perfect-2}  is derived from the following equality, with the variable $\gamma^t_i$ replacing the first fractional term and a new variable $\Gamma^t_{ij}$ replacing each fractional term in the summation,
\[
 \re{x^t_i = \frac{x^t_i}{1+\sum_{k \in [N]}u_k x^t_k} + \sum_{j \in [N]}u_j\frac{x_i^tx_j^t}{1+ \sum_{k \in [N]}u_kx^t_k}     \qquad \for \, t \in [M+1]}.
\]
To obtain~\eqref{eq:perfect-3}, we use the no-purchase variable $\rho^t$ to scale the following McCormick envelope of $x_i^tx_j^t$~\citep{mccormick1976computability}:
\[
\max\{0, x^t_i+x^t_j-1 \} \leq  x^t_ix^t_j \leq  \min \{x_i^t,x_j^t\} \qquad \for i\in [N], j\in [N]\setminus \{i\}.
\]
The last constraint~\eqref{eq:perfect-4} is obtained by observing that the relation $x^t_ix^t_i = x^t_i$ holds for binary variables. Unlike the formulation~\eqref{eq:cyclic-policy-conic}, this formulation does not rely on bounds of the no-purchase probability variables. More specifically, in the derivation of~\eqref{eq:cyclic-policy-conic}, we use $\rho^tx^t_i$ as a nonlinear scaling factor, demanding a linearization using bounds of $\rho^t$. In the development of linear inequalities in~\eqref{eq:2cycMILP-boundfree}, we avoid using a nonlinear scaling factor and, instead, directly use $\rho^t$ as a scaling variable to derive~\eqref{eq:perfect-3}. 
 
The following result shows this bound-free formulation is theoretically tighter than the formulation~\eqref{eq:cyclic-policy-conic}. One of the most important properties of an mixed integer program (MIP) formulation is the strength of its natural continuous relaxation that is obtained by ignoring the integrality constraints. This is important because a tighter continuous relaxation often indicates a faster convergence of the branch-and-bound algorithm on which most commercial MIP solvers are built~\citep{vielma2015mixed}. Since~\eqref{eq:total} has a maximization objective, a tighter formulation has a smaller continuous relaxation objective value.  Let $Z_{\textsc{Bound-free}}$ (resp. $Z_{\textsc{Cycle-Conic}}$) be the optimal objective value of the natural continuous relaxation of formulation~\eqref{eq:2cycMILP-boundfree} (resp.~\eqref{eq:cyclic-policy-conic}). The two objective values have the following relationship. 

% Let $F_{\textsc{Bounds-free}}$ (resp. $F_{\textsc{Cycle-Conic}}$) be the feasible region of the natural continuous relaxation of formulation~\eqref{eq:2cycMILP-boundfree} (resp.~\eqref{eq:cyclic-policy-conic}).  Notice that the two feasible regions are expressed in two different spaces. In order to study their containment relation, we project them onto a common space defined by variables $(\rho,\x,\boldsymbol{\gamma})$. Recall that for a set $S$ in the space of $(\x,\y)$ variables, its projection onto the space of $\x$, denoted as $\proj_{\x}(S)$, is defined as $\bigl\{\x\bigm| \exists \y \st (\x,\y) \in S \bigr\}$.

\begin{theorem}\label{them:bound-free}
\re{$Z_{\textsc{Bound-free}} \le Z_{\textsc{Cycle-Conic}}$ if the cross-period constraint is non-overlapping.}
\end{theorem}

\section{Case Study}\label{sec:case:study}
This section presents a case study using transaction data from a corporation cafeteria. We first estimate the history-dependent effect and then use the results to design assortments. Our estimations indicate that customers have a significant but short-term satiation effect. Furthermore, our findings suggest that including history-dependent effects in choice models will lead to greater assortment varieties and higher revenue.

%Our findings demonstrate that incorporating history-dependent effects into choice models can improve the prediction of customer choices and increase assortment variety and revenue.   

\subsection{Data Description and Estimation}\label{sec:estimator}
%\subsubsection{Data Description}
We cooperated with a Chinese company that provides catering services to over $4,000$ corporation cafeterias and were given access to a dataset for one cafeteria that mainly serves workday luncheons for customers who work nearby. The dataset includes $27,144$ transactions for $53$ main dishes, covering $110$ workdays from July to December 2021. The data for each transaction includes the time, dish, and quantity, as well as a card ID if the dish was purchased with a membership card. The dataset also includes information about the price, main ingredient, and flavor of each dish. In what follows, we use the terms dish and product and also the terms menu and assortment interchangeably.

Our dataset indicates that 84.8\% of the transactions came from 418 members, which generated 84.5\% of the total revenue. Thus, we focus on the members' transactions because they are frequent visitors and generate most of the revenue. Our dataset also reveals that the probability of a member dining in the cafeteria more than once for two, three, or four consecutive days was 72.5\%, 87.7\%, and 91.2\%, respectively. This suggests that customers' dish selections could be influenced by past menus.

We predict customers' choices using the MNL model with history-dependent effects. The utility function for product $i$ on day $t$ is 
\[ U_{i}^t =  \boldsymbol{\theta}^T X_{i} +  \theta_{r} r_i + \sum_{m\in [M]} x_i^{t-m} \beta^m
 + \xi_{i}^t, \]
where $X_{i}$ denotes covariates that include a fixed effect and the main ingredient and flavor of product $i$; $r_i\in [0,1]$ is a normalized price; $\boldsymbol{\theta}$ and $\theta_{r}$ are coefficients of $X_{i}$ and price, respectively; $\boldsymbol{\theta}^T X_{i} +  \theta_{r} r_i$ is the base utility $\beta_i^0$ defined in Section~\ref{sec:model}; and $\xi_i^t$s are i.i.d. standard Gumbel random variables. The history-dependent effect is captured by the term $\sum_{m\in [M]} x_i^{t-m} \beta^m$, in which $M$ is the memory length. Our model is a standard MNL model if $M=0$. Note that we set the same history-dependent coefficient for all products because not all products have consecutive offering records, and therefore, we cannot estimate the product-wise history-dependent effect. 

We utilize members' transactional data to calibrate the MNL choice model with the utility function presented above. We split the transaction data into two sets: \re{data from July to October for training and data from November to December for test}. We first use five-fold cross-validation on the first data set to estimate the utility function, and then we optimize weekly menus based on the estimators. Next, we conduct the five-fold cross-validation again, this time on the second data set, to select the value of $M$ with the best out-of-sample performance, and we set this as the ``ground truth" model for evaluation, which is a common approach in policy evaluation~\citep {cao2022network}. To address the fact that our dataset lacks no-purchase records, we employ the expectation-maximization (EM) approach, a widely used technique to deal with censored data~\citep{van2017expectation, csimcsek2018expectation}, to calibrate our model. 
 
Table~\ref{tab:est:res} summarizes the estimation and five-fold validation results for $M\in\{0,1,2,3,4\}$, using the first data set. Column $\theta_r$ presents price coefficient estimators, and columns $\beta^1$ to $\beta^4$ give history-dependent effects. The columns labeled $\chi^2$ and Log-likelihood provide the average out-of-sample value in the five-fold cross-validation (normalized by sample size). The percentage values in parentheses indicate the relative improvement of each model with a history-dependent effect compared to the standard MNL model. 
\begin{table} 
\TABLE 
{\re{Estimators and cross-validation results under different memory length}\label{tab:est:res}}
{
\begin{tabular}{ccccccccc}
\hline 
&   \multicolumn{5}{c}{\text{Estimators}}   & &\multicolumn{2}{c}{\text{Out-of-sample Performance}} \\
\cline{2-6} \cline{8-9}
$M$  & $\theta_r$ &  $\beta^1$ &  $\beta^2$ &  $\beta^3$ &  $\beta^4$ & & $\chi^2$ & \text{Log-likelihood}   \\ \hline 
0 & $-2.558^{***}$  & \text{-}   & \text{-}  & \text{-}  & \text{-}   & & 0.0860     & -1.7796    \\
 & (0.0767)  &    &   &   &   &  &    \\
1   & $-2.550^{***}$  & $-0.475^{***}$   & \text{-}   & \text{-}   & \text{-} & & 0.0828  ($3.67\%$)  & -1.7769 ($0.15\%$)    \\
 & (0.0767)  & (0.0498)  &   &         &       &     &    &\\
2  &$ -2.520^{***}$  & $-0.493^{***}$         & $-0.366^{***}$         & \text{-}              & \text{-}   &  & 0.0822  ($4.38\%$)  & -1.7754 ($0.24\%$)      \\
 & (0.0767)  & (0.0499)   & (0.0568)   & &  &  &  &     \\
3   & $-2.540^{***}$ & $-0.492^{***}$         & $-0.358^{***}$         & $-0.289^{***}$         & \text{-} &   & 0.0819  ($4.71\%$)   & -1.7749   ($0.27\%$)   \\
  & (0.0768)  & (0.0500)  & (0.0568)   & (0.0708)  &   &   &    &    \\
4   & $-2.540^{***}$         & $-0.495^{***}$         & $-0.361^{***}$         & $-0.286^{***}$         & $-0.146^{***}$        & & 0.0826  ($3.91\%$)   & -1.7753    ($0.25\%$)     \\
 & (0.0768)   & (0.0500)    & (0.0569)   & (0.0709)          & (0.0647)    &    &  &   \\ \hline         
\end{tabular}
}
{Note: The numbers in parentheses under the estimators report the standard errors of the estimations. ``-" indicates that a model does not include the corresponding predictor. ${***}p<0.001$. The columns ``$\chi^2$" and ``Log-likelihood" show the average value of the five-fold cross-validation (normalized by sample size) followed by a percentage in parentheses, which is the relative improvement of each model with a history-dependent effect compared to the standard MNL model.
}
\end{table}
Table~\ref{tab:est:res} shows significant negative history-dependent effects, and the scale of the effects decreases as the time interval between the current and historical offerings increases.  
The estimation results are consistent for different values of $M$. For instance, $\beta^1$ is consistently about $-0.49$ for all $M$. This suggests that the attraction value of a product will shrink by about $\exp(-0.49) = 61\%$ if it was offered yesterday. In other words, if a product is sold on two consecutive days, we need to lower its price by at least RMB \re{$ (0.492/2.540)\times 15 \approx 2.9$ when $M=3$} to compensate for the negative history-dependent effect, where we normalize prices by the maximum price RMB $15$ in the estimation. These results demonstrate the importance of considering history-dependent effects in assortment planning. 

The average out-of-sample validation results are reported in the last two columns of Table~\ref{tab:est:res}. They suggest that history-dependent effects can improve prediction power (i.e., lower $\chi^2$ and higher log-likelihoods). \re{When $M=3$, the model presents the lowest $\chi^2$ and the highest log-likelihood.} This result suggests that a small memory length \re{($M=3$)} can capture the history-dependent effect in customers' choice.
 
\subsection{Assortment Planning}\label{sub:ast:plan}
In this section, we use the formulation proposed in Section~\ref{sec:reformulation} to optimize weekly menus using the above estimators. We design weekly workday menus because the cafeteria does not frequently operate on weekends, and in practice, the cafeteria indeed plans for a week. We use the estimated parameters in Table~\ref{tab:est:res} as the input for our conic models. To keep the menu sizes comparable, we impose a daily cardinality constraint of $6$, which is the median assortment size in the dataset. 

\subsubsection{Assortment Patterns}\label{sec:menu:pattern}
\begin{figure} 
 \FIGURE 
{\begin{subfigure}[b]{0.32\textwidth}
         \centering\caption{\re{Offering times}}
         \label{fig:offertime}
\includegraphics[width=\linewidth]{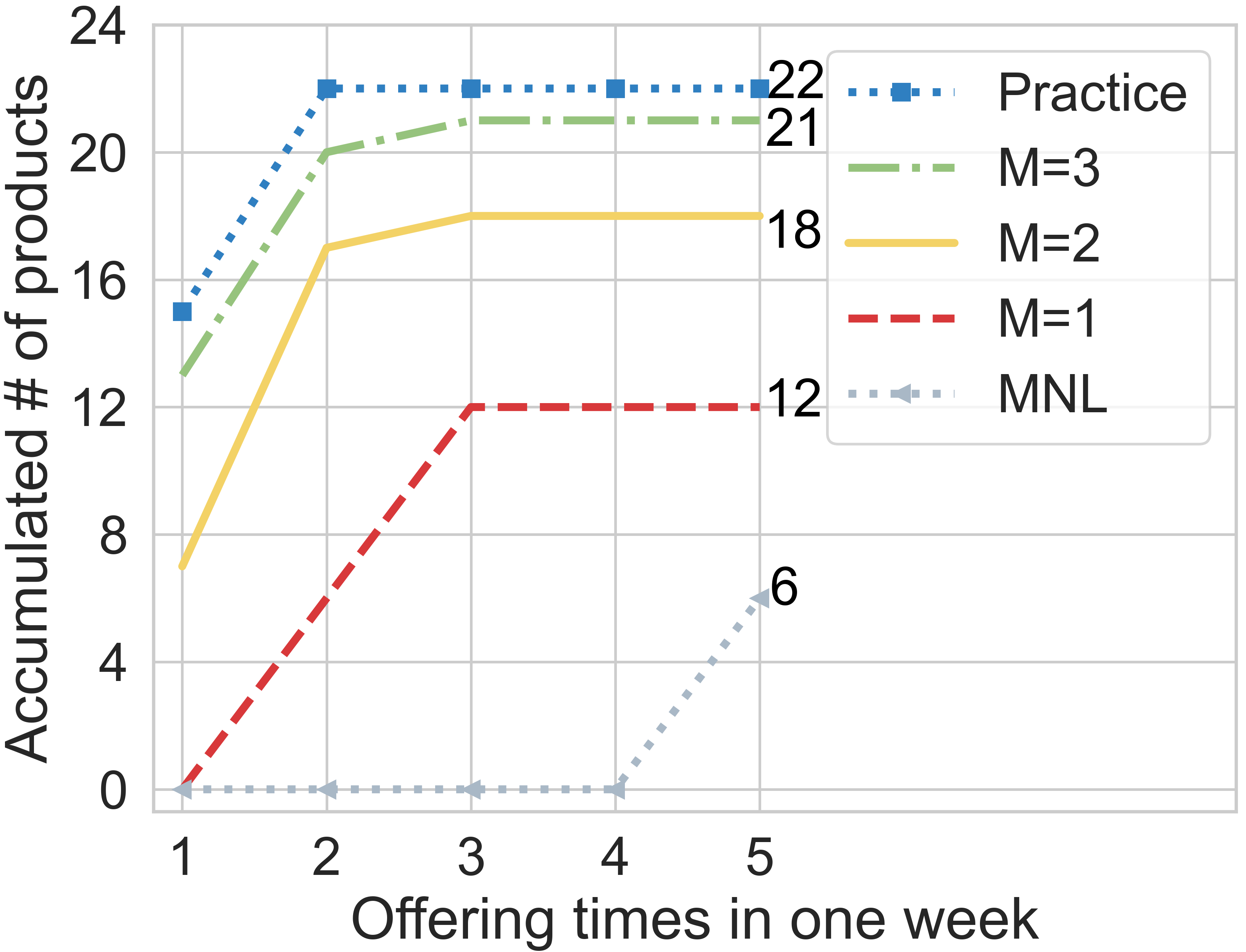}  
     \end{subfigure}
     \hfill 
      \begin{subfigure}[b]{0.32\textwidth}
         \centering
    \caption{\re{Price}}
         \label{fig:price}
\includegraphics[width=\linewidth]{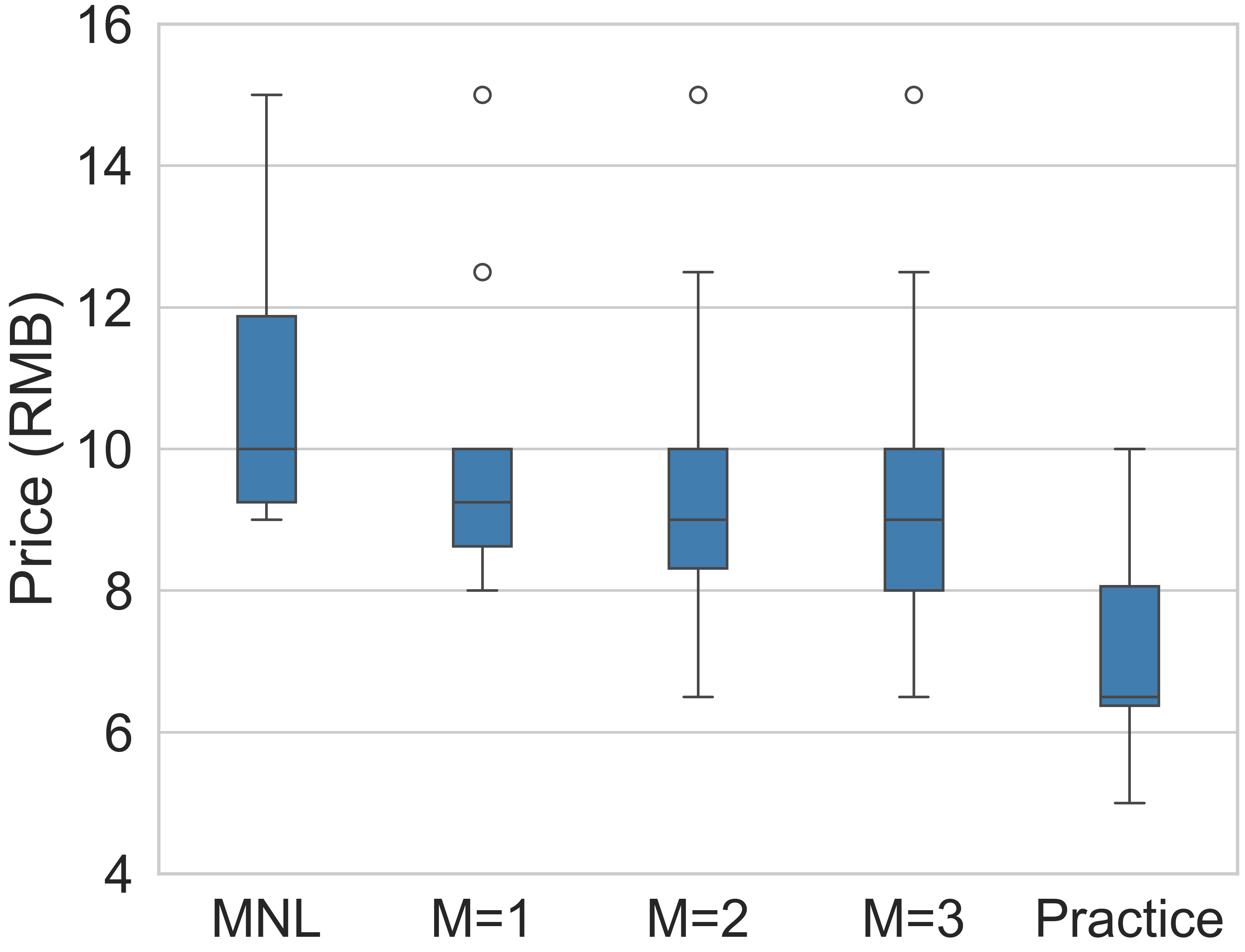}  
     \end{subfigure}
     \hfill
     \begin{subfigure}[b]{0.32\textwidth}
         \centering
         \caption{\re{Base attraction value}}
         \label{fig:attra}
\includegraphics[width=\linewidth]{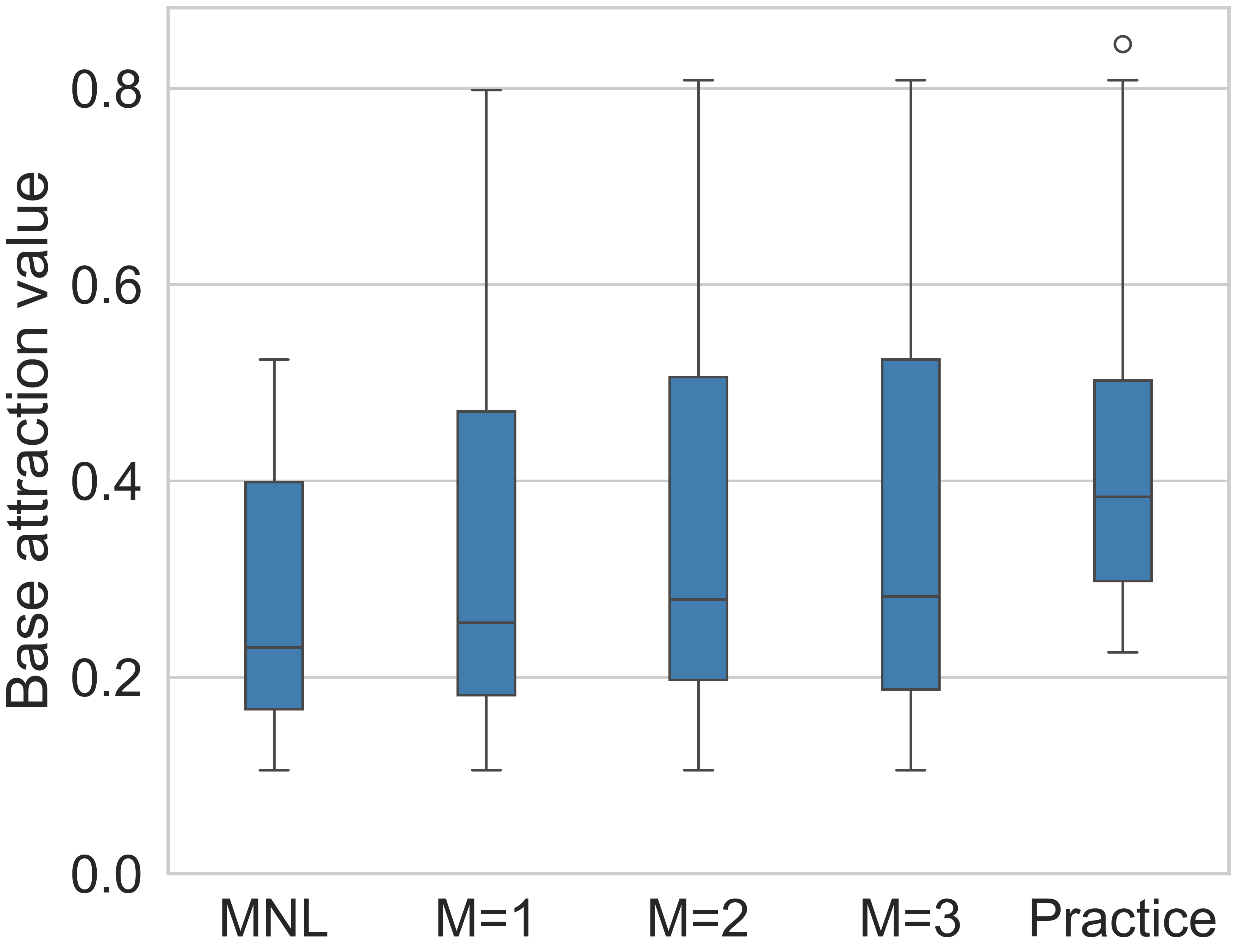} 
     \end{subfigure}}
 {\re{Information of products offered in assortments obtained by different models}\label{fig:menu:freq}}
    {}    
\end{figure}

Figure~\ref{fig:menu:freq} shows the accumulated number of products for different offering times, the prices, and the base attraction values obtained from the assortments under the MNL model, under history-dependent MNL models with different memory lengths, and under an assortment planning used in practice. The current practice assortment was randomly chosen from the second test data set.

\captionsetup[subfigure]{font=small}  
As shown in Figure~\ref{fig:offertime} and~\ref{fig:price}, the MNL-based model tends to offer fewer products with higher prices. Moreover, the number of products that are offered increases with the memory length $M$, because it is beneficial to offer more products and avoid repeating the same product within the memory length when satiation effects last for a long period. When the memory length increases, the assortments include more products with lower prices and higher base attraction values, as shown in Figure~\ref{fig:price} and~\ref{fig:attra}.

\begin{figure}[hbtp]
\FIGURE 
    {\includegraphics[width=0.95
\linewidth]{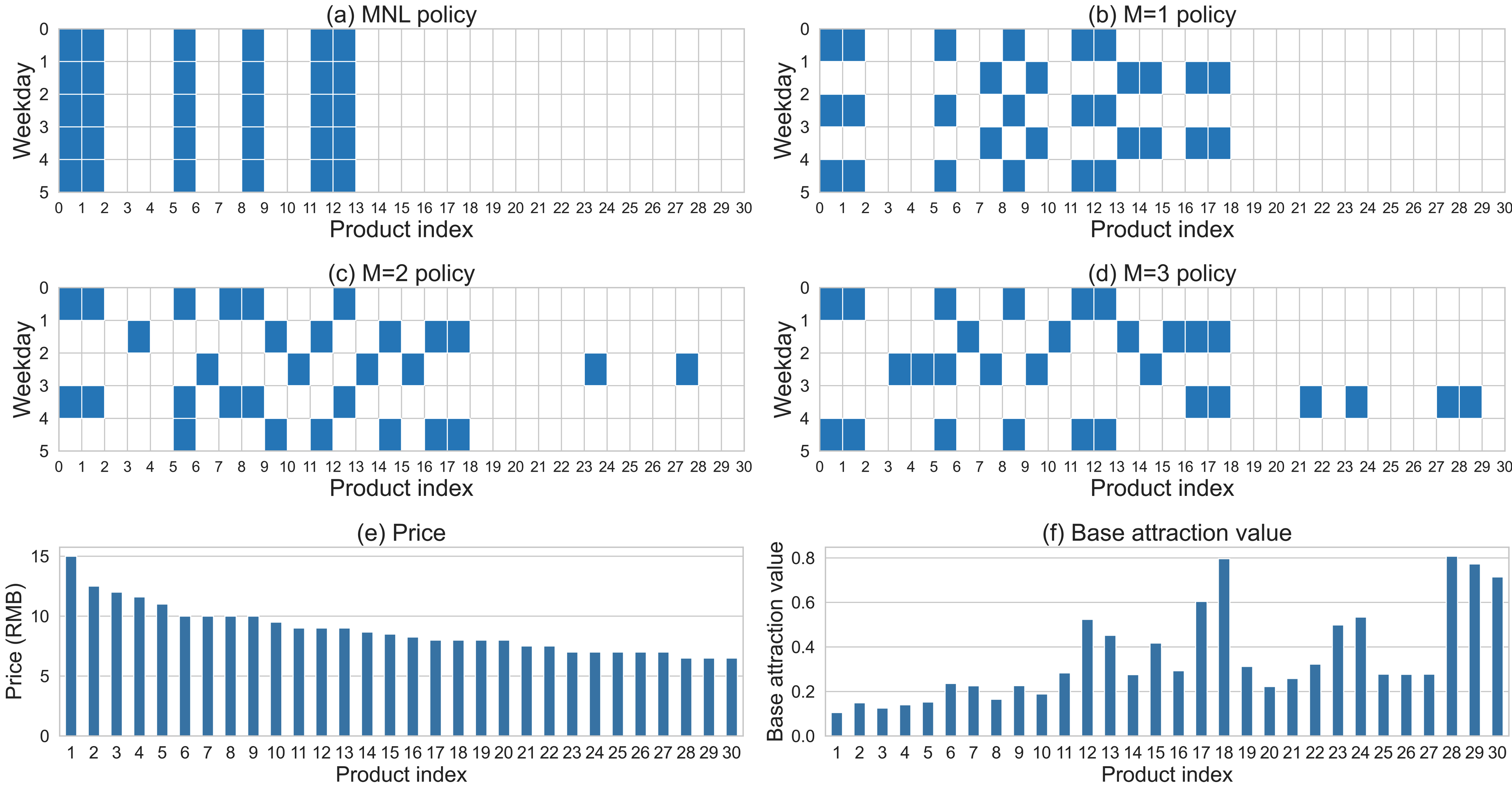}}
{\re{Assortment details of four different policies }
    \label{fig:menu:pattern}}
    {(a)-(d) show the assortments of four polices. The horizontal axes represent products indexed. The vertical axes represent the days of a week. The products offered are colored blue. (e) and (f) show the price and base attraction value of each product.}
\end{figure}
%(We omit the last 23 dishes that are not offered in any assortments.)

Figure~\ref{fig:menu:pattern} presents assortment details for examining how the memory length affects assortment patterns. Figure~\ref{fig:menu:pattern} (a)-(d) show the assortments obtained from MNL-based models with $M\in\{0,1,2,3\}$. The horizontal axes represent products indexed from the highest to the lowest price. The vertical axes represent the days of a week. The products offered are colored blue. The subfigures show that the interval between two adjacent offerings of the same dish grows when $M$ increases, and optimal assortments represent a cyclic-style structure. The interval for two adjacent offerings of most products is the memory length $M$, which suggests that the cycle length is $M+1$, aligning with the $(M+1)$-cyclic policy in Section~\ref{sec:two:cycle}. Moreover, we observe that products with higher prices or lower base attraction values are more sensitive to the satiation effect and thus will be offered at most once within the span of the memory length. For instance, the high-price products 1 to 5 strictly repeat once in intervals of $M$ days to avoid the satiation effect, possibly because the revenue loss caused by the utility deterioration will be scaled up by the high price. On the other hand, products with relatively high base attraction value may possibly be offered more than once within the memory period, e.g., product 17,18. The price and base attraction value for each product are shown in Figure~\ref{fig:menu:pattern}~(e) and~(f) for reference.

%Recall that $\x = (x_i^t)$ represents an assortment planning, a menu in our case study. Each matrix in Figure \ref{fig:menu:pattern} denotes a menu, and positive entries ($x_i^t =1)$ are colored. Dishes are sorted in descending order of price. We also present each dish's price and base attraction value in Figure \ref{fig:menu:pattern} (e) and (f), respectively. For instance, the high-price dishes 1-4 strictly repeat in an interval of $M$ days to avoid the satiation effect because the revenue loss caused by the utility deterioration will scaled up by the high price. On the other hand, dishes 6 and 12 are still provided on two consecutive days when $M=2$ due to their relatively high base attraction value (as shown in Figure~\ref{fig:menu:pattern} (f)), making the two dishes able to resist satiation effects.  

\subsubsection{Performance Comparison}\label{sec:menu:rev}
In this section, we focus on two metrics to evaluate the assortment performance: variety and revenue, as these are concerns for the corporation cafeteria. We adopt the Herfindahl-Hirschman Index (HHI), a measure of market concentration \citep{rhoades1993herfindahl}, to quantify the variety of an assortment planning. The HHI value is $\sum_{i\in[n]}(\frac{h_i}{H})^2$, where $n$ is the total number of products, $h_i$ is the number of times product $i$ is offered, and $H=\sum_{i\in[n]} h_i$. A low HHI value reflects a large variety. To compute the out-of-sample revenue, \re{we randomly sample 20 subsets from the second data set. In each subset, we estimate the ``ground truth" model with $M=2$ that was selected by cross-validation on the second data set and calculate the revenue based on the ``ground truth" demand model. The average revenue across the 20 data sets is used as the out-of-sample revenue to avoid random noise.} Besides the assortments discussed in Section~\ref{sec:menu:pattern}, we additionally compute assortments obtained by the $(M+1)$-cyclic policy with non-overlapping constraints that is presented in Section~\ref{sec:two:cycle}. %we first set set the memory period the model with the best prediction performance for the second data set, as the ``ground truth" model ($M=2$). \re{To make the result more robust, we then randomly sample 20 subsets of the second data set and estimate the ``ground truth" model to perturb its parameters. We compute the average revenue across the 20 data sets as the out-of-sample revenue.} Besides the assortments discussed in Section~\ref{sec:menu:pattern}, 

\begin{figure} 
\FIGURE
{\includegraphics[width=0.7
\linewidth]{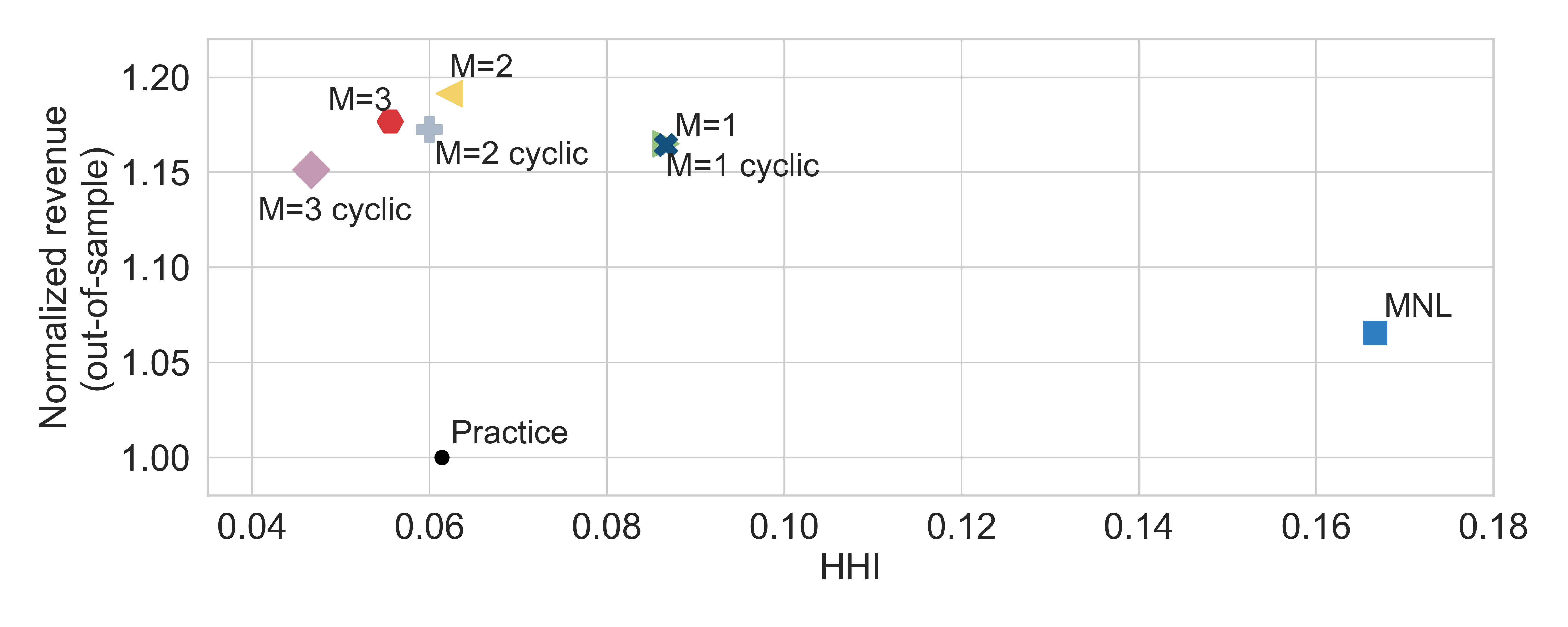}}
{\re{Assortment variety and average revenue in the out-of-sample test}
    \label{fig:menu:rev}}
{The horizontal axis represents the Herfindahl-Hirschman Index value of each policy. A small HHI value means a wide variety. The vertical axis shows the average revenue of a policy across 20 test sets. We normalize the revenue of the current practice to one.}
\end{figure}

Figure~\ref{fig:menu:rev} illustrates the relationship between the variety and the revenue of an assortment plan. The horizontal axis represents the HHI value, and the vertical axis shows the average normalized revenue. We normalize the out-of-sample revenue of different assortment plans with the revenue of current practice as a benchmark. Generally, assortments with a low HHI value can help avoid the negative impacts of satiation effects and maintain high product utility. Figure~\ref{fig:menu:rev} shows that the assortments generated by our model ($M>0$) have a lower HHI and higher revenue than the assortments obtained with the MNL model. For instance, \re{the average revenue obtained using our model with $M=3$ is approximately $10.4\%$ greater than the revenue under the MNL model}. The $(M+1)$-cyclic policy also performs well, \re{with cyclic assortments experiencing only a slight revenue loss compared to the optimal assortments.} We also note that the manually decided assortments in practice have low HHI values but carry substantial revenue losses, as in practice, managers may decide the assortments based on many factors that are not observed. These results demonstrate that a wide variety of products with a low frequency can mitigate the satiation effect. Our model achieves this by efficiently balancing products' utility, revenue, and history-dependent effects.

\section{Numerical Study}\label{sec:numerical}
This section presents our numerical studies on synthetic data. Section~\ref{sec:numerical:speed} illustrates the efficiency of our formulations regarding computation time and optimality gap compared to other formulations. \re{Section~\ref{sec:heuristic} shows that our formulation generates significantly more revenue than benchmark heuristics. In Section \ref{sec:numerical:cycle}, we demonstrate that the bound-free model~\eqref{eq:2cycMILP-boundfree} can scale up to solve large-size instances.} Our numerical studies were implemented in the \texttt{Julia}~\citep{bezanson2017julia}, using modeling language \texttt{JuMP}~\citep{Lubin2023}, on a Virtual Machine with a 32-core Intel Xeon (Skylake, IBRS) @2.39 GHz processor. All MILP models were solved using \texttt{Gurobi 11.0.0}~\citep{gurobi}, and MIECP models were solved using \texttt{Mosek 10.2}~\citep{mosek}.

\subsection{\re{Formulation Efficiency}}\label{sec:numerical:speed}
We consider four alternate formulations for solving~\eqref{eq:total} without constraints and \re{set the bounds of the no-purchase probability variable as $\rho^t_L =  \frac{1}{({1 + \sumj \exp(\beta_j^0 + \sum_{m\in I_j} \beta_j^m ))}}$ and $\rho^t_{U} = 1$ for the first three models}:
\begin{itemize}
    \item [] $\mathsf{Conic}$: This is the mixed-integer exponential cone formulation given by~\eqref{eq:Conic}.
    \item [] $\mathsf{Env}$: This is the mixed-integer linear program obtained by replacing the exponential cone constraint in \eqref{eq:Conic} with the convex envelope based constraints~\eqref{eq:conv-pers-linear1} and~\eqref{eq:conv-pers-linear2}, as detailed in Remark~\ref{rmk:M=2}.
    \item [] $\mathsf{ML}$: This is the mixed-integer linear formulation based on a multilinear extension of attraction value functions, and its derivation is detailed in the ecompanion~\ref{sub:multi:detail}. 
\item [] $\mathsf{SCIP}$: This refers to solving the binary nonlinear model~\eqref{eq:total} by the global optimization solver \texttt{SCIP}.
\end{itemize}

\re{We consider the satiation effects in this numerical study}. For a fixed number of products $N$ and memory length $M$, the base utility of each product $\beta_i^0$ is drawn from a uniform $[-1,1]$ distribution, and its revenue $r_i$ from a uniform $[1,10]$ distribution. The no-purchase utility is $0$. For each $i\in[N]$ and $m\in[M]$, we set the history-dependent effect $\beta_i^m$ from two uniform distributions $U[-1,0]$ and $U[-2,-1]$ to represent a \textit{weak} and \textit{strong} satiation effect, denoted as W and S, respectively. We vary $N\in\{30,50,70\}$, $M\in \{1,2\}$, and $T\in\{5,10\}$. For each combination of $N$, $M$, and a satiation level (W or S), we generate five instances. \re{In~\ref{additional:numerical}, we provide additional computational results on settings where history-dependent effects are mixed, memory lengths are large, and both cross-periods and cross-products constraints are present. To solve our MIECP models with large memory lengths, we propose a cutting-plane algorithm, which is discussed in~\ref{ec:sec:double}.}

We solve each instance of \eqref{eq:total} using one of the formulations $\mathsf{Conic}$, $\mathsf{Env}$, $\mathsf{ML}$, and $\mathsf{SCIP}$. We report the computation time and the end gap, which is defined as  $G_{\textsf{end}} = 100\% \times \frac{R_{\mathsf{U}} -R_{\mathsf{IP}}}{R_{\mathsf{IP}}}$, where $R_{\mathsf{U}}$ and $R_{\mathsf{IP}}$ is the best upper bound and the best integer solution at termination with a time limit of 3600 seconds and an optimality gap tolerance of $0.5\%$, respectively. 

 \begin{table} 
\centering
\caption{\re{Computation time and end gap of four formulations under satiation effects}}\label{tab:speed}
\begin{tabular}{@{\hskip2pt}c@{\hskip2pt}c@{\hskip2pt}c@{\hskip2pt}ccccccccccccc}
\hline 
  &            &            &            & \multicolumn{3}{c}{\textsf{Conic}}                     & \multicolumn{3}{c}{\textsf{Env}}                        & \multicolumn{3}{c}{\textsf{ML}}                        & \multicolumn{3}{c}{\textsf{SCIP}}                      \\
  \cline{5-16} 
\textbf{T}  & \textbf{M} & \textbf{N} & $\boldb$ &  $\# \mathsf{sol}$ & $T_{\mathsf{opt}\small{(s)}}$ & $G_{\mathsf{end}\small{(\%)}} $ &  $\# \mathsf{sol}$ & $T_{\mathsf{opt}\small{(s)}}$ & $G_{\mathsf{end}\small{(\%)}} $ &  $\# \mathsf{sol}$ & $T_{\mathsf{opt}\small{(s)}}$ & $G_{\mathsf{end}\small{(\%)}} $ &  $\# \mathsf{sol}$ & $T_{\mathsf{opt}\small{(s)}}$ & $G_{\mathsf{end}\small{(\%)}} $  \\ \hline
{5}  & 1          & 30         & W          & 5              & 36            & 0            & 5              & 6             & 0            & 5              & 267           & 0            & 0              & -            & 648            \\
{5}  & 1          & 30         & S          & 5              & 90            & 0            & 5              & 96            & 0            & 5              & 527           & 0            & 0              & -            & 920            \\
{5}  & 1          & 50         & W          & 3              & 65            & 0.64         & 5              & 657           & 0            & 3              & 473           & 1.04         & 0             & -            & 1380        \\
{5}  & 1          & 50         & S          & 3              & 1766          & 0.66         & 4              & 1249          & 0.68         & 1              & 3492          & 1.51         & 0              & -          & 1527         \\
{5}  & 1          & 70         & W          & 3              & 453           & 1.75         & 4              & 102           & 2.17         & 3              & 855           & 2.26         & 0              & -           & 1798         \\
{5}  & 1          & 70         & S          & 1              & 1695          & 0.87         & 2              & 757           & 0.68         & 1              & 2523          & 1.87         & 0              & -          & 2191        \\ \hline 
{5}  & 2          & 30         & W          & 5              & 267           & 0            & 5              & 95            & 0            & 3              & 1339          & 0.99         & -              & -          & 913            \\
{5}  & 2          & 30         & S          & 3              & 527           & 2.39         & 3              & 363           & 1.5          & 0              & -             & 11.76        & 0              & -           & 1046         \\
{5}  & 2          & 50         & W          & 3              & 689           & 0.95         & 5              & 160           & 0            & 1              & 2760          & 1.57         & 0              & -           & 1415         \\
{5}  & 2          & 50         & S          & 0              & -             & 2.35         & 1              & 2133          & 1.64         & 0              & -             & 49.13        & 0              & -             & 1756        \\
{5}  & 2          & 70         & W          & 0              & -             & 1.81         & 3              & 811           & 1.6          & 0              & -             & 16.35        & 0              & -             & 2369        \\
{5}  & 2          & 70         & S          & 1              & 1020          & 3.41         & 1              & 57            & 2.37         & 0              & -             & 79.22        & 0              & -          & 2459         \\ \hline 
{10} & 1          & 30         & W          & 0              & -             & 0.92         & 4              & 422           & 0.86         & 0              & -             & 0.94         & 0              & -             & 744         \\
{10} & 1          & 30         & S          & 0              & -             & 2.14         & 1              & 2712          & 1.77         & 0              & -             & 2.86         & 0              & -             & 1005        \\
{10} & 1          & 50         & W          & 0              & -             & 1.35         & 2              & 1994          & 1.4          & 0              & -             & 1.48         & 0              & -             & 1477         \\
{10} & 1          & 50         & S          & 0              & -             & 3.82         & 0              & -             & 3.12         & 0              & -             & 3.88         & 0              & -             & 1577         \\
{10} & 1          & 70         & W          & 0              & -             & 2.44         & 1              & 313           & 2.06         & 0              & -             & 3.27         & 0              & -             & 1907        \\
{10} & 1          & 70         & S          & 0              & -             & 3.08         & 0              & -             & 1.9          & 0              & -             & 4.43         & 0              & -             & 2305         \\ \hline 
{10} & 2          & 30         & W          & 0              & -             & 2.37         & 0              & -             & 1.4          & 0              & -             & 4.15         & 0              & -             & 1001        \\
{10} & 2          & 30         & S          & 0              & -             & 7.52         & 0              & -             & 5.8          & 0              & -             & 165.07       & 0              & -             & 1151         \\
{10} & 2          & 50         & W          & 0              & -             & 2.9          & 0              & -             & 2            & 0              & -             & 63.85        & 0              & -             & 1563          \\
{10} & 2          & 50         & S          & 0              & -             & 8.73         & 0              & -             & 6.05         & 0              & -             & 608.16       & 0              & -             & 1891         \\
{10} & 2          & 70         & W          & 0              & -             & 4.66         & 0              & -             & 3.05         & 0              & -             & 283.46       & 0              & -             & 2756        \\
{10} & 2          & 70         & S          & 0              & -             & 6.51         & 0              & -             & 4.64         & 0              & -             & 750.15       & 0              & -             & 2763  \\\hline       
\end{tabular} 
\end{table}
Table~\ref{tab:speed} summarizes the computation results. Given a parameter combination, we split $5$ instances into two groups: solved and unsolved instances. The solved instances achieved the optimality gap within 3600s, but the unsolved instances failed. $\#\textsf{sol}$ indicates the number of solved instances. The average computation time of solved instances is denoted as $T_{\textsf{opt}}$ and the average end gap of unsolved instances as $G_{\textsf{end}}$. Overall, formulation $\textsf{Env}$ performs best in terms of the number of solved instances, the computation time, and the end gap. Our MILP formulation $\textsf{Env}$ slightly outperforms its conic counterpart, $\textsf{Conic}$, because the former leverages the convex envelope of attraction value functions. Not surprisingly, our formulations $\textsf{Env}$ and $\textsf{Conic}$ dominate $\textsf{ML}$ and $\textsf{SCIP}$ in all aspects. In particular, our formulations solve more instances to global optimality, with fewer computational times, than $\textsf{ML}$ and $\textsf{SCIP}$ do. For unsolved instances, the end gap of our formulations is significantly smaller than that of $\textsf{ML}$ and $\textsf{SCIP}$. 
%\re{We also test the performance of our formulations under constraints and show the results in the ecompanion~\ref{ec:sec:constraint} due to space considerations.} 

% \re{We also test the performance of $\textsf{Env}$ under constraints in the ecompanion~\ref{ec:sec:constraint}. The result shows that our model still performs well under the cross-product constraint but more sensitive to the cross-period constraint.}
\begin{table} 
\centering
\caption{\re{Root gap ($\%$) of continuous relaxation of four formulations under satiation effects}}\label{tab:rootgap}
\begin{tabular}{cccccccccccccccccc}
\cline{1-8} \cline{11-18}
\textbf{T} & \textbf{M} & \textbf{N} & $\boldb$ & \textsf{Conic} & \textsf{Env} & \textsf{ML} & \textsf{SCIP} & \textbf{} & \textbf{} &  \textbf{T} & \textbf{M} & \textbf{N} & $\boldb$ & \textsf{Conic} & \textsf{Env} & \textsf{ML} & \textsf{SCIP} \\  
\cline{1-8} \cline{11-18}
5          & 1          & 30         & W          & 3.86           & 3.22        & 156.73      & \text{-}          &  \text{  }   &      & 10 & 1 & 30 & W & 4.38  & 3.65  & 179.23  & \text{-}  \\
5          & 1          & 30         & S          & 6.81           & 5.89        & 301.08      & \text{-}         &         &     & 10 & 1 & 30 & S & 7.83  & 6.56  & 343.77  & \text{-}  \\
5          & 1          & 50         & W          & 2.89           & 2.5         & 295.99      & \text{-}          &         &     & 10 & 1 & 50 & W & 3.3   & 2.82  & 336.43  & \text{-}  \\
5          & 1          & 50         & S          & 5.64           & 5.04        & 517.19      & \text{-}          &        &      & 10 & 1 & 50 & S & 7.11  & 5.64  & 588.26  &\text{-}  \\
5          & 1          & 70         & W          & 3.95           & 3.55        & 400.53      & \text{-}          &         &     & 10 & 1 & 70 & W & 5.17  & 4.77  & 455.25  & \text{-}  \\
5          & 1          & 70         & S          & 3.84           & 3.45        & 713.79      & \text{-}          &          &    & 10 & 1 & 70 & S & 5.05  & 4.2   & 810.67  & \text{-}  \\
\cline{1-8} \cline{11-18}
5          & 2          & 30         & W          & 4.57           & 4.57        & 224.85      & \text{-}          &          &    & 10 & 2 & 30 & W & 5.64  & 5.61  & 275.29  & \text{-}  \\
5          & 2          & 30         & S          & 13.17          & 13.17       & 587.57      & \text{-}        &          &    & 10 & 2 & 30 & S & 15.4  & 15.32 & 786.57  & \text{-} \\
5          & 2          & 50         & W          & 4.77           & 4.77        & 430.75      & \text{-}          &           &     & 10 & 2 & 50 & W & 5.9   & 5.76  & 539.85  &\text{-}   \\
5          & 2          & 50         &S          & 10.1           & 10.06       & 921.6       & \text{-}          &      &     & 10 & 2 & 50 & S & 13.34 & 12.09 & 1419.53 & \text{-} \\
5          & 2          & 70         & W          & 5.47           & 5.3         & 717.98      & \text{-}          &         &     & 10 & 2 & 70 & W & 6.94  & 6.46  & 940.42  & \text{-}  \\
5          & 2          & 70         & S          & 6.85           & 6.43        & 1243.68     & \text{-}          &          &    & 10 & 2 & 70 & S & 9.02  & 7.85  & 1890.09 & \text{-} \\
\cline{1-8} \cline{11-18}
\end{tabular}
\end{table}

Besides the computation time and end gap, it is also helpful to compare the continuous relaxation of different formulations. We compute the root gap by $ 100\% \times \frac{R_{\textsf{Rlx}} -R_{\textsf{IP}}}{R_{\textsf{IP}}}$,
% \[
% G_{\textsf{root}} = 100\% \times \frac{R_{\textsf{Rlx}} -R_{\textsf{IP}}}{R_{\textsf{IP}}}, 
% \]
where $R_{\mathsf{Rlx}}$ is the objective value of the continuous relaxation of a given formulation. A small root gap indicates a tight formulation. Table~\ref{tab:rootgap} reports the average root gap of \textsf{Conic}, \textsf{Env}, and \textsf{ML} under different parameter configurations. Overall, the root gap of our formulations \textsf{Conic} and \textsf{Env} is significantly smaller than that of \textsf{ML}. It indicates the tightness of our formulations. Moreover, the root gap of our formulations is stable and does not scale with the problem size. For instance, the root gap of \textsf{Env} is no more than $7\%$ when $M=1$ and no more than $16\%$ when $M=2$.

\subsection{\re{Performance of Heuristic Policies}}\label{sec:heuristic}
\re{We compare the revenue generated by our formulation against two heuristic policies without constraints.}
\begin{itemize}
    \item [] \re{$\mathsf{Sequential-RO}$: The sequential RO-policy shown in Algorithm~\ref{alg:greedy}. It can serve as a heuristic policy for settings with satiation effects.} 
    \item [] \re{$\mathsf{Sequential-LOSPO}$: For each period, this policy drops one satiation product that causes the greatest revenue reduction from the base assortment obtained from the $\mathsf{Sequential-RO}$ policy. We refer to it as \textit{sequential leave-one-satiation-product-out} (Sequential-LOSPO) policy.}
\end{itemize} 

\re{We consider instances with mixed addiction-satiation effects. A product has an addiction (resp. satiation) effect with probability $\theta$ (resp. $1-\theta$), where $\theta\in \{0,0.1,0.2\}$. When $\theta=0$, all products have satiation effects. For a product with addiction (resp. satiation) effect, its $\{\beta_i^m\}_{m\in [M]}$ are uniformly sampled from $[0,1]$ (resp. $[-2,-1]$). We fix $M=2$ and generate five instances for each combination of $N,T$, and~$\theta$.}

\re{We compare revenue obtained by the two heuristic policies and our formulation. We define the relative revenue gap as $100\% \times \frac{R_{\mathsf{Env}} - R_{\mathsf{heuristic}}}{R_{\mathsf{Env}}}$, where $R_{\mathsf{heuristic}}$ is revenue generated by one of two heuristic policies. $R_{\mathsf{Env}}$ is computed by $\mathsf{Env}$ within 3600 seconds. Table \ref{tab:rev} shows that the revenue gap of both heuristics is not negligible, especially when $\theta$ is small. For instance, the smallest relative gap is at least more than 10\% if $\theta=0$, suggesting that our formulation has more advantage in scenarios where the satiation effect is dominant. Although the revenue gap shrinks when $\theta$ increases, our formulation still holds value as it has an optimality guarantee and is easily extended to constrained cases.}

\begin{table}[hbtp]
\centering 
\caption{\re{Revenue gap ($\%$) of heuristic policies under different history-dependent effects with $M=2$}}\label{tab:rev}
\begin{tabular}{cccccccccc}
\hline 
& &  \multicolumn{3}{c}{$\mathsf{Sequential-RO}$} & & \multicolumn{3}{c}{$\mathsf{Sequential-LOSPO}$}   \\ \cline{3-5} \cline{7-9}  
\textbf{N}  & \textbf{T}  &$\theta=0$      & $\theta=0.1$    & $\theta=0.2$  &  &$\theta=0$      & $\theta=0.1$    & $\theta=0.2$     \\  \hline 
{30} & 5  & 18.60 & 13.02 & 2.99 &   & 15.09 & 11.77 & 2.60 \\  
30  & 10 & 31.37 & 19.43 & 4.18 & & 28.47 & 18.65 & 3.81 \\
50 & 5  & 14.39 & 5.18  & 2.02 &  & 12.90 & 4.63  & 1.84 \\
50 & 10 & 25.72 & 7.73  & 2.82 & & 24.20 & 7.27  & 2.66 \\
70  & 5  & 12.31 & 4.00  & 3.59 &   & 10.91 & 3.80  & 3.40 \\
70   & 10 & 21.86 & 5.74  & 5.03 &   & 20.61 & 5.59  & 4.88 \\ \hline 
\end{tabular}
\end{table} 

% \begin{tabular}{cccccccccccc}
% \cline{1-5} \cline{7-11}
%  & &  \multicolumn{3}{c}{$\mathsf{Sequential-RO}$}   &   & & \multicolumn{3}{c}{$\mathsf{LOSPO}$}   \\ \cline{1-5} \cline{7-11}
% \textbf{N}  & \textbf{T}  &$\theta=0$      & $\theta=0.1$    & $\theta=0.2$     &  & \textbf{N}  & \textbf{T}  &$\theta=0$      & $\theta=0.1$    & $\theta=0.2$     \\ \cline{1-5} \cline{7-11}
% 30         & 5  & 22.93 & 16.09 & 3.14  &  & 30 & 5  & 49.54 & 31.03 & 11.89   \\
% 30         & 10 & 46.02 & 27.08 & 4.47  &  & 30 & 10 & 66.25 & 44.08 & 15.83   \\
% 50         & 5  & 16.84 & 5.83  & 2.12  &  & 50 & 5  & 49.42 & 21.58 & 3.65    \\
% 50         & 10 & 34.76 & 9.52  & 3.03   &  & 50 & 10 & 68.00 & 30.49 & 4.44    \\
% 70         & 5  & 14.07 & 4.23  & 3.79  &  & 70 & 5  & 41.43 & 12.00 & 9.25    \\
% 70         & 10 & 28.05 & 6.26  & 5.49   &  & 70 & 10 & 56.97 & 16.07 & 11.75   \\
% \cline{1-5} \cline{7-11}
% \end{tabular}

\subsection{\re{Computing $(M+1)$-Cyclic Policies}}\label{sec:numerical:cycle}
In this section, we demonstrate the efficiency and scalability of the bound-free model~\eqref{eq:2cycMILP-boundfree} for computing \re{$(M+1)$-cyclic policies under the non-overlapping constraint.} Note that due to the presence of $\boldsymbol{\Gamma}$ variables, the number of variables in~\eqref{eq:2cycMILP-boundfree} is $\mathcal{O}(N^2)$. To address this issue, we propose a projected-cutting-plane algorithm to obtain a projected approximation of \eqref{eq:2cycMILP-boundfree}. The cutting plane algorithm \re{first sets a base formulation without using $\boldsymbol{\Gamma}$ variables, that is, constraints~\eqref{eq:perfect-2}-\eqref{eq:perfect-4},} and solves its continuous relaxation to obtain a continuous solution. Then, the algorithm calls a separation oracle to generate cuts violated by the continuous solution, \re{which are selected from the projection of constraints \eqref{eq:perfect-2}-\eqref{eq:perfect-4} onto the space of the base formulation, {i.e.}, the space of $(\x,\boldsymbol{\rho},\boldsymbol{\gamma})$ variables}. Adding new cuts to the base model leads to a tighter formulation. Running the separation oracle once is a round of cut generation. We can introduce new cuts over $K$ rounds and refer to the resulting formulation as BF-$K$, where ``BF" stands for ``bound-free". \re{We refer to~\ref{sub:sepa:oracle} for a detailed description of the projected cutting-plane algorithm.}  %Throughout the iterations, the formulation always has the same integer feasible region with the one of~\eqref{eq:2cycMILP-boundfree}, as the base one is. Moreover, the cutting plane algorithm generates a hierarchy of tighter formulations which converges to the base model tightened with constraints \eqref{eq:perfect-2}-\eqref{eq:perfect-4} in terms of the optimal objective value of the continuous relaxation~\citep{kelley1960cutting}. 

% and solve its continuous relaxation to obtain a continuous solution. The algorithm then uses a separation oracle to generate new cuts. The new cuts are the projection of constraints \eqref{eq:perfect-2}-\eqref{eq:perfect-4} in~\eqref{eq:2cycMILP-boundfree}, violated by the fractional solution, to the space of the base formulation. Adding new cuts into the base model leads to a tighter formulation. Running the separation oracle once is a round of cut generation. We can add new cuts for $K$ rounds and the final formulation as {BF}-$K$, where BF denotes ``bound-free". \re{Finally, we add the binary constraint back and solve it for an optimal integer solution.}

% Indeed, one can repeat this process (i.e., a large $K$) until the fractional solution of {BF}-$K$ does not violate \eqref{eq:perfect-2}-\eqref{eq:perfect-4}, then add binary constraints to obtain an integer optimal solution. 

\re{Clearly, a larger $K$ leads to a tighter formulation at the cost of adding cuts. To select a proper $K$ that balances the formulation tightness and the formulation size, we compute the optimality gap closed after the $K^\text{th}$ round as follows 
\[
\textsf{GClosed}_{K}  = 100\% \times  \frac{ \textsf{Opt}_0 -\textsf{Opt}_K}{  \textsf{Opt}_0- \textsf{Opt}_{\infty}},
\]
where $\textsf{Opt}_K$ is the optimal objective value of the continuous relaxation of {BF}-$K$, and the index $0$ (resp. $\infty$) refers to the base (resp. limiting) model. Table \ref{tab:margin} reports the gap closed by {BF}-$K$ for $K\in\{1,2,3\}$ on instances with $M\in \{1,2\}$ and $N\in \{50,100\}$. {BF}-$1$ closes about $80\%$ of the optimality gap, suggesting setting $K=1$ is good enough.}

% \re{
% how much gap between formulation~\eqref{eq:2cycMILP-boundfree} and the base formulation regarding the optimal objective value of their continuous relaxation is reduced by {BF}-$K$, defined as follows.
% \[
% \textsf{GRed}_{K}  = 100\% \times  \frac{\textsf{Opt}_K - \textsf{Opt}_0}{\textsf{Opt}_{\infty} - \textsf{Opt}_0},
% \]
% where $\textsf{Obj}_K$, $\textsf{Obj}_0$, and $\textsf{Obj}_{\infty}$ are the optimal objective value of the continuous relaxation of {BF}-$K$, the base model, and~\eqref{eq:2cycMILP-boundfree}, respectively. We observe that the marginal reduced gap is decreasing with $K$. For instance, in our numerical study ($M\in \{1,2\}$ and $N\in \{50,100\}$), Table \ref{tab:margin} shows that {BF}-$1$ has the most marginal reduced gap and the subsequent reduced gaps are little. After adding one-round cuts, the gap between the continuous relaxation of the tight formulation~\eqref{eq:2cycMILP-boundfree} and the base model shrinks about $80\%$. Hence, we set $K=1$ and use {BF}-$1$ to compute $(M+1)$-cyclic policy.}
 \captionsetup[subtable]{font=small} 
\begin{table}[!ht] 
 \caption{\re{The computational performance of~\eqref{eq:2cycMILP-boundfree} with a projected-cutting-plane implementation}}\label{tab:bf:model}
\begin{subtable}{0.48\textwidth}
\centering \caption{\re{Performance of the cutting-plane algorithm}}
\label{tab:margin}
\begin{tabular}{ccccc}
\hline 
\textbf{N}   &  \textbf{M}     & $\textsf{GClosed}_{1}$ &$\textsf{GClosed}_{2}$ &$\textsf{GClosed}_{3}$  \\ \hline 
50 & 1 & 79.25 & 91.19 & 96.21 \\  
50          & 2 & 82.24 & 93.14 & 95.90 \\
100         & 1 & 89.82 & 95.61 & 97.58 \\
100         & 2 & 77.24 & 90.30 & 93.48 \\ \hline 
\end{tabular}
\end{subtable}
\hfill
\begin{subtable}{0.48\textwidth}
\centering \caption{\re{Performance of bound-free model {BF}-$1$}}
\label{tab:bf1}
\begin{tabular}{cccccc}
\hline 
\textbf{N}   &  \textbf{M}  &  $\# \mathsf{sol}$ & $T_{\mathsf{opt}\small{(s)}}$ & $G_{\mathsf{end}\small{(\%)}} $     \\ \hline 
50 & 1 & 5 & 0.57   & 0.30    \\   
50          & 2 & 5 & 70.74  & 0.46 \\
100         & 1 & 5 & 8.19   & 0.42    \\
100         & 2 & 5 & 378.59 & 0.49  \\ \hline 
\end{tabular}
\end{subtable}
\end{table}
\re{Table~\ref{tab:bf1} reports the performance of {BF}-$1$. The result demonstrates the scalability and efficiency of the bound-free model compared to~\eqref{eq:Conic}. BF-1 solves all instances to optimality within 400 seconds. In contrast, under the number of products $N=50$ and memory length $M=2$, both $\textsf{Conic}$ and $\textsf{Env}$ fail to solve some of the instances as shown in Table~\ref{tab:speed}.}

\section{Conclusion}\label{sec:conclusion}

This paper studies a multi-period assortment planning problem with history-dependent customer utility. We prove that the problem is NP-hard under negative history-dependent effects. \re{To address this for general history-dependent effects}, we develop an MIECP formulation to find a global optimal solution. We also prove the optimality of the sequential policy for the addiction case and the asymptotic optimality of cyclic policies while relating the cycle length to customer memory length.

% We first show that the problem is solvable in polynomial time if and only if customers' utility increases with historical repeated offerings. The optimal policy is a sequential revenue-ordered policy. However, the problem becomes NP-hard if the historical effect is the opposite, even with a planning horizon of two periods. We develop an MIECP formulation to find a global optimal solution. Our numerical analyses indicate that our formulation efficiently achieves optimality and generates higher revenue than heuristic policies, especially under strong satiation effects. We further explore cyclic policies when the problem has a large time horizon. We theoretically prove that a cyclic policy with a cycle equal to the memory length plus one is asymptotically optimal. Additionally, we develop a neat, bound-free reformulation for non-overlapping cyclic policies when the memory length equals one. 

Our research can be extended in several directions. \re{First, it is important to develop efficient approximation algorithms, such as fully polynomial-time approximation schemes, to address the history-dependent assortment planning problem, especially in cases with large memory lengths and long planning horizons.} Second, extending the model to online settings for personalized assortments \re{or dynamic settings to adjust assortments to quickly respond to} stochastic demand. Third, investigating joint assortment, pricing, or inventory problems under history-dependent effects. Finally, exploring history-dependent effects in other choice models for assortment optimization.

% Our research can be extended in several directions. The first direction would be to explore more efficient algorithms than MIECP, i.e., the fully polynomial-time approximation scheme. Second, it would be worth to explore the joint assortment and pricing or inventory problem under history-dependent effects. The third direction would be to extend our model to the online settings --- providing personalized assortments when customers stochastically revisit the online store. Finally, it would be interesting to examine how to incorporate history-dependent effects into different choice models and optimize assortment decisions.

% \theendnotes

% \subsection{Extension: customer heterogeneity}
% customers may have different memory-span of historical assortments. To denote such heterogeneity, we assume there are $K$ customers types or segments. Each customer type is associated with a nonnegative weight $\lambda_k$. $k$-th type of customers could memorize the past $k$ periods assortments. Hence we have the following optimization model. 
% \begin{align} \label{model:diffM} 
%    \max \limits_{S_t\in \X, \forall t\in \T} & \sum_{t \in [T]} \sum_{k=1}^K \lambda_k R(S_t| S_{t-1},\dots,S_{t-k}).
%  \end{align} 

%\ACKNOWLEDGMENT{%
% Enter the text of acknowledgments here
%}% Leave this (end of acknowledgment)
\bibliographystyle{informs2014_r2} % outcomment this and next line in Case 1
\bibliography{AssortmentPaper_v2} 

% the ecompanion here
% Options are (1) APPENDIX (with or without general title) or 
%             (2) APPENDICES (if it has more than one unrelated sections)
% Outcomment the appropriate case if necessary
%
\clearpage
% \begin{APPENDIX}{title}
% \end{APPENDIX}
%
%   or 
%

\clearpage

\ECSwitch

\ECDisclaimer
%%%%%%%%%%%%%%%%%%%%%%%%%%%%%%%%%%%%%%%%%%%%%%%%%%%%%%%%%%

%%% Main head for the e-companion
\ECHead{Additional Discussions, Missing Proofs, Formulations, and Numerical Studies}

\setcounter{section}{8} % 重置 section 计数器
\renewcommand{\thesection}{EC.\the\numexpr\value{section}-8\relax} % 自定义编

% \section*{Additional Discussions, Missing proofs, and Formulations}

\setcounter{table}{8} % 重置表格计数器
\setcounter{figure}{8} % 重置图片计数器
% 重新编号表格和图片
\renewcommand{\thetable}{EC.\the\numexpr\value{table}-8\relax} % 自定义表格编号格式
\renewcommand{\thefigure}{EC.\the\numexpr\value{figure}-8\relax} % 自定义图片编号格式

\section{Meituan's Grocery Flash Sale Channel}\label{sec:ec:meituan}
\begin{figure}[hbtp]
\FIGURE 
{
 \begin{subfigure}[b]{0.32\textwidth}
      \caption{\centering{Main page} }
    \label{fig:main} 
    \includegraphics[width= 0.92 \linewidth]{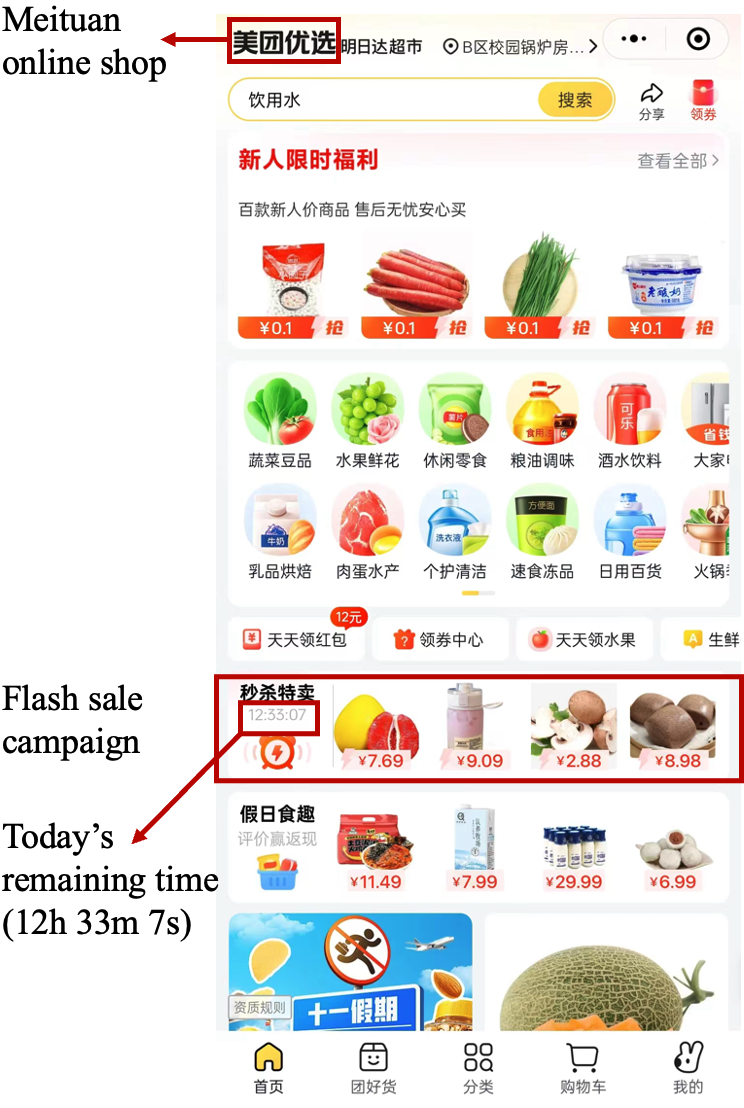}
  \end{subfigure}
  \hfill
 \begin{subfigure}[b]{0.32\textwidth}
  \caption{  \centering{Product page} }
    \label{fig:product} 
    \includegraphics[width= 0.97 \linewidth]{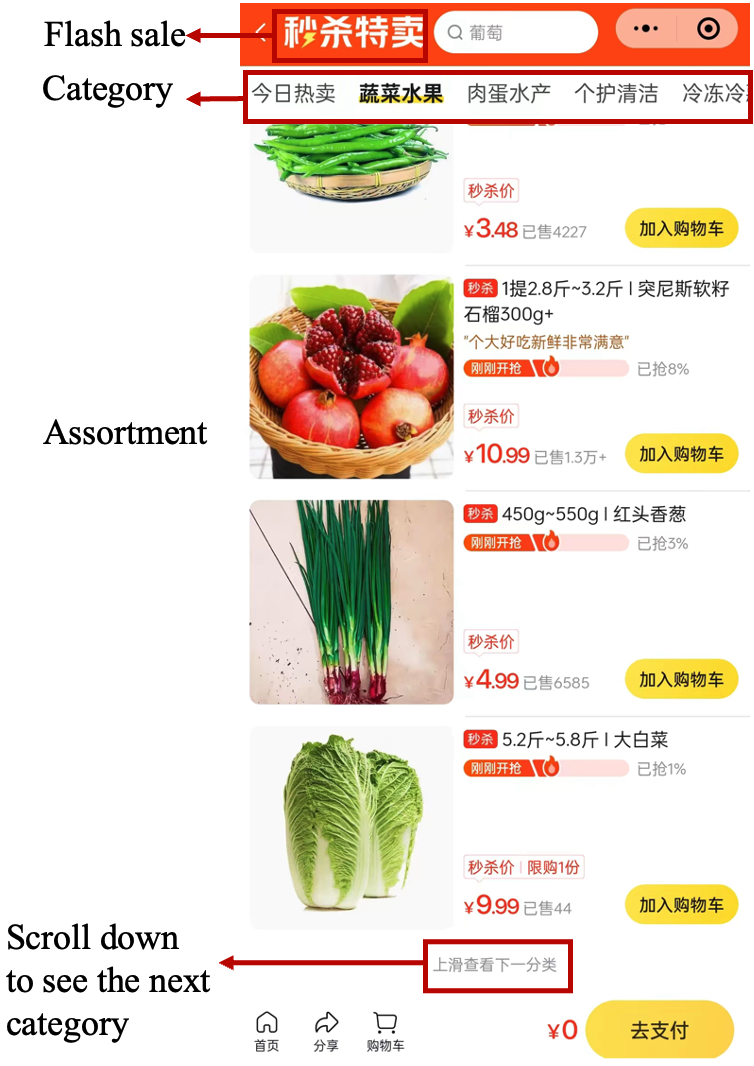}
  \end{subfigure}
  \hfill
 \begin{subfigure}[b]{0.32\textwidth}
   \caption{ \centering{Payment page} }
    \label{fig:payment} 
    \includegraphics[width=  \linewidth]{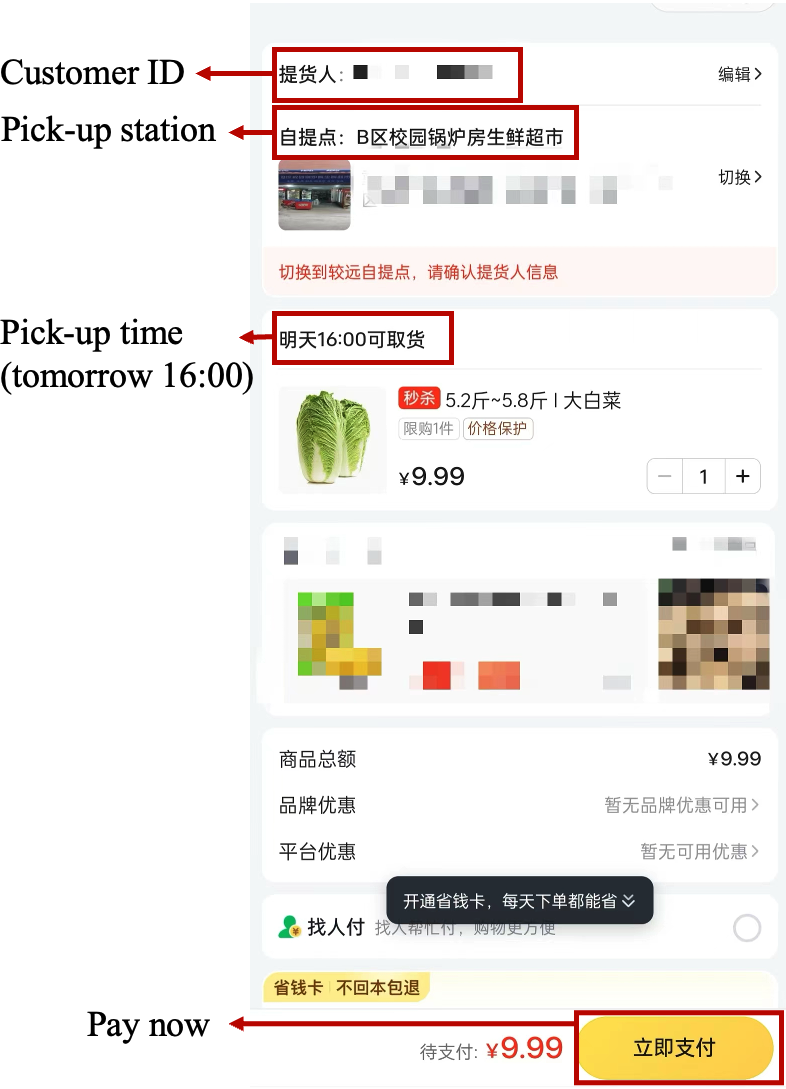}
  \end{subfigure} 
}  
 {Screenshots of Meituan's grocery flash sale channel \label{fig:meituan}}
 {}
\end{figure}

{Figure~\ref{fig:meituan} details Meituan's flash sale channel. Figure~\ref{fig:main} is the main page of the Meituan online store. The main page lists the entry of the grocery flash sale channel and highlights its remaining time in the current campaign. The flash sale usually offers assortments on a daily basis. Hence, the remaining time in our example is no more than 24 hours. Figure~\ref{fig:product} is the screenshot of the flash sale's products. In each category, only a limited number of products are offered. Please note that the products are the same for all customers.}

{Meituan aggregates all orders during a flash sale campaign and forwards the order information to suppliers once the current campaign ends. Suppliers will then prepare the products and deliver them to the pick-up stations the next day. Based on our discussion with Meituan, they must determine the products for flash sales for a particular planning period ahead of time because Meituan needs lead time to negotiate contracts with suppliers. As shown in Figure~\ref{fig:payment}, the payment page details the delivery and pick-up information, such as the customer's chosen pick-up station and the expected pick-up time. This time indicates when the product is expected to reach the pick-up station. If a customer opts for home delivery, the products will initially be sent to the pick-up station, after which the station's staff will deliver them to the customer's address. This case shows that inventory is not the key concern of the platform because all the products are delivered from suppliers to the pick-up station the next day.}

\section{Proofs of Section \ref{sec:model}}
\subsection{Proof of  Proposition~\ref{prop:nphard}}
\begin{repeattheorem}[  Proposition~\ref{prop:nphard}.]
\eqref{eq:total} is NP-hard even when the planning horizon is two, the memory length is one, and history-dependent effects are negative.   
\end{repeattheorem}
\noindent \noindent {\bf Proof.} To prove the NP-hardness of~\eqref{eq:total}, we will reduce an NP-hard problem---the \textit{3/4-Partition} problem, to a special case of~\eqref{eq:total}. The 3/4-Partition problem is defined as follows:
\begin{itemize}
    \item \textbf{Input:} A set $ \{c_1,c_2,\ldots,c_N\}$ of non-negative integers.
    \item \textbf{Output:} True if and only if there exists a subset $S \subseteq [N]$ such that $\sum_{i\in S}c_i = (3/4) \sum_{i\in [N]} c_i$.
\end{itemize} We will show the NP-hardness of the 3/4-Partition problem by showing that any instance of a \textit{Partition} problem, a well-known NP-hard problem, can be reduced to an instance of the 3/4-Partition problem. The Partition problem is defined as follows:
\begin{itemize}
    \item \textbf{Input:} A set $ \{w_1,w_2,\ldots,w_N\}$ of non-negative integers.
    \item \textbf{Output:} True if and only if there exists a subset $S \subseteq [N]$ such that $\sum_{i\in S}w_i = \sum_{i\in [N]/S} w_i$.
\end{itemize}
Notice that $\sum_{i\in S }w_i = \sum_{i\in [N]/S} w_i$ if and only if $\sum_{i\in S}w_i =(1/2) \sum_{i\in [N]}w_i$. Given an instance of Partition problem, we solve it by solving a 3/4-Partition problem as follows. Define $c_i = w_i$ for $i\in [N]$ and $c_{N+1} = \sum_{i\in[N]}w_i$. Then, solve 3/4-Partition with input $(c_1, \ldots, c_{N+1})$. If the output is true and returns an index set $S$ such that $\sum_{i \in S}c_i = \frac{3}{4}\sum_{i \in [N+1]}c_i$ then the output for Partition is true, that is, $\sum_{i \in S \setminus \{N+1\}} w_i = \sum_{i \in [N+1] \setminus S}w_i$ is a partition. The correctness follows from the fact that the set $S$ contains $N+1$. If the output is false then the output for Partition is also false since otherwise adding $c_{N+1}$ into one of two sets yields a 3/4 partition.

Next, we introduce a special case of~\eqref{eq:total}. In such a special case, the planning horizon $T = 2$, the memory length $M = 1$, and all products have identical revenue, that is, $r_i = r>0$ for $i\in[N]$. In addition, we assume that the history-dependent effects are homogeneous, that is, for every $i \in [N]$, $\exp(\beta_i^1) = k$  for some $k \in (0,1)$. Then,  it turns out that~\eqref{eq:total} becomes 
 \begin{equation}\label{eq:nphard-special}
\max \limits_{\x^1,\x^2\in \{0,1\}^N} \frac{r}{2} \left\{ 
\frac{\sumi x_i^1\exp(\beta^0_i) }{ 1+ \sumi x_i^1 \exp(\beta^0_i)}  + \frac{\sumi x_i^2 \exp\bigl(\beta^0_i + x_i^1 \beta_i^1\bigl)}{1+\sumi x_i^2 \exp\bigl(\beta^0_i + x_i^1\beta_i^1 \bigl)} \right\}.     \tag{\textsc{Special}}
 \end{equation}

We simplify~\eqref{eq:nphard-special} by showing that the optimal assortment in the second period is $x_i^2=1$ for $i\in[N]$. Since all products have the same revenue, we can focus on maximizing the purchase probability. Given an assortment $\x^1\in\{0,1\}^N$, the purchase probability at the first period is a constant, and the purchase probability at the second period increases with $\sumi x_i^2 \exp\bigl(\beta^0_i + x_i^1 \beta_i^1\bigl)$. Since $\exp\bigl(\beta^0_i + x_i^1 \beta_i^1\bigl)>0$ for each $i\in[N]$, then the optimal solution at the second period is $x_i^2=1$ for $i\in [N]$. Let $\nu_i = \exp(\beta_i^0)$ be the base attraction value and $V= \sumi  \exp(\beta_i^0) = \sumi \nu_i$ denote the total base attraction value. Thus, solving~\eqref{eq:nphard-special} equals to solving the following model:
\begin{equation}\label{eq:first:reduce-3}
    \max   \limits_{\x^1 \in \{0,1\}^N}    \frac{r}{2} \left\{ 
                    \frac{\sumi x_i^1 \nu_i }{ 1+ \sumi x_i^1 \nu_i}  + \frac{ V- \sumi x_i^1 \nu_i + k \sumi x_i^1 \nu_i}{1+ V-\sumi x_i^1 \nu_i  + k  \sumi x_i^1 \nu_i } \right\}.
\end{equation}  
  
Next, we show that the 3/4-Partition problem can be reduced to our special case~\eqref{eq:nphard-special}. First, define $\nu_i = a c_i$ for $i\in[N]$, $k=1/4$, and $C=\sum_{i\in [N]}c_i$, where 
\[a = \frac{4( 1-\sqrt{1-k} )}{ C (3(1-k)+3\sqrt{1-k}-4)} = \frac{16-8\sqrt{3}}{C(6\sqrt{3}-7)}>0.\] 
Hence, $V=\sum_{i\in[N]} \nu_i  = a C$. Second, set the target average revenue as 
\[\frac{r}{2}(\frac{3aC}{4+3aC} + \frac{7aC}{16+7aC}). \]
The Partition problem indeed has a solution if and only if there exists a $\x^1 \in \{0,1\}^N$ such that 
\begin{align*}
 &  \frac{r}{2} \left\{ 
                    \frac{\sumi x_i^1 \nu_i }{ 1+ \sumi x_i^1 \nu_i}  + \frac{ V- \sumi x_i^1 \nu_i + k \sumi x_i^1 \nu_i}{1+ V-\sumi x_i^1 \nu_i  + k  \sumi x_i^1 \nu_i } \right\} \\
                    = &   \frac{r}{2} \left\{ 
                    \frac{ a \sumi x_i^1 c_i }{ 1+ a \sumi x_i^1 c_i}  + \frac{ aC- (1-k)a \sumi x_i^1 c_i  }{1+ aC- (1-k)a\sumi x_i^1 c_i   } \right\} \\
                    \le & \max \limits_{y\in [0,C]} \frac{r}{2} \left\{\frac{ay}{1+ay} + \frac{aC- (1-k) a y}{1+aC-(1-k) ay}\right\} = \frac{r}{2}(\frac{3aC}{4+3aC} + \frac{7aC}{16+7aC}).
\end{align*}
This holds because if we view the objective function in the last row as a one-dimension function, denoted as 
\[h(y):= \frac{r}{2}\left(\frac{ay}{1+ay} + \frac{aC-(1-k)ay}{1+aC-(1-k)ay}\right),\] $h(y)$ is concave in $y$ and achieves a unique maximum at $y=3C/4$. \hfill \Halmos

\def \bnu {\boldsymbol{\nu}}
\def \bmu {\boldsymbol{\mu}}

\section{\re{Proofs in Section \ref{eq:solutions}}}

\subsection{\re{Proof of Proposition~\ref{prop:conc-alpha}}}
\re{To prove this proposition, we will invoke the following technical lemma. 
\begin{lemma}\label{EC:lemma:aff}
    Let $f:  \{0,1\}^n \to \bbmr$ and $L:\bbmr^n \to \bbmr^n$ be an invertible affine transformation such that $L(\{0,1\}^n) = \{0,1\}^n$, and consider a composition $(f \mcirc L)(\x) := f(L(\x))$ for every $\x \in \{0,1\}^n$. Then, 
    \[
    \conc(f)(\x) = \conc(f \mcirc L)\bigl( L^{-1}(\x) \bigr) \quad \for \x \in [0,1]^n.
    \]
\end{lemma}
}
\noindent {\bf Proof.} \re{
Let $g(\x):= \conc(f \mcirc L)\bigl( L^{-1}(\x)\bigr)$ for $\x \in [0,1]^n$. 
Clearly, $ g(\x) = f(\x)$ for every $\x \in \{0,1\}^n$, and $g(\cdot)$ is a concave function. Thus, 
we show that  $g(\x) \geq \conc(f)(\x)$ for every $\x\in [0,1]^n$ since the concave envelope is the smallest concave function on $[0,1]^n$ that overestimates  $f(\cdot)$ on $\{0,1\}^n$.}

\re{
To show the opposite direction, we need the dual definition of the concave envelope~\cite{rockafellar1997convex}, that is,
\[
\conc(f)(\x) = \max\biggl\{\sum_{v \in V} f(v)\lambda_v \biggm| \x = \sum_{v \in V}v\lambda_v,\ \sum_{v \in V} \lambda_v = 1,\ \lambda_v \geq 0 \; \for v \in V \biggr\},
\]
where $V:=\{0,1\}^n$. Now, consider a point $\bar{\x} \in [0,1]^n$ and let $\bar{\y} = L^{-1}(\bar{\x})$. We will argue that $\conc(f\mcirc L)(\bar{\y}) \leq \conc(f)(\bar{\x})$. By the dual definition of concave envelope, there exists a convex multiplier $\{\bar{\lambda}_v\}_{v \in V}$ such that $\bar{\y} =  \sum_{v \in V} v \bar{\lambda}_v $ and $\conc(f \mcirc L)(\bar{\y}) = \sum_{v \in V} (f \mcirc L)(v) \bar{\lambda}_v $. It follows readily that $\bar{\x} = L (\bar{\y}) = \sum_{v \in V} L(v) \bar{\lambda}_v$, and thus
\[
\conc(f)(\bar{\x}) \geq  \sum_{v \in V}f\bigl(L(v)\bigr)  \cdot \bar{\lambda}_v = \sum_{v\in V} (f \mcirc L)(v) \cdot \bar{\lambda}_v = \conc(f \mcirc L)(\bar{\y}),
\]
where the first inequality holds due to the dual definition of $\conc(f)(\cdot)$ and $\bar{\x}$ can be expressed as $\sum_{v \in V} L(v) \bar{\lambda}_v$. 
\hfill \Halmos
}
\begin{repeattheorem}[Proposition~\ref{prop:conc-alpha}]
 \re{Constraint~\eqref{eq:conc-pers} is equivalent to the following system of linear inequalities.}
\re{\begin{equation*} 
     y^t_i \leq \alpha_i(\h^\sigma_{i,0} ) ( \gamma^t_i - \tilde{z}_{i\sigma(1)}^t) + \sum_{k \in [M]} \alpha_i(\h^\sigma_{i,k} ) (   \tilde{z}_{i\sigma(k)}^t - \tilde{z}_{i\sigma(k+1)}^t)   \qquad  \text{ for } \sigma \in \Omega\\
\end{equation*}where $\tilde{z}_{i\sigma(M+1)}^t = 0$, $\tilde{z}^t_{i\sigma(k)} =    z^t_{i\sigma(k)}$ if $\sigma(k) \notin  I_i$, and $\tilde{z}^t_{i\sigma(k)} =  \gamma_i^t -  z^t_{i\sigma(k)}$ if $ \sigma(k) \in I_i$. } 
\end{repeattheorem}

\noindent {\bf Proof.}  \re{We start with the satiation case. Since $\boldsymbol{\beta}_i$ are non-positive, by Corollary 3.14 in~\cite{tawarmalani2013explicit},  $\alpha_i(\cdot)$ is supermodular over $(x_i^{t-1},\dots, x_i^{t-M})$.  Then, it follows from Theorem 3.3 in \cite{tawarmalani2013explicit} that the concave envelope of $\alpha_i(\cdot)$ is given as:
\begin{equation}\label{eq:Lovasz}
    \conc(\alpha_i)(x_i^{t-1},\ldots,x_i^{t-M}) = \min \limits_{\sigma \in \Omega}  \biggl\{    \alpha_i(\boldsymbol{w}^\sigma_{0}  ) ( 1 - x_i^{t-\sigma(1)}) + \sum_{k \in [M]} \alpha_i(\boldsymbol{w}^\sigma_{k}  ) (  x_{i}^{t-\sigma(k)} -  x_{i}^{t-\sigma(k+1)}) \biggr\},
\end{equation} 
where $\boldsymbol{w}^\sigma_{0} = \boldsymbol{0}$, $  \boldsymbol{w}^\sigma_{k} = \boldsymbol{w}^\sigma_{k-1} + \e_{\sigma(k)}$ for $k\in [M]$, and $x_i^{t-\sigma(M+1)}=0$. It follows readily that its perspective function can be represented as follows under the presence of constraint~\eqref{eq:relax2}:
\begin{equation*}
    \pers (\conc(\alpha_i) )(\gamma_i^t,\z_i^t) = \min \limits_{\sigma \in \Omega}  \biggl\{    \alpha_i(\boldsymbol{w}^\sigma_{0}  ) ( \gamma_i^t - z_{i\sigma(1)}^{t}) + \sum_{k \in [M]} \alpha_i(\boldsymbol{w}^\sigma_{k}  ) (  z_{i\sigma(k)}^{t} - z_{i\sigma(k+1)}^{t}) \biggr\},
\end{equation*}
where $z_{i\sigma(M+1)}^t=0$}

\re{For the general case of mixed satiation and addiction effects, we transform the variable $x_i^{t-m}$ to $1-x_i^{t-m} $ if $\beta_i^m >0$, and obtain a transformed attraction value function, defined as follows:
\begin{align*}
    \bar{\alpha}_i(w^1,\dots, w^M) & := \exp\biggl( \beta_i^0 + \sum_{k\in  I_i} \beta_i^k (1-w^k)  + \sum_{k \in [M]\setminus  I_i} \beta_i^k w^k  \biggr) \\
    & = \exp\biggl(\beta_i^0 + \sum_{k\in  I_i}\beta_i^k + \sum_{k\in [M]} (-\vert \beta_i^k \vert)w^k \biggr).
\end{align*} 
Now, the coefficients $-\vert \beta^k_i\vert$  are non-positive, and  by Corollary 3.14 in~\cite{tawarmalani2013explicit}, the transformed function  $\bar{\alpha}_i(\cdot)$ is supermodular. Thus, its concave envelope can be described using~\eqref{eq:Lovasz}. }

\re{Next, we utilize $\conc(\bar{\alpha}_i)(\cdot)$ to characterize the concave envelope of the original attraction value function  $\alpha_i(\cdot)$. Define $\tilde{x}_i^{t-m}  = x_i^{t-m}$ if $m\in [M]\setminus  I_i$, $\tilde{x}_i^{t-m}  =1- x_i^{t-m}$ if $m \in  I_i$, and $\tilde{x}_i^{t-\sigma(M+1)}=0$. We then build the following connections:
\begin{align*}
    \conc(\alpha_i)(x_i^{t-1},\ldots,x_i^{t-M}) & = \conc(\bar{\alpha}_i)(\tilde{x}_i^{t-1} ,\dots,\tilde{x}_i^{t-M} ) \\
    & =  \min \limits_{\sigma \in \Omega}  \biggl\{    \bar{\alpha}_i(\boldsymbol{w}^\sigma_{0}  ) ( 1 - \tilde{x}_i^{t-\sigma(1)}) + \sum_{k \in [M]} \bar{\alpha}_i(\boldsymbol{w}^\sigma_{k}  ) (  \tilde{x}_i^{t-\sigma(k)} - \tilde{x}_i^{t-\sigma(k+1)}) \biggr\}\\
    & = \min \limits_{\sigma \in \Omega}  \biggl\{    \alpha_i(\h^\sigma_{i,0} ) ( 1 - \tilde{x}_i^{t-\sigma(1)}) + \sum_{k \in [M]} \alpha_i(\h^\sigma_{i,k} ) (  \tilde{x}_{i}^{t-\sigma(k)} -  \tilde{x}_{i}^{t-\sigma(k+1)}) \biggr\}
\end{align*}
The first equality follows from Lemma~\ref{EC:lemma:aff}, the second equality holds by~\eqref{eq:Lovasz}, and the last equality holds since it follows from the definition of $\h$ in~\eqref{eq:historyassortment} that  $\bar{\alpha}_i(\boldsymbol{w}^{\sigma}_k) = \alpha_i(\h^{\sigma}_{ik}) $. }

\re{Finally, we scale the concave envelope to obtain 
\[\pers (\conc(\alpha_i) )(\gamma_i^t,\z_i^t) = \min \limits_{\sigma \in \Omega}  \biggl\{  \alpha_i(\h^\sigma_{i,0} ) ( \gamma^t_i - \tilde{z}_{i\sigma(1)}^t) + \sum_{k \in [M]} \alpha_i(\h^\sigma_{i,k} ) (   \tilde{z}_{i\sigma(k)}^t - \tilde{z}_{i\sigma(k+1)}^t)\biggr\}.\]
This completes the proof. }\hfill \Halmos

\subsection{Proof of Theorem~\ref{them:decomposition}}\label{thm:proof:miecp}
\begin{repeattheorem}[Theorem~\ref{them:decomposition}]
 \eqref{eq:Conic} is a mixed-integer exponential cone formulation of problem~\eqref{eq:total}.    
\end{repeattheorem}

\noindent {\bf Proof.} By lifting, we equivalently decompose~\eqref{eq:choice} into the choice constraints~\eqref{eq:extended}:
 \begin{align*}
     & \gamma_{i}^t = \rho^t x_{i}^t && (2a) \\
& z_{im}^t = \gamma_{i}^t x^{t-m}_{i} \quad \for \, m \in [M]  && (2b)\\
&y^t_i = \pers(\alpha_{i})( \gamma^t_i , \z_{i}^t)&& (2c).
 \end{align*}
Thus, proving the equivalence of ~\eqref{eq:Conic} and~\eqref{eq:total} suffices to show that constraints \eqref{eq:relax1}~\eqref{eq:relax2}~\eqref{eq:conv-pers} and~\eqref{eq:pers-extension} construct an equivalent representation of the choice constraints~\eqref{eq:extended}. Clearly, since $\x$ is binary,~\eqref{eq:relax1} (resp. ~\eqref{eq:relax2}) models constraint \eqref{eq:extended-1} (resp. \eqref{eq:extended-2}). Next, we show constraints~\eqref{eq:conv-pers} and~\eqref{eq:pers-extension} provide an exact representation of~\eqref{eq:extended-3}. 
% By the definition of perspective function, equation \eqref{eq:extended-3} is 
% \[ y_i^t \begin{cases}
% \gamma_i^t \exp( \frac{\beta_i^0 \gamma_i^t+ \summ \beta_i^m Z_{im}^{t}}{\gamma_i^t})    &  \gamma_i^t >0 \\
%   0     & \gamma_i^t =0 \bigr.
% \end{cases}\] 
% Given $\x \in \{0,1\}^{N\times T}$ and constraints \eqref{eq:relax1} and~\eqref{eq:relax2}, the above equation equals to 
% \[  y_i^t =  \begin{cases}
%  \exp(\beta_i^0) \rho^t x_{i}^t \exp( \summ \beta_i^m x_i^{t-m})   & x_i^t =1 \\
%   0     & x_i^t =0. \end{cases} \]
By $\x \in \{0,1\}^{N\times T}$ and constraints \eqref{eq:relax1} and~\eqref{eq:relax2}, constraint~\eqref{eq:conv-pers} equals to 
\begin{align} \label{eq:convex:lower}
    y_i^t \ge   \rho^t x_{i}^t \exp( \beta_i^0+\summ \beta_i^m x_i^{t-m}).
\end{align}
On the other hand, the constraint~\eqref{eq:pers-extension} is the concave envelope of $\alpha_i(\cdot)$ scaled by $\gamma_i^t$. That is, $y_i^t \le \gamma_i^t \cdot \conc(\alpha_i)(x_i^{t-1},\dots,x_i^{t-M}) $. In our case, the concave envelope $\alpha_i$ is the tightest concave extension $\alpha_i$~\citep[see][Theorem 6]{tawarmalani2002convex}. 
Thus, $\conc(\alpha_i)(x_i^{t-1},\dots,x_i^{t-M}) = \alpha_i(x_i^{t-1},\dots,x_i^{t-M})$ for $x_i^{t-m} \in \{0,1\}$ and $m\in[M]$. Therefore, constraint~\eqref{eq:pers-extension} equals to $ y_i^t \le \gamma_i^t \cdot \alpha_i(x_i^{t-1},\dots,x_i^{t-M}) $. Combined with constraint~\eqref{eq:relax1}, constraint~\eqref{eq:pers-extension} equals to 
\begin{align}\label{eq:concave:upper}
    y_i^t \le  \rho^t x_{i}^t \exp(\beta_i^0 +  \summ \beta_i^m x_i^{t-m}).
\end{align}
Constraints~\eqref{eq:convex:lower} and \eqref{eq:concave:upper} suggest that constraints~\eqref{eq:conv-pers} and~\eqref{eq:pers-extension} are exact representation of \eqref{eq:extended-3}.

Recall that exponential cones are three-dimensional convex cones:
\[\mathcal{K}_{\exp} = \bigl\{ (x_1,x_2,x_3)\mid x_1\ge x_2 \exp(x_3/x_2),x_2>0 \bigr\}  \cup \bigl\{ (x_1,0,x_3)\mid x_1\ge 0, x_3\le 0\bigr\}. \]
By introducing new variable $w_i^t = \beta_i^0 \gamma_i^t + \summ \beta_i^m z_{im}^t$, \eqref{eq:conv-pers} can be represented by exponential cones:
\[ (y_i^t,\gamma_i^t,w_i^t) \in \mathcal{K}_{\exp}.\]
Therefore, \eqref{eq:Conic} is an equivalent mixed-integer exponentical cone programming of~\eqref{eq:total}. 

\subsection{Proof of Theorem \ref{thm:seq:rev:opt}}\label{ec:prop:seq:ro}
\begin{repeattheorem} [Theorem \ref{thm:seq:rev:opt}]
Assume the absence of cross-product and cross-period constraints. Then, the sequential-revenue-ordered policy solves~\eqref{eq:total} if the history-dependent effects are non-negative, that is, $\boldb\ge \boldsymbol{0}$.
\end{repeattheorem}
Before proving the result, we present some preliminary results that will be used in our proof. Given a vector of attraction values $\bnu \in \bbmr_+^N$, consider a single-period assortment optimization problem under the MNL model given as follows:
\begin{equation}\label{single-ast-mnl}
    \max \limits_{S} \biggl\{\MNL(S , \bnu )  := \sum_{i\in S } r_i \nu_i /\bigl(1 +\sum_{i\in S }  \nu_i  \bigr) \biggm| S \subseteq [N]  \biggr\}.
\end{equation}We assume that products are indexed
such that $r_1 \geq \cdots \geq r_N$. For $\kappa \in [N]$, we will use $[\kappa]$ to denote a revenue-ordered (RO) assortment. \cite{talluri2004revenue} show the optimal revenue is achieved by a RO assortment $[\kappa]$ for some cut-off product $\kappa \in[N]$. Without loss of generality, we assume that there exists a unique optimal solution; otherwise, we will use the optimal RO assortment with the largest number of products. 
\begin{lemma}\label{eclemma:monotone}
% [\citealp{talluri2004revenue,liu2008choice,chen2021optimal}]
Let $\bnu \in \bbmr^{N}_+$ and assume that $r_1 \geq \cdots \geq r_N$. Then, we obtain the following results:
\begin{enumerate} 
    % \item \label{eclemma:monotone-1}  For a given attraction value $\bnu \in \bbmr^{N}_+$, assortment $O_k$ is optimal if and only if 
    % \[ 
    % \MNL(O_1,\bnu)\le  \MNL(O_2,\bnu) \le \dots \le  \MNL(O_k,\bnu) > \MNL(O_{k+1},\bnu) \ge \dots \MNL (O_N,\bnu).
    % \]
    \item  \label{eclemma:monotone-2} The RO assortment $[\kappa]$ is  optimal if and only if  $r_\kappa \ge \MNL([\kappa],\bnu)$ and  $r_{\kappa+1} < \MNL([\kappa],\bnu)$.
    
   \item  \label{eclemma:monotone-3} If the RO assortment $[\kappa]$ is optimal then, for every $\bmu$ such that $\bmu \geq \bnu$, $\MNL([\kappa], \bmu ) \ge \MNL([\kappa],\bnu )$.
\end{enumerate}
\end{lemma}
\noindent {\bf Proof of Lemma~\ref{eclemma:monotone}.} 
% Part~\ref{eclemma:monotone-1} follows from Proposition 1 of~\citep{chen2021optimal}, and
Part~\ref{eclemma:monotone-2} follows from~\citep{talluri2004revenue}. To prove Part~\ref{eclemma:monotone-3}, we consider $\bmu$ such that $\mu_j \geq \nu_j$ for some $j \in [N]$ but $\mu_i = \nu_i$ for all $i \neq j$, and observe that 
\begin{equation*}
    \MNL([\kappa], \boldsymbol{\mu}) -  \MNL([\kappa], \bnu)  =\begin{cases}
        \frac{\mu_j- \nu_j}{1+ \sum_{i\in [\kappa]} \mu_i}  \bigl(r_j - \MNL([\kappa],\bnu) \bigr)   &  \text{ if }  j \in  [\kappa] \\
        0 & \text{otherwise}.
    \end{cases}
\end{equation*} 
Due to Part~\ref{eclemma:monotone-2}, we obtain that  $r_j \ge \MNL([\kappa], \bnu) $ for $j \in [\kappa]$. This completes the proof.  \hfill \Halmos 

The next result establishes that assortments generated by our policy are nested. 
\begin{lemma}\label{eclemmanested}
    Let $[\kappa^1], \ldots, [\kappa^T]$ be a sequence of RO assortments generated by the sequential-revenue-ordered policy. Then, $\kappa^1 \geq \kappa^2 \cdots  \geq \kappa^{T}$. 
\end{lemma}
\noindent {\bf Proof of Lemma~\ref{eclemmanested}.}  Let $\bnu = (\bnu^1, \ldots, \bnu^T)$ be corresponding attraction values, that is, 
\[
\nu^t_{i} := \alpha_i \bigl( \ind(i\in  [\kappa^{t-1}]) ,\dots, \ind(i\in [\kappa^{t-M}]) \bigr) \qquad \for t \in [T] \text{ and } i \in [N]. 
\] 
We prove that $N \geq \kappa^1 \geq \cdots \geq \kappa^t \geq 1$ holds for every $t \in [T]$. Clearly, it holds for $\tau = 1$. Assume that it holds for some $1 \leq \tau <T$.
Now, we observe that 
\[
r_{\kappa^{\tau}+1} < \MNL([\kappa^\tau], \bnu^\tau) \leq \MNL([\kappa^\tau], \bnu^{\tau+1}) \leq \MNL([\kappa^{\tau+1}], \bnu^{\tau+1}),
\]
where the first inequality follows from Part~\ref{eclemma:monotone-2} of Lemma~\ref{eclemma:monotone}, the second inequality holds due to  Part~\ref{eclemma:monotone-3} of Lemma~\ref{eclemma:monotone} since $\boldsymbol{\beta} \geq 0$ and $\kappa^1 \geq  \cdots \geq \kappa^\tau$ imply that $\bnu^{\tau +1} \geq \bnu^\tau$, and the last inequality follows from the optimality of the RO assortment $[\kappa^{\tau+1}]$. Thus, $r_{\kappa^{\tau}+1} < \MNL([\kappa^{\tau+1}], \bnu^{\tau+1})$. Therefore, by Part~\ref{eclemma:monotone-2} of Lemma~\ref{eclemma:monotone}, we can conclude that $\kappa^{\tau +1 } \leq \kappa^{\tau}$. \hfill \Halmos

Now, we are ready to prove Theorem \ref{thm:seq:rev:opt}.
\def \zu { Z_{\mathsf{U}}}
\def \ru { R_{\mathsf{U}}}
\def \ku { k_{\mathsf{U}}}
\def \kro { R_{\mathsf{RO}}}
\def \zro { Z_{\mathsf{RO}}}
\def \v { \boldsymbol{v}}
\def \vu { \boldsymbol{v}_{\mathsf{U}}}
\def \vro { \boldsymbol{v}_{\mathsf{RO}}}

\noindent {\bf Proof.} Let $Z_*$ and $Z_{\mathsf{RO}}$ be the revenue generated from solving~\eqref{eq:total} and from using the sequential revenue-ordered policy. Clearly, our policy generates a feasible plan and thus $Z_{\mathsf{RO}} \leq Z_*$. To show the reverse, we consider the revenue generated from solving the following optimization problem,
\[
Z_{\mathsf{U}} :=  \frac{1}{T}\sumt \max_{S^t} \biggl\{ \MNL(S^t,  \bnu_{\mathsf{U}}^t ) \biggm| S^t \subseteq [N] \biggr\},
\]
where $\nu^t_{\mathsf{U},i}$ is the attraction value of product $i$ in period $t$ obtained by assuming that  the product $i$ is offered in all past periods, that is, 
\[
\re{\nu^t_{\mathsf{U},i} := \alpha_i\bigl( \ind({t-1>0}), \ldots, \ind({t-M}>0) \bigr). }
\]
Moreover, for each $t \in [T]$, let $[\kappa^t_{\mathsf{U}}]$ denote the RO assortment that maximizes $\MNL(\cdot, \bnu^t_{\mathsf{U}})$.  To complete the proof, we will show that $Z_* \leq Z_{\mathsf{U}} \leq Z_{\mathsf{RO}}$.

Now, we show that $Z_* \leq Z_{\mathsf{U}}$. Let $S^1_{*}, \ldots, S^T_*$ be an optimal solution of~\eqref{eq:total}. Let $\bnu_* = (\bnu^1_*, \ldots \bnu^T_*)$ denote the corresponding attraction values, that is,  for $t \in [T]$ and $i\in [N]$, 
\[
\nu^t_{*,i} := \alpha_i \bigl( \ind(i\in S^{t-1}_*) ,\dots, \ind(i\in S^{t-M}_*) \bigr),
\]
where $S^s := \emptyset$ for $s \leq 0$. Then, 
\[
Z_* = \frac{1}{T} \sum_{t \in [T]}\MNL(S^t_*, \nu^t_*). 
\]
Now, for $t \in [T]$, let $\kappa^t_* \in [N]$ such that the RO assortment $[\kappa^t_*]$ maximizes $\MNL(\cdot, \bnu^t_*)$. Then, $Z_* \leq Z_{\mathsf{U}}$ follows by observing that 
\[
\begin{aligned}
\MNL(S^t_*,\bnu^t_*) \leq  \MNL\bigl([\kappa^t_*],\bnu^t_*\bigr) \leq  \MNL\bigl([\kappa^t_*],\bnu^t_{\mathsf{U}}\bigr) \leq  \MNL\bigl([\kappa^t_{\mathsf{U}}],  \bnu_{\mathsf{U}}^t \bigr) ,    
\end{aligned}
\]
where the first (resp. last) inequality holds since the right-hand-side is the optimal revenue that can be achieved under the attraction value $\bnu^t_*$ (resp. $\bnu^t_{\mathsf{U}}$), and the second inequality follows from Part~\ref{eclemma:monotone-3} of Lemma~\ref{eclemma:monotone} since $\boldsymbol{\beta}\geq 0$ implies $\bnu^t_{\mathsf{U}} \geq \bnu^t_{*}$. 

Last, we show that $Z_{\mathsf{U}} = Z_{\mathsf{RO}}$. Let $[\kappa^1_{\textsf{RO}}], \ldots, [\kappa^T_{\textsf{RO}}]$ be the sequence of RO assortments given by our policy, and let  $\bnu^t_{\textsf{RO}}$ be the corresponding attraction value, that is, 
\[
\nu^t_{\textsf{RO},i} := \alpha_i \bigl( \ind(i\in  [\kappa^{t-1}_{\mathsf{RO}}]) ,\dots, \ind(i\in [\kappa^{t-M}_{\mathsf{RO}}]) \bigr),
\]
where $[\kappa^s_{\mathsf{RO}}] := \emptyset$ for $s \leq 0$. In general, $\boldsymbol{\beta} \geq 0$ implies that $\bnu^t_{\mathsf{RO}} \leq \bnu^t_{\mathsf{U}}$. However, due to the nested structure on $\kappa^{1}_{\mathsf{RO}}, \ldots , \kappa^{T}_{\mathsf{RO}}$ in Lemma~\ref{eclemmanested}, for each $t \in [T]$, we obtain that $\nu^t_{\textsf{RO},k} = \nu^t_{\textsf{U},k}$ for $ k \in [\kappa^t_{\mathsf{RO}}]$ and, thus
\begin{equation}\label{eq:them1proof-1}
 \MNL([k], \bnu^t_{\mathsf{RO}}) = \MNL([k], \bnu^t_{\mathsf{U}}) \qquad \for k \in [\kappa^t_{\mathsf{RO}}].   
\end{equation}
Moreover, since $[\kappa^t_{\mathsf{RO}}]$ maximizes  $\MNL(\cdot, \bnu^t_{\mathsf{RO}})$, Part~\ref{eclemma:monotone-2} in Lemma~\ref{eclemma:monotone} shows that 
\begin{equation}\label{eq:them1proof-2}
    r_{\kappa^t_{\mathsf{RO}}} \geq \MNL([\kappa^t_{\mathsf{RO}}], \bnu^t_{\mathsf{RO}}) > r_{\kappa^t_{\mathsf{RO}}+1}.
\end{equation}
By~\eqref{eq:them1proof-1} and~\eqref{eq:them1proof-2}, we obtain $r_{\kappa^t_{\mathsf{RO}}} \geq \MNL([\kappa^t_{\mathsf{RO}}], \bnu^t_{\mathsf{U}}) > r_{\kappa^t_{\mathsf{RO}}+1}$. By invoking Part~\ref{eclemma:monotone-2} in Lemma~\ref{eclemma:monotone},  we obtain that  $ \kappa^t_{\mathsf{U}} = \kappa^t_{\mathsf{RO}}$. Therefore, we can conclude that 
\[
\MNL([\kappa^t_{\mathsf{RO}}], \bnu^t_{\mathsf{RO}}) = \MNL([\kappa^t_{\mathsf{RO}}], \bnu^t_{\mathsf{U}}) = \MNL([\kappa^t_{\mathsf{U}}], \bnu^t_{\mathsf{U}}), 
\]
showing that $Z_{\mathsf{U}} = Z_{\mathsf{RO}}$. \hfill \Halmos

\subsection{Proof of Proposition~\ref{prop:hold:condition}}\label{ec:union:ro}
\begin{repeattheorem} [Proposition~\ref{prop:hold:condition}]
Assume the absence of cross-product and cross-period constraints. Then, there exists an optimal assortment $\x = (\x^1, \ldots, \x^T)$ of~\eqref{eq:total} such that $\bigl\{i \in [N] \bigm| \sum_{t \in T}x^t_i \ge 1 \bigr\}$ is revenue-ordered.   
\end{repeattheorem}

\noindent {\bf Proof.} Let $\x = (\x^1, \ldots, \x^T)$ be an optimal solution, and let $\bnu^t = (\nu_i^t)_{i\in[N]} $ be the attraction value vector in period $t$, that is, \re{$\nu_i^t  = \alpha_i ( x_i^{t-1} ,\dots,  x_i^{t-M} )$}. For each $t \in [T]$, let $S^t$ denote $\{i \in [N] \mid x^t_i =1\}$, and define $S := \cup_{t \in [T]} S^t$. Suppose that $S$ is not revenue-ordered. We assume that products are indexed such that $r_1\ge \dots \ge r_N$, then there exists a pair of products, $i$ and $j$ with $i <j$, such that $i \notin S$ but $j \in S$. We consider two possible cases on the revenue of the $i^{\text{th}}$ product. 

Suppose that there exists $\tau \in [T]$ such that $\MNL(S^\tau, \bnu^\tau)\leq r_i$, where $\MNL(\cdot, \cdot)  $ is defined as in~\eqref{single-ast-mnl}. Let $\rev_{\textsf{HAP}}(Z^1, \ldots,Z^T)$ be the objective value of~\eqref{eq:total} if $Z^t$ is the assortment for the $t^{\text{th}}$ period. It follows readily that
\[
\begin{aligned}
\rev_{\textsf{HAP}}(S^1, \ldots, S^{\tau-1}, S^\tau \cup \{i\}, S^{\tau +1},   \ldots,S^T) &= \frac{1}{T}\biggl( \sum_{t \in [T] \setminus \{\tau\}}\MNL(S^t, \bnu^t) + \MNL(S^\tau \cup \{i\},\bnu^\tau)\biggr) \\   
& \geq \frac{1}{T} \sum_{t \in [T]}\MNL(S^t, \bnu^t)  = \rev_{\textsf{HAP}}(S^1,   \ldots,S^T),
\end{aligned}
\]
where the first equality holds since $i \notin S^t$ for all $t \in [T]$. The inequality holds because $\MNL(S^\tau \cup \{i\}, \bnu^\tau ) = \lambda r_i + (1-\lambda) \MNL(S^\tau,\bnu^\tau) \geq \MNL(S^\tau,\bnu^\tau)$, where $\lambda = \nu^\tau_i/(1+\sum_{k \in S^\tau} \nu^\tau_k + \nu^\tau_i )$. If $\MNL(S^\tau, \bnu^\tau) =  r_i$, then we obtain a new optimal solution with the same revenue as that of $S^1,\dots, S^T$. If $\MNL(S^\tau, \bnu^\tau) < r_i$, the condition is not possible because the above inequality contradicts with the optimality of $S^1,\dots, S^T$.

Now, assume that $\MNL(S^t,\v^t)>r_i$ for all $t \in [T]$. However, this is not possible as we argue next.  Let $\tau$ be the last period such that $j \in S^\tau$. Then, we have $\MNL(S^t,\bnu^t)>r_j$ for all $t\in [T]$. We obtain a contradiction to the optimality of $S^1, \ldots, S^T$ as follows.
\[
\begin{aligned}
\rev_{\textsf{HAP}}(S^1, \ldots, S^{\tau-1}, S^\tau \setminus \{j\}, S^{\tau +1},   \ldots,S^T) &=    \frac{1}{T}\biggl( \sum_{t \in [T] \setminus \{\tau\}}\MNL(S^t, \bnu^t) + \MNL(S^\tau \setminus \{j\},\bnu^\tau)\biggr) \\
& >  \frac{1}{T} \sum_{t \in [T]}\MNL(S^t, \bnu^t)  = \rev_{\textsf{HAP}}(S^1,   \ldots,S^T),
\end{aligned}
\]
where the first equality holds due to $j \notin S^t$ for $t \geq \tau+1$. The inequality holds because $\MNL(S^\tau, \bnu^\tau)  = \lambda r_j + (1-\lambda) \MNL(S^\tau \setminus\{j\},  \bnu^\tau)< \lambda \MNL(S^\tau, \bnu^\tau) +(1-\lambda)\MNL(S^\tau \setminus\{j\},  \bnu^\tau)$, where $\lambda = \nu^\tau_j/(1+\sum_{k \in S^\tau}\nu^\tau_k)$. The above inequality contradicts with the optimality of $S^1,\dots, S^T$. \hfill \Halmos

\section{Proofs of Section \ref{sec:two:cycle}}
 We derive a graph representation of~\eqref{eq:total} for proving cyclic policies. With each instance of~\eqref{eq:total}, we associate a directed graph $D = (V,A)$ with nodes $V$ and arcs $A$, referred to as \textit{assortment graph}, as follows. A node in the assortment graph is a tuple of historical assortments, $(\x^{t-1},\dots,\x^{t-M})$. That is, a node $v = (v^1,\dots,v^M)$ in $V$ has $M$ elements, and each one is an assortment. The first element denotes the latest offered assortment, and the last one is the earliest one in memory. The number of nodes is $2^{N M}$. We say a node $\mu = (\mu^1, \ldots, \mu^M)$ is a predecessor of $v = (v^1, \ldots, v^M)$ if $\mu^{k} = v^{k+1}$ for $k = 1, 2, \ldots, M-1$. Then, an ordered pair $(\mu,v)$ is an arc going from $\mu$ to $v$ if $\mu$ is a predecessor of $v$. Therefore, the arc set in the assortment graph is given as follows:
\[
A := \{(\mu,v) \mid   \text{ if } \mu \text{ is predecessor of } v  \text{ and }\mu, v \in V \}.
\]
The number of arcs is $2^{NM+N}$. The arc $(\mu,v)$ records the latest assortment $v^1$ and drops the earliest one $\mu^M$, which defines the updating process of history assortments within memory. Then, we can define an arc's weight as 
\[
w(\mu,v) = \sumi r_i \pi_{i}(v^1, \mu^1, \ldots, \mu^{M})
\]
where $\pi_{i}(\cdot)$ is the purchase probability of $i$ given current assortment $v^1$ and historical assortments $(\mu^1, \ldots, \mu^{M})$ defined in Section~\ref{sec:sub:model}. Here, we omit the superscript of period $t$ for simplicity. 
% To highlight the dependence on the predecessor, we define $R(v^1 , \mu ) = \sumi r_i \pi_{i}(v^1, \mu^1, \ldots, \mu^{M})$. Hence, $w(\mu, v) = R(v^1,\mu )$. 

\subsection{Proof of Theorem~\ref{thm:general:cyclic}}
\begin{repeattheorem}[Theorem~\ref{thm:general:cyclic}]
  Given an instance of~\eqref{eq:asymptotic}, there exists a positive integer $L^*$ such that $L^*$-cyclic policy is optimal to~\eqref{eq:asymptotic}.
\end{repeattheorem}

\noindent {\bf Proof.} We show that finding an optimal assortment planning for~\eqref{eq:asymptotic} equals finding the maximum mean cycle in the assortment graph, which is the cycle with the largest mean weight. 

% This implies that Since the maximum mean cycle of a graph with non-negative weights of arcs exists \citep{karp1978characterization}, ~\eqref{eq:asymptotic}.   

First, we show an assortment planning of~\eqref{eq:asymptotic} maps to a path in the assortment graph. Suppose $\x^1,\x^2,\dots$ is an assortment planning of~\eqref{eq:asymptotic}, it defines a path in the assortment graph, that is, $v_{\boldsymbol{0}},v(1),v(2),\dots$, where $v_{\boldsymbol{0}} = (\boldsymbol{0},\dots, \boldsymbol{0})$ and 
\[
v(t)^k  =\begin{cases}
     \x^{t-k+1} & \text{ if }t-k \ge 0 \\
      \boldsymbol{0} & \text{ otherwise},
\end{cases} \quad \text{ for } k\in [M] \text{ and } t\in \Z_{+}.
\]
% \begin{align*}
%      v(t)^k  & = \bigl \{ \begin{array}{cc}
%        \x^{t-k+1}   &  t-k+1\ge 1 \\
%        \boldsymbol{0}   & t-k+1 < 1
%      \end{array}\bigr. \text{ for } k\in [M], t\in \mathbb{N}_{+}.
% \end{align*} 
It suggests $v(t)$'s first element records the assortment decision in period $t$, and the following $M-1$ elements record the historical assortments before $\x^t$. Then, the mean weight of the path equals to the average revenue of the assortment planning. That is,  
 \[  \lim_{T\to \infty} \frac{1}{T} \sum_{t\in[T]} w\bigl(v(t-1),v(t)\bigr) = \lim_{T\to \infty}  \frac{1}{T} \sum_{t \in [T]} \sum_{i \in [N]} r_i\pi_i^t(\x^t,\x^{t-1}\dots, \x^{t-M}). \]  
 
Next, we show that to solve~\eqref{eq:asymptotic}, it suffices to find the maximum mean cycle denoted as $C^*$. Since any assortment planning is a path in the assortment graph, solving~\eqref{eq:asymptotic} equals finding an infinite-length path with the largest mean weight and the starting node $v_{\boldsymbol{0}}$. We call such a path an optimal path. Because the assortment graph has finite nodes and arcs, an optimal path will inevitably repeat cycles in the assortment graph. Given an infinite-length path, we can repeatedly replace each cycle in the path with the maximum mean cycle. Finally, we obtain a path that repeats the maximum mean cycle forever and has a larger mean weight than the initial one. We do not need to consider the path's weight from the starting node $v_{\boldsymbol{0}}$ to the maximum mean cycle because it will vanish when taking the average over an infinite horizon $T$. 
% Since the assortment graph has finite cycles and non-negative weights on all arcs, we can definitely find the maximum mean cycle. 

Next, we show how to recover an assortment planning from the maximum mean cycle. Suppose the maximum mean cycle is $C^* = v(1),v(2),\dots,v(L^*-1), v(0)$ where $L^*$ is the length of the maximum mean cycle. We can construct a cyclic policy as follows.
\[ \x^t = v^1( t\text{ mod } L^*) \for t \in \mathbb{N}_+. \]
It suggests $\x^t$ is the first element of $v( t\text{ mod } L^*)$. Therefore, we prove that an optimal solution of~\eqref{eq:asymptotic} is a cyclic policy. The average revenue generated by the cyclic policy is identical to the mean weight of the maximum mean cycle $C^*$, and the length of the cyclic policy is $L^*$. If we know $L^*$, we can compute the cyclic assortment planning via~\eqref{eq:cyclic-policy-conic}. \hfill \Halmos

\subsection{Proof of Theorem~\ref{thm:nonoverlap:condition} }
\begin{repeattheorem}[Theorem~\ref{thm:nonoverlap:condition}]
 Assume that $M < N$. An $(M+1)$-cyclic policy is optimal for problem~\eqref{eq:asymptotic} if the cross-period constraint is non-overlapping. 
\end{repeattheorem}

\noindent {\bf Proof.} We prove this result by showing any $L$-cyclic policy with length $L>M+1$ will be dominated by the $(M+1)$-cyclic policy in terms of the expected revenue. Recall that we do not consider cyclic polices with length $L<M+1$ because they contain empty sets.

Given an $L$-cyclic policy satisfying the non-overlapping condition, $\x^1,\dots,\x^L$, the purchase probability of product $i$ in period $t$ is 
\[\pi_i^t(\x^t,\x^{t-1},\dots,\x^{t-M}) = \frac{x_i^t \exp(\beta_i^0) }{ 1+ \sumi x_i^t \exp(\beta_i^0) },\]
because $x_i^t + x_i^{t-1} + \dots + x_i^{t-M} \le 1 $ for $i\in [N]$ and $t\in [L]$. Hence, we can simplify the average revenue of such a cyclic policy as follows 
\[ \rev_L  =  \frac{1}{L}\sum_{t \in [L]} R(\x^t) =  \frac{1}{L} \sum_{t \in [L]} \sumi r_i x_i^t \exp(\beta_i^0)/(1+ \sumi x_i^t \exp(\beta_i^0) ) . \]

For any $L$-cyclic policy with length $L>M+1$, we show its average revenue is a convex combination of cyclic policies with length $M+1$. We rewrite $\rev_L$ as 
\begin{align*}
    \rev_L  
     & =\frac{ (M+1)R(\x^1) + \dots +(M+1)R(\x^L) }{(M+1)L} \\
     & = \frac{ \bigl(R(\x^1)+R(\x^2)+ \dots + R(\x^{M+1}) \bigr) +    \dots  +\bigl(R(\x^L)+ R(\x^1) +\dots + R(\x^{M}) \bigr)   }{(M+1)L} \\
     & = \frac{1}{L} \frac{R(\x^1)+R(\x^2)+\dots + R(\x^{M+1}) }{M+1} +  \dots + \frac{1}{L} \frac{R(\x^L)+ R(\x^1)+ \dots + R(\x^{M}) }{M+1} \\
     & \le \rev^*_{M+1}.
\end{align*}
The first equation holds by multiplying both the denominator and numerator by $M+1$. The second equation holds by re-arranging orders of $R(\x^1),\dots, R(\x^L)$. The third equation shows that we construct $L$ cyclic policies with length $M+1$ by abstracting $M+1$ adjacent assortments in the cycle $\x^1,\dots, \x^L$. Each cyclic policy is a feasible solution of~\eqref{eq:cyclic-policy-conic} with length $M+1$ since they are non-overlapping. Then, we have the last inequality because $\rev^*_{M+1}$ is generated by the optimal solution of~\eqref{eq:cyclic-policy-conic} with length $M+1$. Therefore, the $(M+1)$-cyclic policy has the largest revenue under the non-overlapping condition.

\subsection{\re{Proof of Theorem~\ref{them:bound-free}}}
\begin{repeattheorem}[Theorem~\ref{them:bound-free}]
\re{$Z_{\textsc{Bound-free}} \le Z_{\textsc{Cycle-Conic}}$ if the cross-period constraint is non-overlapping.}  
\end{repeattheorem}

\noindent {\bf Proof.} 
In the following, we will show that $Z_{\textsc{Bound-free}} \le  Z_{\textsc{McCormick}} \le Z_{\textsc{Cycle-Conic}}$, where $Z_{\textsc{McCormick}}$ is the optimal objective value of the continuous relaxation of the following integer programming:
\re{\begin{subequations} \label{ec:mccormick}
    \begin{alignat}{3}
        \max \quad & \frac{1}{M+1}  \sum_{t \in [M+1]} \sumi  r_i  u_i \gamma_{i}^t  \notag \\ 
        \st \quad & x_i^t + \summ x^{\tau(m \mid t)}_i  \leq 1   && \for i \in [N] \label{ec:same-1} \\
               &  \rho^t + \sum_{i \in [N]}u_i \gamma^t_i = 1 && \for t\in [M+1]  \label{ec:same-2} \\ 
   & \gamma^t_i  \leq  \rho^t_Ux^t_i &&\for i \in [N] \text{ and } t\in [M+1] \label{ec:mccormick-1}\\
   &  \gamma_i^t  \leq   \rho^t+\rho^t_Lx^t_i -\rho^t_L &&\for i \in [N] \text{ and } t\in [M+1]\label{ec:mccormick-2} \\
   &  \gamma^t_i \geq  \rho^t + \rho^t_Ux_i^t - \rho^t_U &&\for i \in [N] \text{ and } t\in [M+1]\label{ec:mccormick-3}\\
   &  \gamma^t_i  \geq  \rho^t_Lx^t_i &&\for i \in [N] \text{ and } t\in [M+1] \label{ec:mccormick-4} \\ 
  &  \x^t \in \mathcal{X} \cap \{0,1\}^{N }  && \for t\in [M+1], 
    \end{alignat}
\end{subequations}}where constraints \eqref{ec:mccormick-1}-\eqref{ec:mccormick-4} are exactly \eqref{eq:relax1}, obtained by using the McCormick envelopes to linearize $\gamma_i^t = \rho^t x_i^t$.

First, we prove $Z_{\textsc{Bound-free}} \le Z_{\textsc{McCormick}}$. Since the bound-free formulation~\eqref{eq:2cycMILP-boundfree} and formulation~\eqref{ec:mccormick} have identical objectives, it suffices to show that any feasible solution of the continuous relaxation of~\eqref{eq:2cycMILP-boundfree} is also feasible to that of~\eqref{ec:mccormick}. More specifically, we need to show that constraints \eqref{ec:mccormick-1}-\eqref{ec:mccormick-4} are implied by~\eqref{eq:perfect-1}-\eqref{eq:perfect-4} of the bound-free formulation. Constraints~\eqref{eq:perfect-2}, $\Gamma^t_{ij} \leq \gamma_i^t$ for $i\neq j$ from~\eqref{eq:perfect-3}, and $\Gamma^t_{ii} =  \gamma_i^t$ from~\eqref{eq:perfect-4} imply that 
\[
x^t_i \leq \biggl(1+\sum_{j \in N}u_j\biggr) \gamma_i^t,
\]
which is~\eqref{ec:mccormick-4}.  Constraints~\eqref{eq:perfect-1}~\eqref{eq:perfect-2}, $\Gamma^t_{ij} \leq \gamma_j^t$ for $i\neq j$ from~\eqref{eq:perfect-3}, and $\Gamma^t_{ii} = \gamma_i^t$ from~\eqref{eq:perfect-4} imply that 
\[
x^t_i \leq \gamma^t_i + \sum_{j \in [N]}u_j \gamma^t_j = \gamma^t_i + 1 - \rho^t,
\]
which is~\eqref{ec:mccormick-3}. Also, constraint~\eqref{eq:perfect-2}, $\Gamma^t_{ij} \geq 0 $ for $i \neq j$ and $\Gamma^t_{ii} = \gamma^t_i$ imply
\[
x^t_i \geq \gamma^t_i, 
\]
which is~\eqref{ec:mccormick-1}. Constraint~\eqref{eq:perfect-1}~\eqref{eq:perfect-2}, $\Gamma^t_{ij} \geq \gamma_i^t + \gamma_j^t -\rho^t$ for $i\neq j$, $\Gamma_{ii}^t = \gamma_i^t$, and $\gamma_i^t \le \rho^t$ imply
\[
\begin{aligned}
    x^t_i&\geq \gamma^t_i + u_i \gamma^t_i+ \sum_{j \neq i}u_j(\gamma^t_i+\gamma^t_j-\rho^t) =(1+\sum_{j \in [N]}u_j)(\gamma^t_i - \rho^t)+ (1+u^t_i)\rho^t +\sum_{j\neq i}u_j\gamma^t_j \\
    & = (1+\sum_{j \in [N]}u_j)(\gamma^t_i - \rho^t)+ u^t_i\rho^t + 1-u^t_i \gamma^t_i \geq(1+\sum_{j \in [N]}u_j)(\gamma^t_i - \rho^t)+1,
\end{aligned}
\]
which is \eqref{ec:mccormick-2}.

% \begin{align*}\label{eq:middle}
%     \max \quad & \sum_{t =1,2} \frac{1}{2} \sumi  r_i  u_i \gamma_{i}^t  \tag{\textsc{Middle}} \\ 
%     \st \quad & x^1_i+x^2_i \leq 1   && \for i \in [N]  \\
%            &  \rho^t + \sum_{i \in [N]}u_i \gamma^t_i = 1 && \for t=1,2 \\
%           &  \eqref{eq:relax1} \text{ and }  \x^t \in \mathcal{X} \cap \{0,1\}^{N }  && \for  t=1,2 , 
% \end{align*}
% where \eqref{eq:relax1} is the McCormick envelopes for linearizing $\gamma_i^t = \rho^t x_i^t$.

Next, we prove $ Z_{\textsc{McCormick}} \le Z_{\textsc{Cycle-Conic}}$. Note that in this case, \eqref{eq:cyclic-policy-conic} becomes  
\re{\begin{subequations}
    \begin{alignat}{3}
         \max  \quad & \frac{1}{M+1}  \sumi \sum_{t\in [M+1]} r_iy_{i}^t  \notag  \\
    \st \quad & x_i^t + \summ x^{\tau(m \mid t)}_i  \leq 1   && \for i \in [N] \notag \\
    &\rho^t + \sumi y_{i}^t = 1   && \for\, t \in [M+1]  \notag \\ 
    &\max \bigl\{ 0,\ \rho^t_U x^{\tau(m|t)}_i + \gamma^t_i - \rho^t_U \bigr\} \leq z^m_{it} \leq  \min \bigl\{ \gamma_{i}^t , \rho^t_U x^{\tau(m|t)}_i  \bigr\}   && \for\, t \in [M+1], m\in [M] \label{conic:mc-1}\\  
   & ~\eqref{eq:relax1},~\eqref{eq:conv-pers},~\eqref{eq:modify:h}    && \for t\in [M+1], i\in [N] \notag\\
   &  \x^t \in \mathcal{X} \cap \{0,1\}^{N }  && \for  t\in [M+1], \notag
    \end{alignat}
\end{subequations} }where~\eqref{eq:relax1} and~\eqref{conic:mc-1} are McCormick envelopes linearizing $\gamma_i^t = \rho^t x_i^t$ and $z_{im}^t = \gamma_i^t x_{i}^{\tau(m|t)}$, respectively. ~\eqref{eq:conv-pers} and~\eqref{eq:modify:h} are convex and concave extensions of $y_i^t$, respectively. Letting $(\x, \boldsymbol{\rho}, \boldsymbol{\gamma})$ be a feasible solution of the continuous relaxation of~\eqref{ec:mccormick}, and define $y^t_i = \gamma^t_iu_i$ and $z^t_{im} = 0$ for all \re{$t \in [M+1]$} and $i \in [N]$. To complete the proof, we only need to show that $(\x, \boldsymbol{\rho}, \boldsymbol{\gamma},\boldsymbol{y},\boldsymbol{z})$ is feasible to~\eqref{eq:cyclic-policy-conic}. Clearly, constraints~\eqref{eq:relax1}, ~\eqref{eq:conv-pers} and~\eqref{eq:modify:h} are satisfied. \re{Next, we verify constraint~\eqref{conic:mc-1} is satisfied. Clearly, $\min \bigl\{ \gamma_{i}^t , \rho^t_U x^{\tau(m|t)}_i  \bigr\} \ge 0 = z^t_{im} $. On the other hand, 
\begin{align*}
     \rho^t_U x^{\tau(m|t)}_i + \gamma^t_i - \rho^t_U & \le  \rho^t_U x^{\tau(m|t)}_i + x_i^t \rho^t_U - \rho^t_U  \\
    & =  \rho^t_U (x^{\tau(m|t)}_i+ x_i^t -1) \\
    & \le  0 =z^t_{im},
\end{align*} 
where the first inequality holds because of \eqref{ec:mccormick-1} and the second inequality holds because of~\eqref{ec:same-1}.}\Halmos

\subsection{Example of the Unstable Maximum Mean Cycle Length}\label{eg:cycle:change}

The following example shows that a slight fluctuation in one parameter could result in a different maximum mean cycle length if there are no non-overlapping constraints.
\begin{example} 
Consider a case where a firm sells three products to customers with a memory length of one. There are no cross-period or cross-product constraints. The outside option has utility 0. Table~\ref{tab:example} provides a detailed parameter setup.
\begin{table}[hbtp]
    \centering
    \caption{Product parameters}
    \label{tab:example}
    \begin{tabular}{cccc}\hline 
        Product & Price & Base utility  & One-period history-dependent effect \\ \hline 
       1  & 6 & 1.3 & -0.5\\
       2  & 16 & -0.6 & -0.8\\
       3  & 17 & 0.2 & -0.9\\ \hline 
    \end{tabular}
\end{table} 
The two-cyclic policy has larger revenue than that of the three-cyclic policy. However, once the base utility of product $3$ increases from $0.2$ to $0.3$, the average revenue of the two-cyclic policy is smaller than that of the three-cyclic policy. \hfill   \Halmos
\end{example}

\section{Additional Material to Section \ref{sec:numerical}}

% \subsection{\re{Convex Envelope of the Attraction Value Function}}\label{ec:sec:env}

% \subsection{\re{Derive Bounds of No-Purchase Probabilities}}\label{ec:sec:rho:bound}
% \re{In all our numerical studies, we use the same global upper bound and lower bound of the no-purchase probability variable $\rho^t$ for each $t\in [T]$. We use the minimum (resp. maximum ) attraction value of each product $i$ to compute the global upper bound $\rho^t_U$ (resp. lower bound $\rho^t_L$). We have:

\subsection{Multilinear Extension Based Formulation}\label{sub:multi:detail}
In this subsection, we introduce the multilinear extension based formulation of~\eqref{eq:total}. For the nonlinear structure~\eqref{eq:choice}, we first use the multilinear extension to present the attraction value function $\alpha_i(\cdot)$ \citep{o2014analysis}, and then recursively apply McCormick envelopes to derive a linear representation. It is well-understood that both the quality and the size of the resulting linear representations depend on the recursive sequence, and finding an optimal recursive sequence amounts to solving a difficult combinatorial optimization problem~\citep{speakman2017quantifying,khajavirad2023strength}. \cite{speakman2017quantifying} show that the relaxation quality depends on the sequence of recursion, but we do not exploit the best sequence in this paper. Even though the number of introduced variables and constraints increases exponentially with the memory length, this formulation is solvable under a small $M$.

We represent the attraction value function via a sum of multiple multilinear functions. Given memory length $M\in \Z_{+}$, $\{S_k\mid S_k\subseteq [M],k\in [2^M]\}$ is the collection of all subsets of $[M]$. For every $S_k,k\in[2^M]$, let $\chi^{S_k}$ be its indicator vector, that is, the $m^{\text{th}}$ coordinate of $\chi^{S_k}$ is $1$ if and only if $m\in S_k$. We can rewrite the attraction value function as 
\[ \alpha_i(x_i^{t-1},\ldots,x_i^{t-M}) = \sum_{k\in[2^M]} \bigl( \Pi_{m\in S_k} x_i^{t-m}\bigr) \cdot a_{ik} \quad \for \, i\in [N], \]
where $a_{ik} = \alpha_i(\chi^{S_k})- \sum_{j:S_j \subset S_K} a_{ij}$. 
% For instance, when $M=2$, 
% \[\alpha_i(x_i^{t-1},x_i^{t-2}) = a_{i1} + x_i^{t-1}a_{i2} + x_i^{t-2}a_{i3} + x_i^{t-1}x_i^{t-2}a_{i4}\quad \for \, i\in [N], \]
% where 
% \begin{align*}
%     a_{i1} &= \alpha_i(0,0) \\
%     a_{i2} &= \alpha_i(1,0)-\alpha_i(0,0) \\
%     a_{i3} & = \alpha_i(0,1)-\alpha_i(0,0) \\
%     a_{i4} &= \alpha_i(1,1)+ \alpha_i(0,0) - \alpha_i(1,0) - \alpha_i(0,1).
% \end{align*}

Then, introduce $\theta_{ik}^t$ to present the multilinear term, that is, 
\[\theta_{ik}^t = x_i^t \cdot \bigl( \Pi_{m\in S_k} x_i^{t-m}\bigr) = \Pi_{m\in S_k \cup \{0\}} x_i^{t-m} \quad \for \,i\in [N], t\in[T], k\in[2^M]. \]
We recursively use McCormick envelopes to linearize $\theta_{ik}^t$. The idea is to sequentially introduce artificial variables and reduce the number of variables in the multilinear term. First, select any two components $p,q\in S_k\cup \{0\}$. Then, introduce a variable $h_i^t(p,q)$ that corresponds to the bilinear product $x_i^{t-p} x_i^{t-q}$ and the multilinear term changes to 
\[ \theta_{ik}^t = h_i^t(p,q) \cdot \Pi_{m\in S_k \cup \{0\}\setminus\{p,q\}} x_i^{t-m}.\]
We can use McCormick envelopes to linearize $h_i^t(p,q)=x_i^{t-p} x_i^{t-q}$,
\begin{align*}
    \max\{0, x_i^{t-p} + x_i^{t-q} -1 \} \le h_i^t(p,q) \le \min\{x_i^{t-p},x_i^{t-q}\}.
\end{align*}
This procedure can be recursively applied to the remaining parts until $\theta_{ik}^t$ is completely linearized. We use a notation $\mathcal{H}(i,t,k)$ to present the corresponding McCormick envelopes. 

Note that $\theta_{i,k}^t \in\{0,1\}$, we also apply the McCormick envelopes to linearize  $\rho^t \theta_{ik}^t$ denoted by $w_{ik}^t$. That is,
\begin{subequations}\label{model:multi:mc}
    \begin{alignat}{3}
         w_{ik}^t &  \leq \min \bigl\{\rho^t_L \theta_{ik}^t + \rho^t-\rho^t_L ,\ \rho^t_U \theta_{ik}^t  \bigr\}   \\
 w_{ik}^t & \geq \max \bigl\{ \rho^t_L \theta_{ik}^t,\ \rho^t_U \theta_{ik}^t + \rho^t - \rho^t_U \bigr\}  .
    \end{alignat}
\end{subequations}
Finally, we obtain a mixed-integer linear formulation of~\eqref{eq:total} as follows. 
\begin{align*}  
     \max \quad & \frac{1}{T}\sumt 
 \sumi r_iy_{i}^t \notag \\ \tag{\textsc{Multilinear}} 
        \st \quad &  \x^t \in \{0,1\}^{N } \cap \mathcal{X} \text{ and }  \x \in \mathcal{P} && \for \, t \in [T]  \notag \\
    &\rho^t + \sumi y_{i}^t = 1  \quad && \for\, t \in [T] \notag \\
    & y_i^t = \sum_{k\in[2^M]} w_{ik}^t a_{ik} \quad  && \for \, t \in [T] \text{ and } i \in [N] \\
    & \theta_{ik}^t \in \mathcal{H}(i,t,k) \quad  && \for \, t \in [T] , i \in [N] \text{ and } k \in [2^M] \\
    &\eqref{model:multi:mc}&& \for \, t \in [T] , i \in [N] \text{ and } k \in [2^M] .
\end{align*}
 
\subsection{\re{A Projected Cutting-Plane Implementation for Bound-Free Formulation~\eqref{eq:2cycMILP-boundfree}}}\label{sub:sepa:oracle}
\re{The bound-free formulation~\eqref{eq:2cycMILP-boundfree} has $\mathcal{O}( N^2)$ variables due to the presence of $\boldsymbol{\Gamma}$ variables. We design a cutting-plane procedure to implement~\eqref{eq:2cycMILP-boundfree} for computation efficiency. Recall that~\eqref{eq:2cycMILP-boundfree} is reformulation of~\eqref{eq:2cycMNL}. In our implementation, we first reformulate~\eqref{eq:2cycMNL} as a base model. Next, we use a projected separation oracle to cut off infeasible solutions for $K$ rounds.}

\re{We linearize~\eqref{eq:2cycMNL} to a base model as follows.} First, following the idea in Section~\ref{sec:reformulation}, we introduce variable $\rho^t$ to denote $1/({1+ \sumj u_j x^t_j})$ with constraint $\rho^t + \sumj u_j x^t_j \rho^t =1$. Next, let $\gamma_j^t$ denote the bilinear product $x^t_j \rho^t$ and apply McCormick envelopes to linearize the bilinear term. We also notice that, for $t \in [M+1]$,
\begin{equation*}
    \rho^t (1+\sumj x_j^t u_j) \ge 1,
\end{equation*} 
which is a valid convex constraint and is representable via the second-order cone. Let $w^t = 1+\sumj x_j^t u_j $ for $t \in [M+1]$ and add the above constraint into our formulation. Finally, we obtain a mixed-integer second-order conic formulation of~\eqref{eq:2cycMNL} as follows.
\re{\begin{equation}\label{eq:2cycMNL:linear}
\begin{aligned}
    \max \quad & \frac{1}{M+1}  \sum_{t\in [M+1]}  \sumi  r_i u_i \gamma_i^t \\
    \st \quad  & x_i^t + \summ x^{\tau(m \mid t)}_i  \leq 1 &&   \for i \in [N] \\
    & \rho^t + \sumi u_i \gamma_i^t =1 &&   \for t \in [M+1]\\
    & \gamma_{i}^t   \leq \min \bigl\{\rho^t_L x_{i}^t + \rho^t-\rho^t_L ,\ \rho^t_U x_{i}^t  \bigr\}  && \for i \in [N], t\in [M+1]  \\
 & \gamma_{i}^t \geq \max \bigl\{ \rho^t_L x_{i}^t,\ \rho^t_U x_{i}^t + \rho^t - \rho^t_U \bigr\}&& \for i \in [N], t\in [M+1] \\ 
 & \rho^t + w^t \ge ||\rho^t - w^t, 2||&& \for t\in [M+1]  \\
 & w^t = 1+\sumj x_j^t u_j && \for t\in [M+1] \\ 
    &\x^t \in \{0,1\}^{N} &&   \for t \in [M+1],
\end{aligned}\tag{\textsc{Base}}
\end{equation}}
where $\rho^t_L$ and $\rho^t_U$ are lower and upper bound of $\rho^t$ for $t\in [M+1]$, respectively. Here, although we need the bound of $\rho^t$ when implementing \eqref{eq:2cycMILP-boundfree}, we do not require it to be very tight, and we use it only for getting an equivalent formulation of~\eqref{eq:2cycMNL}.  

\re{We tighten the base formulation through a projected separation oracle (Algorithm~\ref{alg:sepa:oracle}). Given a solution $(\hat{\x},\hat{\boldsymbol{\gamma}},\hat{\boldsymbol{\rho}})$, the separation oracle checks whether the projection of~\eqref{eq:perfect-2} in the space of the base formulation is violated for each $i\in[N]$ and $t\in[M+1]$. More specifically, we replace $\Gamma_{ij}^t$ in \eqref{eq:perfect-2} with the bounds in~\eqref{eq:perfect-3} and~\eqref{eq:perfect-4}, then we have  
\begin{equation}\label{eq:bound:xit}
    \gamma_i^t(1+u_i)  + \sum_{j\in [N],j\neq i} u_j\min\{0, \gamma_i^t + \gamma_j^t -\rho^t \}\le  x_i^t \le \gamma_i^t + \sum_{j\in [N]} u_j \min\{\gamma_i^t, \gamma_j^t\}. 
\end{equation}
If the current solution violates \eqref{eq:bound:xit}, we add the corresponding constraint to the base model. We detail each step in Algorithm~\ref{alg:sepa:oracle}. }

\begin{algorithm} [htbp]
    \caption{Projection-based separation oracle for~\eqref{eq:2cycMNL:linear}}\label{alg:sepa:oracle}
    \begin{algorithmic}[1]
        \STATE Input:  $(\hat{\x},\hat{\boldsymbol{\gamma}},\hat{\boldsymbol{\rho}})$ and base utility $\boldsymbol{u}$ 
         \STATE  Cuts $\gets \emptyset$, $ A\gets \emptyset$, $B\gets \emptyset$, and $C\gets \emptyset$
        \FOR{$i\in [N]$ and $t\in [M+1]$} 
        \FOR{$j \in [N]$}
        \IF{ $\hat{\gamma}_i^t + \hat{\gamma}_j^t \ge \hat{\rho}_t$ and $i\neq j$}
        \STATE{ $A\gets A \cup \{j\}$}
        \ENDIF 
         \IF{ $\hat{\gamma}_j^t \ge \hat{\gamma}_i^t $}
        \STATE $B\gets B \cup \{j\}$
        \ELSE 
        \STATE $C\gets C \cup \{j\}$
        \ENDIF 
        \ENDFOR 
        \STATE lower-bound $\gets (1+u_i)\cdot \hat{\gamma}_i^t + \sum_{j\in A}u_j(\hat{\gamma}_i^t + \hat{\gamma}_j^t - \hat{\rho}_t)$  
        \STATE upper-bound $\gets \hat{\gamma}_i^t + \sum_{j\in B} u_j\hat{\gamma}_i^t +  \sum_{j\in C}u_j\hat{\gamma}_j^t$
        \IF{ $\hat{x}_i^t \le \text{lower-bound} $ }
        \STATE Cuts $\gets x_i^t \ge (1+u_i)\gamma_i^t  + \sum_{j\in A}u_j( {\gamma}_i^t +  {\gamma}_j^t -  {\rho}_t)$
        \ELSIF{$\hat{x}_i^t \ge \text{upper-bound}  $ }
        \STATE  Cuts $\gets x_i^t \le \gamma_i^t + \sum_{j\in B} u_j \gamma_i^t + \sum_{j\in C} u_j \gamma_j^t $ 
        \ENDIF 
        \ENDFOR
        \STATE Output: Cuts
    \end{algorithmic}
\end{algorithm} 

\re{Running the separation oracle once is a round of cut generation. We add new cuts for K rounds and the final formulation as BF-K, where BF denotes “bound-free”. In each round $k\in [K]$, the cutting-plane algorithm solves the continuous relaxation of current formulation {BF}-$(k-1)$ and obtain an optimal solution $(\hat{\x},\hat{\boldsymbol{\gamma}},\hat{\boldsymbol{\rho}})$. Then, the algorithm calls the separation oracle (Algorithm~\ref{alg:sepa:oracle}) to generate new cuts. The algorithm finally solves {BF}-$K$ with binary constraints and obtains an optimal integer solution. We present details of this implementation process in Algorithm~\ref{alg:cutting}.}  
\begin{algorithm}[hbtp]
    \caption{A projected-cutting-plane implementation of bound-free formulation~\eqref{eq:2cycMILP-boundfree}}\label{alg:cutting}
    \begin{algorithmic}[1]
        \STATE Input:  Formulation~\eqref{eq:2cycMNL:linear} (denoted as {BF}-$0$) and a positive integer $K$ 
        \STATE $k\gets 1$
      \WHILE{$k \le  K$}
      \STATE $(\hat{\x},\hat{\boldsymbol{\gamma}},\hat{\boldsymbol{\rho}}) \gets$ an optimal solution to the continuous relaxation of {BF}-$(k-1)$
      \STATE Cuts $\gets $ Call the separation oracle, Algorithm~\ref{alg:sepa:oracle}, to separate$(\hat{\x},\hat{\boldsymbol{\gamma}},\hat{\boldsymbol{\rho}})$ 
        \STATE {BF}-$k$ $\gets$ Add Cuts into {BF}-$(k-1)$formulation 
        \STATE $k\gets k+1$
        \ENDWHILE
        \STATE Solve {BF}-$K$ with binary constraints of $\x$
        \STATE Output: An optimal solution of {BF}-$K$, $(\x^*,\boldsymbol{\gamma}^*,\boldsymbol{\rho}^*)$
    \end{algorithmic}
\end{algorithm}

\subsection{\re{A Cutting-Plane Implementation for Large Memory Length}}\label{ec:sec:double}
\re{One of the main difficulties of solving~\eqref{eq:Conic} is that for each $i\in [N]$ and $t\in [T]$, the number of linear inequalities used to describe the concave envelope of the attraction value function in~\eqref{eq:pers-extension} is $M!$, which grows exponentially as $M$ becomes large. In the following, we will implement constraints in~\eqref{eq:pers-extension} in the so-called lazy fashion. More specifically, we use the \texttt{Callback} routine of \texttt{Gurobi} to implement~\eqref{eq:pers-extension} as lazy constraints. Constraints in~\eqref{eq:pers-extension} are removed from~\eqref{eq:Conic} and placed in a lazy constraint pool. Then, \texttt{Gurobi} solves the relaxed problem and checks whether inequalities in the pool are violated at each integer solution generated in the branch-and-bound tree. If a violated inequality is found, it is added to the node, the integer solution is cut off, and the node is resolved. We elaborate on our implementation as follows.} 
% \re{We design a cutting-plane procedure to implement~\eqref{eq:Conic} with large memory length. In our implementation, we first reformulate \eqref{eq:Conic} to an MILP formulation to utilize the efficient Callback procedure of Gurobi for MILP. Next, we solve a sequence of restricted models in a constraint generation loop through a separation oracle to add new cuts. We apply the separation oracle for both the continuous relaxation of the restricted model and its integral model.}

\re{Since \texttt{Gurobi} does not support the exponential function, we linearize the exponential conic constraint~\eqref{eq:conv-pers} using the subgradient inequalities of the continuous relaxation of $\alpha_i(\cdot)$ and treat the subgradient inequalities as a lazy constraint as well. Note that the continuous relaxation $\tilde{\alpha}_i(\cdot)$ is a closed convex function over $[0,1]^M$. Thus, by~\cite{rockafellar1997convex}, for $(x_i^{t-1},\dots, x_i^{t-M}) \in [0,1]^M$
\begin{equation*} 
  \tilde{\alpha}_i(x_i^{t-1},\dots, x_i^{t-M}) = \max_{\w}\biggl\{ \tilde{\alpha}_i({\boldsymbol{w}} ) + \summ (x_i^{t-m} - {w}^m) \beta_i^m \tilde{\alpha}_i({\boldsymbol{w}}) \biggm| {\w} \in [0,1]^M \biggr\}. 
\end{equation*}
where $\beta_i^m  \tilde{\alpha}_i({\boldsymbol{w}})$ is the $m^{\text{th}}$ partial derivate of $\nabla \tilde{\alpha}_i({\boldsymbol{w}})$. Therefore, using the variable $\gamma_i^t$ to scale the subgradient inequalities, the exponential conic constraint~\eqref{eq:conv-pers} is equivalent to 
\begin{equation}\label{eq:conv-env:grad}
     y_i^t \ge \gamma_i^t \alpha_i({\boldsymbol{w}}) + \summ (z_{im}^{t} - \gamma_i^t {w}^m) \beta_i^m   \alpha_i({\boldsymbol{w}}) \qquad \for {\boldsymbol{w}}\in [0,1]^M.
\end{equation}
}

\re{To implement~\eqref{eq:pers-extension} and~\eqref{eq:conv-env:grad} as lazy constraints in \texttt{Gurobi}, we use Algorithm~\ref{alg:sepa:lovasz} to solve the separation problem. Given a point (solution) $( \hat{\y}, \boldsymbol{\hat{\gamma}},\hat{\z})$ in the space of these two constraints, for each product $i\in [N]$ and period $t\in [T]$, we check whether a constraint in~\eqref{eq:pers-extension} and~\eqref{eq:conv-env:grad} is violated. For constraint~\eqref{eq:pers-extension}, we obtain a permutation $\hat{\sigma}$ by sorting $\{g_{im}^t\}_{m\in [M]}$ in an decreasing order, that is, $g_{i\hat{\sigma}(1)}^t\ge  \dots \ge  g_{i\hat{\sigma}(M)}^t$, where $g_{im}^t = \hat{z}_{im}^t$ if $m\notin I_i$ and $g_{im}^t = \hat{\gamma}_i^t - \hat{z}_{im}^t$ if $m\in I_i$. The constraint in~\eqref{eq:pers-extension} corresponding to the permutation $\hat{\sigma}$ is 
\begin{equation}\label{eq:violate:conc}
  y^t_i \leq \alpha_i(\h^{\hat{\sigma}}_{i0} ) ( \gamma^t_i - \tilde{z}_{i{\hat{\sigma}}(1)}^t) + \sum_{k \in [M]} \alpha_i(\h^{\hat{\sigma}}_{ik} ) (   \tilde{z}_{i{\hat{\sigma}}(k)}^t - \tilde{z}_{i{\hat{\sigma}}(k+1)}^t).
\end{equation} 
 The given point $( \hat{\y}, \boldsymbol{\hat{\gamma}},\hat{\z})$ satisfies~\eqref{eq:pers-extension} if and only if it satisfies the constraint~\eqref{eq:violate:conc}. For constraint~\eqref{eq:conv-env:grad}, we directly check the subgradient constraint obtained at historical assortments $(\hat{w}^1, \dots, \hat{w}^M) = (\hat{z}_{i1}^{t}/ \hat{\gamma}_i^t,\dots, \hat{z}_{iM}^{t}/ \hat{\gamma}_i^t)$:
\begin{equation}\label{eq:violate:conv}
    y_i^t \ge \gamma_i^t \alpha_i(\hat{w}^{1},\dots, \hat{w}^{M}) + \summ (z_{im}^{t} - \gamma_i^t \hat{w}^{m}) \beta_i^m   \alpha_i(\hat{w}^{1},\dots, \hat{w}^{M}).
\end{equation}
The given point $( \hat{\y}, \boldsymbol{\hat{\gamma}},\hat{\z})$ satisfies~\eqref{eq:conv-env:grad} if and only if it satisfies the constraint~\eqref{eq:violate:conv}. }

\re{The last element in our implementation is to select a starting formulation before calling \texttt{Gurobi}. We solve a collection of relaxed models to obtain continuous solutions and call the separation oracle to add new cuts. The relaxed models contain a subset of constraints~\eqref{eq:pers-extension} and~\eqref{eq:conv-env:grad}, defined as follows:}
\re{\begin{equation}\label{eq:mp:mip}
\begin{aligned}
\max_{\x,\y,\boldsymbol{\rho},\boldsymbol{\gamma},\z} \quad & \frac{1}{T}  \sumt  \sumi r_iy_{i}^t  \\
    \st \quad & \x^t \in [0,1]^{N } \cap \mathcal{X} \text{ and }   \x \in \mathcal{P}  && \for \, t \in [T]  \\
    &\rho^t + \sumi y_{i}^t = 1  \quad && \for\, t \in [T]  \\ &\eqref{eq:relax1},~\eqref{eq:relax2}  \quad && \for \, t \in [T] \text{ and } i \in [N]\\ 
      & \eqref{eq:pers-extension}  && \for \sigma \in \Omega_{it}^k, t \in [T] \text{ and } i \in [N] \\
    & \eqref{eq:conv-env:grad} && \for \tilde{\boldsymbol{w}} \in W^k_{it} , t \in [T] \text{ and } i \in [N],
\end{aligned}     \tag{\textsc{Rlx$_k$}}
\end{equation}}\re{where $\Omega_{it}^k\subset \Omega$ (resp. $W^k_{it}\subset [0,1]^M$) is a subset of permutations (resp. nodes). Given a continuous solution, we call the separation oracle (Algorithm~\ref{alg:sepa:lovasz}) to add new cuts and update $\textsc{Rlx}_k$ to $\textsc{Rlx}_{(k+1)}$. We add $K$ rounds of cuts and set the final $\textsc{Rlx}_K$ with binary constraints as a starting model for \texttt{Callback}. We present details of the implementation process in Algorithm~\ref{alg:cut-plane:lovasz}.}

\begin{algorithm}[h]
    \caption{\re{Separation oracle for~\eqref{eq:pers-extension} and~\eqref{eq:conv-env:grad}}}\label{alg:sepa:lovasz}
    \begin{algorithmic}[1]
        \STATE Input: $(\hat{\x},\hat{\y},\hat{\boldsymbol{\rho}},\hat{\boldsymbol{\gamma}},\hat{\z}) $
        \STATE Cuts$\gets \emptyset$
        \FOR{ $i \in [ N]$ and $ t\in [T]$}
        \STATE Generate $\hat{\sigma}$ by sorting $\{g_{i\hat{\sigma}(m)}^t\}_{m\in [M]}$ such that $g_{i\hat{\sigma}(1)}^t\ge  \dots \ge  g_{i\hat{\sigma}(M)}^t$, where $g_{im}^t = \hat{z}_{im}^t$ if $m\notin I_i$ and $g_{im}^t = \hat{\gamma}_i^t - \hat{z}_{im}^t$ if $m\in I_i$.
        \IF{ \eqref{eq:violate:conc} is violated}
        \STATE Cuts$\gets$~\eqref{eq:violate:conc}
        \ENDIF
        \STATE Find the historical assortments $(\hat{w}^1, \dots, \hat{w}^M) = (\hat{z}_{i1}^{t}/ \hat{\gamma}_i^t,\dots, \hat{z}_{iM}^{t}/ \hat{\gamma}_i^t)$
        \IF{ \eqref{eq:violate:conv} is violated}
        \STATE  Cuts$\gets$~\eqref{eq:violate:conv}
        \ENDIF 
        \ENDFOR 
        \STATE Output: Cuts 
    \end{algorithmic}
\end{algorithm}

\begin{algorithm}[h]
    \caption{\re{A cutting-plane algorithm for \eqref{eq:Conic} with large memory length}}\label{alg:cut-plane:lovasz}
    \begin{algorithmic}[1]
        \STATE Input: Parameters $\r,\{\beta_i^0\}_{i\in[N]}, \boldb,M,N,T$, $\{\Omega_{it}^1,W_{it}^1\}_{i\in [N],t \in [T]}$, and $K$
        \STATE $k\gets 1$   
        \WHILE{ $k \le K$ } 
        \STATE $(\hat{\x},\hat{\y},\hat{\boldsymbol{\rho}},\hat{\boldsymbol{\gamma}},\hat{\z})$ Solve \eqref{eq:mp:mip} 
        \STATE{Cuts $\gets $ Call the separation oracle, Algorithm \ref{alg:sepa:lovasz}, to separate $(\hat{\x},\hat{\y},\hat{\boldsymbol{\rho}},\hat{\boldsymbol{\gamma}},\hat{\z})$ }
        \STATE $\textsc{Rlx}_{k+1} \gets $ Add Cuts into $\textsc{Rlx}_{k}$
        \STATE $k \gets k+1$    
        \ENDWHILE
        \STATE Add binary constraints to the final $\textsc{Rlx}_K$ and solve it through the \texttt{Callback} routine of \texttt{Gurobi} to call the separation oracle Algorithm \ref{alg:sepa:lovasz} 
        \STATE Output: An optimal integer solution $(\x^*,\y^*,\boldsymbol{\rho}^*,\boldsymbol{\gamma}^*,\z^*)$
    \end{algorithmic}
\end{algorithm} 

\newcommand{\norm}[1]{\left\lVert#1\right\rVert}

\section{\re{Supplementary Numerical Studies}}\label{additional:numerical}
\subsection{\re{Performance under Cross-product and Cross-period Constraints}}\label{ec:sec:constraint}
\re{In this numerical study, we add cross-product and cross-period constraints to understand how such constraints impact the computational performance of our formulations. We consider the cardinality constraint for each period:
\begin{equation}\label{eq:constr:card}
    \sum_{i\in [N]} x_i^t \le C \quad \for t \in [T].
\end{equation}
That is, the number of offered products in each period is no more than $C$, a positive integer. For the cross-period constraint, we consider that the number of times provided for a product $i\in [N]$ shall not exceed $K$. That is, 
\begin{equation}\label{eq:constr:offer}
    \sum_{t\in [T]} x_i^t \le K \quad \for i \in [N].
\end{equation}}

\re{We focus on the satiation effects and generate parameters as follows. The revenue and base utility of a product is randomly generated from uniform distributions with ranges of $[1,10]$ and $[-1,1]$, respectively. The history-dependent effect $\beta_i^1$ is uniformly sampled from $[-1,0]$ for each $i\in [N]$. We fix $N=30,T=5$ and $M=1$. We generate 5 synthetic instances for each configuration of $C \in \{5,10,15,30\}$ and $K\in \{2,3,4,5\}$. We use the MILP model, $\mathsf{Env}$, to solve each instance within 3600 seconds and with an optimality gap $0.5\%$.}

\begin{table}[htb]
\centering
  \caption{\re{Computation time, end gap, and root gap of $\mathsf{Env}$ formulation under weak satiation effects with cardinality constraint~\eqref{eq:constr:card} and offering constraint~\eqref{eq:constr:offer} ($N=30,T=5,M=1$)  }}\label{tab:constraint}
\begin{tabular}{ccccccccccccc}
\cline{1-6} \cline{8-13}
 \textbf{C} & \textbf{K} &  $\# \mathsf{sol}$ & $T_{\mathsf{opt} {(s)}}$ & $G_{\mathsf{end} {(\%)}} $  & $G_{\mathsf{root} {(\%)}} $ & &  \textbf{C} & \textbf{K} &  $\# \mathsf{sol}$ & $T_{\mathsf{opt} {(s)}}$ & $G_{\mathsf{end} {(\%)}} $  & $G_{\mathsf{root} {(\%)}} $ \\
\cline{1-6} \cline{8-13}
5  & 2 & 1   & 801.20   & 2.15    & 13.58  & &  15 & 2 & 1   & 796.22   & 1.93    & 13.23   \\
5  & 3 & 4   & 509.65   & 1.73    & 7.49    & &  15 & 3 & 4   & 385.27   & 1.46    & 7.03  \\
5  & 4 & 4   & 369.20   & 1.28    & 5.81    & &  15 & 4 & 4   & 207.41   & 0.86    & 4.86  \\
5  & 5 & 4   & 139.18   & 1.18    & 4.41   & &  15 & 5 & 5   & 86.74    & 0       & 3.26     \\
10 & 2 & 1   & 960.74   & 1.96    & 13.22  & &   30 & 2 & 1   & 871.66   & 1.86    & 13.22   \\
10 & 3 & 4   & 339.98   & 1.42    & 7.03   & &  30 & 3 & 4   & 439.71   & 1.44    & 7.03   \\
10 & 4 & 4   & 187.44   & 0.71    & 4.86  & &  30 & 4 & 4   & 214.48   & 0.93    & 4.86    \\
10 & 5 & 5   & 82.12    & 0       & 3.27  & &  30 & 5 & 5   & 89.89    & 0       & 3.26    \\
\cline{1-6} \cline{8-13}
\end{tabular}
\end{table}

\re{Table~\ref{tab:constraint} summarizes the computation performance of $\mathsf{Env}$ under constraints of~\eqref{eq:constr:card} and~\eqref{eq:constr:offer}. The first two columns are the limitation of cardinality size and offering times, respectively. The third column $\#\textsf{sol}$ indicates the number of instances solved to optimality within 3600 seconds. $T_{\textsf{opt}}$ is the average computation time of solved instances and $G_{\textsf{end}}$ is the average end gap of unsolved instances. The last column $G_{\textsf{root}}$ is the average root gap of all instances. Table~\ref{tab:constraint} shows that our formulation is more sensitive to the cross-period offering constraint~\eqref{eq:constr:offer} than the cross-product cardinality constraint~\eqref{eq:constr:card}. For instance, under the same offering times $K$, the average computation time, end gap, and root gap are stable under different cardinality sizes $C$. However, under the same $C$, if we increase offering times $K$ from $2$ to $4$, the average computation time, end gap, and root gap shrink by more than one-half. It is because our formulation does not characterize the problem structure in the space of cross-period choice probabilities, $(y_i^t,y_i^{t+1})$. Although not all instances are solved to optimality, the end gap is no more than $2.5\%$. It indicates the robustness and efficiency of our formulation. }

\subsection{\re{Performance of Instances with Large Memory Length}}\label{ec:sec:numerical:largem}
\re{This numerical study demonstrates that our \eqref{eq:Conic} formulation combined with the cutting-plane algorithm can solve instances with moderately large memory length.}

\re{We fix the number of products $N=20$ and planning horizon $T=M+1$ for $M\in \{4,5,6\}$. We consider both addiction and satiation effects. Specifically, a product has an addiction (resp. satiation) effect with probability $\theta$ (resp. $1-\theta$), where $\theta\in \{0,0.1,0.2\}$. For products with addiction (resp. satiation) effects, its $\{\beta_i^m\}_{m\in [M]}$ are uniformly sampled from $[0,1]$ (resp. $[-1,0]$). For each configuration of $M,T$, and $\theta$, we randomly generate 10 instances. The revenue and the base utility of a product are from uniform distributions with ranges of $[1,10]$ and $[-1,1]$, respectively. To initialize the cutting-plane algorithm, we set $\Omega_{it}^1$ by randomly generating 2,10, and 20 different permutations over $[M]$ for memory lengths of 4,5, and 6, respectively. We initialize $W_{it}^1$ by adding $\boldsymbol{w}\in \{0,1\}^M$ satisfying $\norm{\boldsymbol{w}}_2^2 \in \{0,1,2,M\}$. When computing the starting model for \texttt{Callback}, we do not fix $K$ rounds to add cuts. Instead, we call the separation oracle until the reduction of the objective value of~\eqref{eq:mp:mip} is less than a small constant $\epsilon$. We set $\epsilon = 10^{-8}$ in the numerical study.}

\begin{table}[htb]
\centering
  \caption{\re{Average computation time of solved instanced and average end gap of unsolved instances}}\label{tab:double}
\begin{tabular}{cccccc}
\hline 
\textbf{M}        & \textbf{T} & $\theta$  &  $\# \mathsf{sol}$  & Time $(s)$  & $G_{\mathsf{end} {(\%)}} $    \\ \hline 
4         & 5 & 0   & 10   & 317.13& 0       \\
4         & 5 & 0.1   & 10 & 2.94  & 0         \\
4         & 5 & 0.2   & 10 & 17.00 & 0         \\ 
\hline 
5         & 6 & 0     & 9  & 273.89 & 1.6     \\
5         & 6 & 0.1   & 9   & 1632.32 & 3.19   \\
5         & 6 & 0.2    & 10  & 170.72 & 0        \\ 
\hline 
6         & 7 & 0   & 6   & 3416.57 & 2.96    \\
6         & 7 & 0.1  & 8  & 963.53   & 3.8    \\
6         & 7 & 0.2   & 7  & 18.16  & 2.75      \\
\hline 
\end{tabular}
\end{table}

\re{Table~\ref{tab:double} summarizes the computation results, where we set a time limit of 7200 seconds and an optimality gap tolerance of 0.5\%. The first three columns show the memory length, planning horizon, and the proportion of products with positive effects. The last three columns present the number of solved instances, the average time of solved instances (including set-up time of the starting model and \texttt{Callback} time), and the average end gap of unsolved instances. Table~\ref{tab:double} shows that the cutting-plane algorithm can solve most instances to optimality within two hours, even when $M=6$. For instance, when $M=6$ and $\theta=0.2$, 7 instances are solved to optimality using only 18.16 seconds on average. While some instances cannot achieve optimality within two hours, their end gaps are small and no more than 4\%.}

% Cross-product constraint:
% \begin{itemize}
%     \item \textit{Maximum number of offerings}: $\sum_{i\in [N]} x_i^t \le C_t$ for $t\in [T]$
%     \item \textit{Minimum number of offerings within a cluster}: $\sum_{i \in S_j} x_i^t \ge C_j$ for $j\in [J], t\in [T]$
%     \item \textit{Inter-no-touch constraints}: $x_i^t + x_j^t\le 1$ 
%     \item \textit{Base-product constraints}: $ x_i^t \ge 1$ for $i\in B_i, t\in  [T]$
% \end{itemize}

% Cross-period constraint:
% \begin{itemize}
%     \item \textit{Non-overlapping constraints}: $\sum_{k \in [K]\cup \{0\}} x_i^{t+k}\le 1$ for $i\in [N], t\in [T]$
%     \item \textit{Maximum number of offerings per product}: $\sumt x_i^t \le T_i$ for $i\in [N]$
% \end{itemize}

% References(outcomment the appropriate case) 

% CASE 1: BiBTeX used to constantly update the references 
%   (while the paper is being written).
% if more than one, comma separated

% CASE 2: BiBTeX used to generate mypaper.bbl (to be further fine tuned)
%\input{mypaper.bbl} % outcomment this line in Case 2

\end{document}